\documentclass[12pt]{amsart}

\usepackage[usenames,dvipsnames,svgnames,table]{xcolor}
\usepackage{graphicx}
\usepackage{inputenc}
\usepackage{amsmath, latexsym, amsfonts, amssymb, amsthm, amscd, stmaryrd, MnSymbol}
\usepackage{hyperref}

\usepackage{caption}
\usepackage{subcaption}
\usepackage{float}

\usepackage[left=3cm, right=2.5cm, top=2.5cm, bottom=2.5cm]{geometry}

\renewcommand{\emptyset}{\font\cmsy = cmsy10 at 10pt
 \hbox{\cmsy \char 59}
}

\newtheorem{theorem}{Theorem}[section]
\newtheorem{lemma}[theorem]{Lemma}
\newtheorem{corollary}{Corollary}
\newtheorem{proposition}{Proposition}
\newtheorem{rem}{Remark}

\newcommand{\p}{\mathbb{P}}

\newcommand{\e}{\mathbb{E}}

\def\indi{\mathbf{1}}

\def\N{\mathbb{N}}

\def\t{\textbf{t}}

\def\indi{\mathbf{1}}
\newcommand{\ud}[1]{\, \mathrm{d}#1}

\newcommand{\Exp}[1]{\mathbb{E}\left[#1\right]}

\begin{document}

\title[The coalescent structure of GW trees in varying environments]{ The coalescent structure of Galton-Watson trees in varying environments}
 
 \author{Simon C. Harris}
\address{Department of Statistics, University of Auckland, Private Bag 92019, Auckland 1142, New Zealand}
\email{simon.harris@auckland.ac.nz}

\author{Sandra Palau}
\address{IIMAS, Universidad Nacional Aut\'onoma de M\'exico. CDMX, 04510, Ciudad de M\'exico, M\'exico}
\email{sandra@sigma.iimas.unam.mx}
 
 \author{Juan Carlos Pardo}
\address{CIMAT A.C. Calle Jalisco s/n. C.P. 36240, Guanajuato, Mexico}
\email{jcpardo@cimat.mx}
\thanks{Part of this work was undertaken during an academic visit of SP and JCP supported by the University of Auckland. SH acknowledges University of Auckland FRDF (3720685) and SP acknowledges support by UNAM-PAPIIT  (IN103924). }

\begin{abstract}

We investigate the genealogy of a sample of $k\geq2$ particles chosen uniformly without replacement from a population alive at large times  in a critical discrete-time Galton-Watson process in a varying environment (GWVE). 
We will show that subject to an explicit deterministic time-change involving only the mean and variances of the varying offspring distributions, the sample genealogy always converges to the same universal genealogical structure; it  has the same tree topology as Kingman's coalescent, and the coalescent times of the $k-1$ pairwise mergers look like a  mixture of independent identically distributed times.
Our approach uses $k$ distinguished \emph{spine} particles and a suitable change of measure  under which (a) the spines form a uniform sample without replacement, as required, but additionally (b) there is $k$-size biasing and discounting according to the population size.  
Our work significantly extends the spine techniques developed in Harris, Johnston, and Roberts \emph{[Annals Applied Probability, 2020]} for genealogies of uniform samples of size $k$  in near-critical continuous-time Galton-Watson processes, as well as a two-spine GWVE construction  in Cardona and Palau \emph{[Bernoulli, 2021]}.
Our results complement recent works by Kersting \emph{[Proc. Steklov Inst. Maths., 2022]} and Boenkost, Foutel-Rodier, and Schertzer \emph{[arXiv:2207.11612]}.
\end{abstract}

\maketitle
\noindent {\sc Key words and phrases}:   Galton-Watson processes in varying environments, Galton-Watson trees in varying environments, coalescent, genealogy, spines

\noindent MSC 2020 subject classifications: 60J80, 60F17, 60G09.

\section{Introduction and main results}

In this manuscript, we are interested in the genealogy of a sample of $k\geq2$ particles chosen uniformly without replacement from a population alive at a large time $n$ in a critical discrete-time Galton-Watson process in a varying environment (GWVE for short). 

GWVE  generalises the classical Galton-Watson processes since the offspring distributions may  vary from generation to generation. This class of branching processes has attracted a lot of attention recently, see for instance \cite{BS2015, BP2017, CardonaPalau, FLL2022, gonzalez2019, kersting2020unifying, kersting2021genealogical, kersting2017discrete}.
Initial research on GWVE was hampered by some special features making difficult to generalise the classical long-term behaviour  of Galton-Watson processes concerning extinction or divergence to infinity and even the classification on supercritical, critical and subcritical regimes.
For a resolution to such issues, see Kersting \cite{kersting2020unifying} or the monograph of Kersting and Vatutin \cite{kersting2017discrete} for further details.

Our study  is motivated by the recent result of Kersting \cite{kersting2021genealogical}  and the two spine construction of Cardona and Palau \cite{CardonaPalau}. The former article deals not only with the description of the generation of the most recent common ancestor of all particles alive conditional on survival (asymptotically) but also with determining the limiting random object of the reduced process of a critical GWVE which turns out to be a time change Yule process. The latter article provides a probabilistic approach to Yaglom's limit for critical GWVE using a two-spine argument (see also Kersting \cite{kersting2020unifying} for an analytic technique). Surprisingly, up to a deterministic time change depending only on the offspring means and variances of the environments,  the limiting objects for the GWVE are the same as for the critical Galton-Watson (GW for short) with constant environment, and our results confirm that this is also the case when we sample  $k$  particles  uniformly without replacement from all particles alive in a GWVE at  large times. 

Let $(\Omega, \mathcal{F}, \mathbb{P})$ be a probability space. A {\it varying environment} is a sequence  $e=(q_1,q_2,\ldots)$  of probability measures on  $\N_0=\{0,1,2,\ldots\}$. Thus, a  {\it Galton-Watson process}  
$Z = (Z_n : n\geq 0)$ 
\emph{in a varying environment}  $e$  is  defined  as 
$$ Z_0=1\quad \mbox{ and } \quad 
Z_n:=\underset{i=1}{\overset{ Z_{n-1}}{\sum}}\chi_{i}^{(n)}, 
\qquad n\geq 1,$$
where   
$(\chi_{i}^{(n)}: i,n\geq 1)$
is a sequence of independent random variables
such that
$$\p(\chi_{i}^{(n)}=k)=q_n(k), \qquad 
\quad k\in \N_{0},\ i,\ n\geq 1.$$ 
In other words, $\chi_{i}^{(n)}$  denotes the number of offspring of the  $i$-th particle in the  $(n-1)$-th generation. We denote by $(Z; \p^{(e)})$ the law of such a process.

Let $f_n$ be the generating function associated with $q_n$, that is 
$$ f_n(s):= \Exp{s^{\chi_{i}^{(n)}}}
=\sum_{k=0}^{\infty} s^k q_n(k), \qquad\textrm{for}
\quad 0 \leq s\leq 1\quad\textrm{and}\quad  n\geq 1.$$ 
By the branching property, it is straightforward to deduce that the generating function of  $Z_n$  can be written  in terms of  $(f_1, f_2,\ldots)$  as 
\begin{equation*}
\e^{(e)}\left[s^{Z_n}\right]=f_1\circ\cdots\circ f_n(s), 
\qquad\textrm{for}\quad 0 \leq s\leq 1\quad\textrm{and}\quad  n\geq 1,
\end{equation*}
where  $f\circ g$  denotes the composition of  $f$  with  $g$.  By differentiating  the previous expression with respect to  $s$ and evaluating at $s=1$,   we obtain
\begin{equation}\label{eq_media_Zn}
\e^{(e)}[Z_n]=\mu_n,\qquad  \mbox{and} \qquad 
\frac{\e^{(e)}[Z_n(Z_n-1)]}{\e^{(e)}[Z_n]^2}
= \sum_{k=0}^{n-1}\frac{\nu_{k+1}}{\mu_{k}}, 
\qquad n\geq 1,
\end{equation}
where  $\mu_0:=1$\  and, for any  $n\geq 1$,
\begin{equation*}\label{defi mu}
\mu_n:= f_1'(1)\cdots f_n'(1), \qquad \mbox{ and }\qquad 
\nu_n:=\frac{f_n''(1)}{f_n'(1)^2}
=\frac{\text{Var}\left[\chi_{i}^{(n)}\right]}{\Exp{\chi_{i}^{(n)}}^2}
+\left(1-\frac{1}{\Exp{\chi_{i}^{(n)}}}\right),
\end{equation*}
where  Var$\left[\cdot\right]$  denotes the variance under $\p$. For more  details, we refer the reader to  \cite{kersting2017discrete}. 

As we mentioned before, we require some knowledge on the long-term behaviour of a specific class of GWVE. GWVEs may behave quite  differently to standard GW processes, such as possessing multiple rates of growth, see e.g. MacPhee and Schuh \cite{macphee1983galton}, but we will not consider  such cases here.  As explained by Kersting \cite{kersting2020unifying}, these strange behaviours can be excluded by the following condition: for every $\epsilon>0$, there is a finite constant $c_\epsilon$ such that for all  $n\ge 1$
\begin{equation}\label{eq_cond_kersting}\tag{$\star$}
\mathbb{E}\left[\left(\chi_{1}^{(n)}\right)^2\mathbf{1}_{\{\chi_{1}^{(n)}> c_\epsilon(1+\mathbb{E}[\chi_{1}^{(n)}])\}}\right]\leq \epsilon \mathbb{E}\left[\left(\chi_{1}^{(n)}\right)^2\mathbf{1}_{\{\chi_{1}^{(n)}>2\}}\right].
\end{equation} 
According to Kersting \cite{kersting2020unifying}, a GWVE is known as \emph{regular} if it satisfies Condition  \eqref{eq_cond_kersting}. In what follows, we always consider regular GWVE.

However, directly verifying Condition \eqref{eq_cond_kersting} can be difficult so,   as  suggested  by Kersting \cite{kersting2020unifying},  one may alternatively use the following mild third moment condition: there exists $ c>0$ such that    
\begin{equation}\label{eq_cond_mild}
f_n'''(1)\leq c f_n''(1)(1+ f_n'(1)),  \mbox{ for any } n\geq 1,
\end{equation} 
which, according to Kersting \cite[Proposition 2]{kersting2020unifying},  implies condition \eqref{eq_cond_kersting}. As it is  explained in \cite{kersting2020unifying},  Condition \eqref{eq_cond_mild} is easier to verify and many common distributions satisfy it, e.g.   binomial, geometric, hypergeometric, Poisson,  and  negative binomial distributions and also  
random variables that are almost surely uniformly bounded  by a constant.

It turns out that under Condition \eqref{eq_cond_kersting}, the behaviour of a GWVE   is basically determined by the following two sequences  $(\mu_n, n\geq 0)$ and $(a_n^{(e)}, n\geq 0)$;  where,  given a varying environment  $e$, 
\begin{equation*}\label{def:an}
a_0^{(e)}=1, \qquad \mbox{and}\qquad 
a_n^{(e)}:= \frac{\mu_n}{2}\sum_{k=0}^{n-1}\frac{\nu_{k+1}}{\mu_{k}},
\qquad n\geq 1.
\end{equation*}
With both sequences in hand, regular GWVEs can be classified into four distinct classes:  {\it subcritical, critical, supercritical and asymptotically degenerate}. In the latter,  the process may freeze in a positive state. For complete details  about such classification, we refer to \cite{kersting2020unifying}. According to \cite[Theorems 1 and 4]{kersting2020unifying},  a regular GWVE is \emph{critical} if and only if
\begin{equation}
\label{def: critical}
\underset{n\rightarrow\infty}{\lim}a_n^{(e)}=\infty \qquad \mbox{and}\qquad \underset{n\rightarrow\infty}{\lim}\frac{a_n^{(e)}}{\mu_n}=\infty.
\end{equation}  
In this case, the process becomes extinct  a.s.  and 
\begin{equation}\label{limit a}
\lim\limits_{n\to \infty}\frac{a_n^{(e)}}{\mu_n} \p^{(e)}(Z_n>0)=1. 
\end{equation}
As it  was noted by Kersting \cite{kersting2020unifying} and Cardona and Palau \cite{CardonaPalau}, for a  critical GWVE  the so-called {\it Yaglom's} limit exists, that is 
\begin{equation}\label{eq: Yaglom}
\left(\frac{Z_n}{a_n^{(e)}}; \p^{(e)}(\ \cdot\ | Z_n > 0)\right) \stackrel{(d)}{\longrightarrow}  \left(\mathbf{e}; \p\right) , \qquad \mbox{ as }\ n\rightarrow \infty,
\end{equation}
where  $\mathbf{e}$ or  $\left(\mathbf{e}; \mathbb{P}\right)$  denotes a standard exponential random variable, under $\p$.  It is important to note that the previous  limit was obtained by Kersting \cite{kersting2020unifying} under condition \eqref{eq_cond_kersting} using analytical arguments, and by Cardona and Palau  \cite{CardonaPalau}  under the third moment condition in \eqref{eq_cond_mild} using a probabilistic approach based on a two-spine decomposition argument, where a {\it spine} is a distinguished (marked) genealogical line.  

The main goal of this article is to show the emergence of an explicit universal limiting genealogy when sampling $k$ particles uniformly at random without replacement at large times $n$ in a critical regular GWVE conditioned to survive, 
as described below. 

On the event $\{Z_n\geq k\}$, pick $k$ particles $U_n^{(1)},\dots, U_n^{(k)}$ uniformly random without replacement from the $Z_n$ particles alive at time $n$. Let $\mathcal{P}_m(n)$ be the partition of $\{1,\dots, k\}$ induced by letting $i$ and $j$ be in the same block if the particles $U_n^{(i)}$ and $U_n^{(j)}$ share a common ancestor at time $m$. 
Let $v_m^k(n)$ be the number of blocks in $\mathcal{P}_m(n)$, that is the number of distinct ancestors of $U_n^{(1)},\dots, U_n^{(k)}$ at time $m$. Denote the last time when there are at most  $i$ blocks as follows
$$B_i^k(n):=\max\{m\geq 0: v_m^k(n)\le i \}, \qquad i=1,\dots, k-1.$$
It will also be convenient to define the corresponding unordered times $(\widetilde{B}_1^k(n),\dots,\widetilde{B}_{k-1}^k(n))$ as a uniformly random permutation of $(B_1^k(n),\dots,B_{k-1}^k(n))$.

 To describe the tree topology of the sample, we let the partition 
$\tt{P}_i(n):=\mathcal{P}_{B_i^k(n)+1}(n)$  for $i=1,2,\dots,k-1$ and $\tt{P}_0(n):=\{1,\dots, k\}$. The former definition follows from the fact that we only observe changes in the partitions one generation after the $B_i^k(n)$, since $B_i^k(n)$ is the last time where there are at most $i$ blocks. Then $\mathcal{H}:= \sigma(\tt{P}_0(n),\dots, \tt{P}_{k-1}(n))$ contains all topological information about the genealogical tree of the sample of $k$ particles. However, we note that $\mathcal{H}$ includes no direct information about the times of the splits, only the shape of the tree.

We  let 
\begin{equation}\label{rhonDef}
\rho_0=0 \qquad \mbox{and} \qquad \rho_n=\frac{\e^{(e)}[Z_n(Z_n-1)]}{\e^{(e)}[Z_n]^2}= \sum_{k=0}^{n-1}\frac{\nu_{k+1}}{\mu_{k}},
\qquad n\geq 1.
\end{equation}
For $n\geq 1$, the family $(\rho_k/\rho_n,{0\leq k\leq n})$ describes an increasing sequence that can be thought as a cumulative probability distribution. We can thus define its right-continuous generalised inverse as
\begin{equation}\label{SntDef}
\tau_n(t)=\max\left\{k\in\{0,\dots,n\}: \frac{\rho_k}{\rho_n}\leq t\right\}, \qquad t\in [0,1].
\end{equation}
This deterministic time change is due to Kersting \cite{kersting2021genealogical} and by means of this function, the distances  between generations in the reduced process (that is, the process generated by particles which are ancestors  of those alive at a given time) are re-scaled. As we now briefly explain, this time change was already hidden in \cite{CardonaPalau} where a two-spine decomposition was used to prove the Yaglom limit \eqref{eq: Yaglom},  although the authors never defined it explicitly.   Denote by $\ddot{Z}_m$ the population size at the $m$-th generation of a tree with two distinguish particles (spines) selected at time $n$, with $m\leq n$. Let $\psi$ be the generation of the most recent common ancestor of the two distinguished particles.
Then, $\ddot{Z}_n$ is made up of  the descendants attached to the longer spine (genealogical line of one particle) plus the descendants attached to the shorter spine (genealogical line of the other  particle up to the time that they coalesce), as well as the two spines themselves. 
In \cite{CardonaPalau}, the term $1-\rho_{\psi+1}/\rho_n$ was exactly the factor to provide the correct (Yaglom) normalisation for the descendants attached to the shorter spine. More precisely, from  \cite[Proposition 3]{CardonaPalau}, we see that for any  $x\in [0,1]$,
$$\underset{n\rightarrow \infty}{\lim}\mathbb{P}\left(\psi+1\leq \tau_n(x)
 \right)= \underset{n\rightarrow \infty}{\lim}\mathbb{P}\left( \frac{\rho_{\psi+1}}{\rho_{n}}\leq x
 \right)=x.
$$
In other words, if we  take a sample of two particles at generation $n$, after rescaling properly,  in the limit the time of their most recent common ancestor is intimately related to the time change  $\tau_n$. The main result of this paper generalise the previous statement  to $k\ge 2$ selected particles.
\begin{theorem}\label{principal}\emph{\textbf{Uniform sampling from a critical regular GWVE at large times.}} 
Let us assume that conditions \eqref{eq_cond_kersting} and \eqref{def: critical} are satisfied. Consider $\{t_1,\dots,t_{k-1}\}\subset (0,1)$ with $t_i\neq t_j$ for any $i\neq j$. Then, 
\begin{equation}\label{eqtheorem}
\begin{split}
\underset{n\rightarrow \infty}{\lim}
\p^{(e)}
&\left(\widetilde{B}_1^k(n)\geq \tau_n(t_1), \dots,\widetilde{B}_{k-1}^k(n)\geq \tau_n(t_{k-1})\mid Z_n\geq k\right)\\
&\hspace{2cm}=\int_0^{\infty}\frac{k}{(1+\lambda)^{2}}\prod_{i=1}^{k-1}\left(1-\frac{1}{1+\lambda (1-t_i)}\right)\ud \lambda\\
&\hspace{2cm}=k\left(\prod_{i=1}^{k-1} \frac{t_i-1}{t_i}-\underset{j=1}{\overset{k-1}{\sum}}\frac{1-t_j}{t_j^2}\prod_{i=1, i\neq j}^{k-1}\frac{1-t_i}{t_j-t_i}\log\left(1-t_j\right)\right).
\end{split}
\end{equation}
Furthermore, the  times $(\widetilde{B}_1^k(n),\dots,\widetilde{B}_{k-1}^k(n))$ are asymptotically independent of the sample tree topology $\mathcal{H}$, and the partition process $(\tt{P}_0(n), \tt{P}_1(n),\dots, \tt{P}_{k-1}(n))$ that describes the tree topology has the following description:
\begin{enumerate}
\item if  \ $\tt{P}_i(n)$ contains blocks of sizes $b_1,..., b_{i+1}$, the probability that the next block to split is block $j$ converges to $(b_j-1)/(k-i-1)$ as $n\rightarrow\infty$;
\item if a block of size $b$ splits, it creates two blocks where the size of first block  is  uniformly chosen between $\{1,\dots, b-1\}$.
\end{enumerate}
\end{theorem}

This new result for the discrete-time critical GWVE is 
a generalisation of
the continuous-time critical GW in a constant environment as first discovered in Harris, Johnston and Roberts \cite{harris2020coalescent}. 
Near the completion of writing up this article, we became aware that 
some closely related results were in progress in Boenkost, Foutel-Rodier, and Schertzer \cite{BFS2022}. We provide a discussion about the differences between our results and approaches at the end of this section.
Whilst we will follow a similar general approach (involving $k$-spines) to \cite{harris2020coalescent}, we introduce a number of innovations and must deal with additional challenges and technicalities that discrete time and a varying environment introduce. 
Of course, the same limiting universal sample genealogy emerges in both situations, just as intuitively expected from Kersting \cite{kersting2021genealogical} where the reduced tree associated to all particles alive at generation $n$ and conditional on survival converges in the Skorokhod topology towards a time changed Yule process, as $n$ goes to infinity. We observe that the time change appearing in Theorem \ref{principal} is the same as it appears in \cite{kersting2021genealogical}. 
The topology of the limiting sample tree here (that we describe going forward in time) agrees with the Kingman $k$-coalescent (usually described going backwards in time, with any pair of block equally likely to merge next). 
As in the middle line in equation \eqref{eqtheorem}, the split times (or, if thinking backwards in time, the coalescent times) can be represented as a mixture of $k-1$ independent identically distributed random variables, similarly as described in Harris et. al. \cite{harris2020coalescent}. For further details, see the construction for critical genealogies (Theorem 4) in the same reference.

In order to understand the genealogy of a uniform sample of size $k$ taken at a large time in a critical GWVE, and then prove Theorem \ref{principal}, we will first generalise the GWVE two-spine construction in \cite{CardonaPalau} to describe a special random tree with $k$-spines.
We introduce a new measure $\mathbf{Q}_n^{(e,k, \theta)}$ under which (a) the $k$ spines form a uniform sample without replacement at time $n$, as required, but additionally (b) there is $k$-size biasing and discounting at rate $\theta$ by the population size at time $n$, that is the measure $\mathbf{Q}_n^{(e,k, \theta)}$ will be constructed  via a change of measure with respect to a reference probability measure (that we will described below) and whose Radon-Nikodym derivative  is proportional to
\[
Z_n(Z_n-1) \times\cdots \times(Z_n-k+1)e^{-\theta Z_n}.
\]
We interpret the latter as a $k$-size biasing together with  discounting at rate $\theta$. The discounting  is necessary since  $Z_n$ may not have finite $k$-th  moments and  it helps to control large population sizes, although in the limit $\theta$ will not really play any role.


Although $\mathbf{Q}_n^{(e,k, \theta)}$ will be constructed to have these properties, we give a complete forward in time construction for the random tree with $k$ spines under $\mathbf{Q}_n^{(e,k, \theta)}$. 
The $k$ size-biasing with discounting is essential as it provides a great deal of structural independence in the $k$-spine construction which ultimately enables explicit calculations (which would otherwise be intractable due to complex dependencies).  
However, this spine construction with $k$ size-biasing alone requires $k$-th moments for the offspring distributions to even exist (an issue encountered in \cite{harris2020coalescent} which there necessitated an ad hoc - and somewhat inelegant - truncation argument).
Introducing discounting by the final population size alongside the $k$-size biasing, whilst it requires some extra analysis, is quite natural and turns out to form an elegant extension to the theory. Importantly, including discounting means that no additional offspring moment assumptions are required in our approach, with only the second moment conditions for criticality. 
Indeed, analogous $k$-size biasing with discounting for $k$-spines was developed  by Harris, Johnston, and Pardo \cite{HJP2022} as an essential tool to analyse sampling from critical continuous-time GW processes with heavy-tailed offspring. In such situations,  unusually large family sizes lead to a phenomena of multiple mergers in the coalescent, in stark contrast to only pairwise mergers found with finite variance offspring. Also see Abraham and Debs \cite{AD2020} for some related $k$-spine changes of measures for classical GW processes.

Combining the special properties of $\mathbf{Q}_n^{(e,k, \theta)}$ together with the Yaglom theorem for GWVE from Kersting \cite{kersting2020unifying} or Cardona-Palau \cite{CardonaPalau}, plus the time-change from Kersting \cite{kersting2021genealogical}  that encodes variation within the environment, we can undo the $k$-size biasing and discounting to reveal how the universal sample genealogy emerges in the limit.

It is worth noting that our description and construction of $\mathbf{Q}_n^{(e,k, \theta)}$ is exact for all times $n$ and works in generality for all GWVE processes. 
In fact, we only rely on having a \emph{critical} GWVE to get the limiting genealogy by using the Yaglom result and to ensure that in the limit the coalescent times are spread out over the entire time period and do not concentrate near $0$ or $n$.  
The limiting universal genealogy also appears automatically from our construction of $\mathbf{Q}_n^{(e,k, \theta)}$. In particular, we do not need to guess any limiting genealogical structure in advance. To be more precise, we would like to emphasise that the techniques here presented are so robust that can be applied  to  other scenarios even when condition \eqref{eq_cond_kersting} fails and that they can be  extended to investigate a suitably defined ``heavy-tailed" GWVE. Indeed, in the so-called ``linear fractional" case, condition \eqref{eq_cond_kersting} may fail but the Yaglom result still holds (as noted in Example 8 in \cite{kersting2020unifying}), and further we can also verify that the limit of the coalescent times are spread out over the entire time period and do not concentrate at $0$ or $n$. 
In particular, Theorem \ref{principal} can still hold in fractional cases even when condition \eqref{eq_cond_kersting} fails (work in progress \cite{HPP2022}).

We conclude our discussion about our main result with few words about the asymptotic independence between the split times and the sample tree topology.  In fact, the  limiting sample tree topology  will be fully  described under $\mathbf{Q}^{(e,k,\theta_n)}_{n}$. Indeed,  it turns out that the limiting topology under $\mathbf{Q}^{(e,k,\theta_n)}_{n}$ remains unchanged when we undo the additional $k$-size biasing and discounting in order to get back to a uniform sample of size $k$ under $\p^{(e)}(\,\cdot\, | Z_n\geq k)$, as required. To see this intuitively, we first observe that undoing the $k$-size biasing and discounting only involves the final scaled population size at time $n$, this being made up of all the sub-populations coming off the spines (plus the $k$ spines themselves). However, in the large time limit, the subpopulations coming off any spine will always be $1$-size biased.\footnote{Except at spine split times where the births are $2$-size biased. However, we will see that the rate of the $2$ and $1$-size biased are asymptotically similar.}
Importantly, in particular at any given moment the subpopulations coming off each branch carrying spines  are independent from the number of spines actually following that branch. 
This means that the rate at which sub-populations coming off the spines contribute to the final scaled population only depends on the number of spine branches alive, not on the sizes of the spine groupings.
Since we only observe spines splitting into two in the limit, this means that the final scaled population under $\mathbf{Q}^{(e,k,\theta_n)}_{n}$ in the limit \emph{only depends on the splitting times of the spines but not on the precise topology (i.e. numbers of spines along each branch)}. Thus, undoing the $k$-size biasing and discounting by the final population will have no effect on any probabilities relating to the topology. That is, in the large $n$ limit, $\p^{(e)}(\, \cdot\, | Z_n\geq k)$ has exactly the same spine topology as $\mathbf{Q}^{(e,k,\theta_n)}_{n}$ as required.

As mentioned, we would like to say a few more words about the differences of our result and that of the closely related work of Boenkost, Foutel-Rodier, and Schertzer  \cite{BFS2022}. Despite that the limiting genealogy  of  $k$ uniformly sampled individuals is the same in both papers, our approaches are quite different and 
complementary to each other. Indeed,  our approach follows a forward in time construction and exact representations about sampling at any fixed time for a given GWVE, from which  the limiting object  emerges directly using the notion of criticality  from Kersting \cite{kersting2020unifying} and the corresponding Yaglom results. On the other hand, the result in \cite{BFS2022} considers a sequence of  what they introduced as near critical GWVEs,  that is, their re-scaled means and  cumulative variances converge towards a c\`adl\`ag process (in the Skorokhod topology) and a non-decreasing continuous function, respectively (the latter guarantees that the variances between generations do not fluctuate too strongly). A near critical sequence of  GWVEs 
 is  critical (on each fixed time horizon)  in the sense  of \cite{kersting2020unifying},  as it is mentioned in Remark 2.2 in \cite{BFS2022},  but we emphasise  that this notion is for sequences of GWVEs rather than a given GWVE.  
Moreover, their approach relies on a spinal decomposition technique allowing them to get a Yaglom limit result for the rescaled size of a sequence of  nearly critical GWVEs conditional on survival under a Lindeberg type condition. In addition, their spinal approach  allows them to get a convergence 
in the Gromov-Hausdorff-Prohorov  topology 
of the genealogical structures of the populations at a fixed time horizons  (where the sequence of trees is seen  as a sequence of metric spaces) 
towards a  limiting metric space which turns out to be  a time-changed version of the Brownian coalescent point process.  Perhaps the main difference between both  works is not only that we are not considering a sequence of  GWVEs,  but also that our approach is quite robust in the sense that it could be applied not only to other scenarios where a Yaglom limit still holds (as in their near critical cases, cf. \cite{harris2020coalescent}), but  also significantly where the limiting genealogy is not necessarily given by binary splittings (work in progress \cite{HPP2022}, cf. \cite{HJP2022}). On the other hand, their approach in \cite{BFS2022} appears to rely heavily on the fact that the limiting genealogy is associated with  the Brownian coalescent point process, i.e. only binary splittings can be expected in the limit. 
Another main difference is that we introduce discounting to avoid $k$-moment assumptions whereas  they use truncation and the method of moments to identify the limiting object.

We also mention that  Conchon-Kerjan, Kious, and Mailler \cite{conchon} recently studied certain critical Galton Watson trees in random environment conditioned to be large and showed that they converge to the Aldous continuum random tree.

The structure of this paper is as follows. Section 2 is devoted to the study of  rooted trees with spines. Here we also introduce Galton-Watson trees in varying environments with spines as well as the change of measure $\mathbf{Q}_n^{(e,k, \theta)}$, under which several functionals of such trees are determined. In particular, we provide the law of the   last time where all spines are together (implicitly the law of the first spine split time) and the offspring distribution of particles  on and off the spines. With these tools in hand, we provide a forward description of such trees. Section 3 is devoted to the scaling limits of critical Galton-Watson trees with $k$ spines under $\mathbf{Q}_n^{(e,k, \theta)}$. Finally, in Section 4 the main result is proved.

\section{Galton-Watson trees with spines} In this section, we introduce Galton-Watson trees in varying environments with spines and study some of their distributional properties. Moreover, we also introduce a change of measure under which the spines form a uniform sample without replacement at a given fixed time. The main result in this section is a forward construction of the spines under the aforementioned new measure at fixed time. 

\subsection{Rooted trees with spines}
Let us recall  the so-called {\it Ulam-Harris labelling}. Let  
$$\mathcal{U}:=\{\varnothing\}\cup\underset{n=1}{\overset{\infty}\bigcup}\N^n$$  
be the set of finite sequences of positive integers with  $\N=\{1,2,\dots\}$. For  $ u\in \mathcal{U}$, we define the \textit{length}  of  $u$  by  $|u|=n$, if   $u=(u_1,\cdots, u_n)\in \N^n$ and  $|\varnothing|=0$.   If  $u=(u_1,\dots,u_n)$  and  $v=(v_1,\dots,v_m)$  belong to $\mathcal{U}$,\ we denote by  $uv:=(u_1,\dots,u_n,v_1,\dots,v_m)$  the \textit{concatenation} of  $u$  and  $v$,\ with the convention  that  $u\varnothing=\varnothing u=u$.  
We say that $v$ is an ancestor of $u$ and write $v\preccurlyeq u$ if there exists $w\in\mathcal{U}$ such that $u=vw$. For $v\preccurlyeq u$, we define the \textit{genealogical line between $v$ and $u$} as $\llbracket v,u\rrbracket:=\{w\in\mathcal{U}: v\preccurlyeq w\preccurlyeq u\}$.
If $v=\varnothing$, we just say that $\llbracket \varnothing,u\rrbracket$ is the \textit{genealogical line}  of $u$.
A {\it rooted tree}  $\textbf{t}$  is a subset of  $\mathcal{U}$  that satisfies 
\begin{itemize}
\item   $\varnothing\in\textbf{t}$.
\item   $\llbracket\varnothing,u\rrbracket\subset\textbf{t}$ for any $u\in\textbf{t}$.
\item  For every $u\in \textbf{t}$  there exists a number $l_u(\textbf{t})$ such that $uj\in \textbf{t}$ if and only if  $1
\leq j\leq l_u(\textbf{t})$.
\end{itemize}
The integer $l_u(\textbf{t})$ represents the number of offspring of the vertex $u \in \textbf{t},$ which will be  called a \textit{leaf} if $l_u(\textbf{t})=0$.
The empty string $\varnothing$  is called {\it the root} of the tree and represents the founding ancestor, which is the only particle in generation 0.
The height of  $\textbf{t}$  is defined by  $|\textbf{t}|:= \sup\{|u|: u \in \textbf{t}\}$ and  the population size at generation $n\geq 0$ is  denoted by 
$$X_n(\textbf{t}):=\mbox{Card}(\{u\in\textbf{t}: |u|=n\}).$$ 
We denote by  $\mathcal{T} $ the set of all rooted ordered trees and by $\mathcal{T}_n=\{\textbf{t}\in \mathcal{T}: |\textbf{t}|\leq n\}$ the restriction of  $\mathcal{T}$  to trees with height  less or equal to $n$. 

\begin{figure}
	\includegraphics[width=9cm]{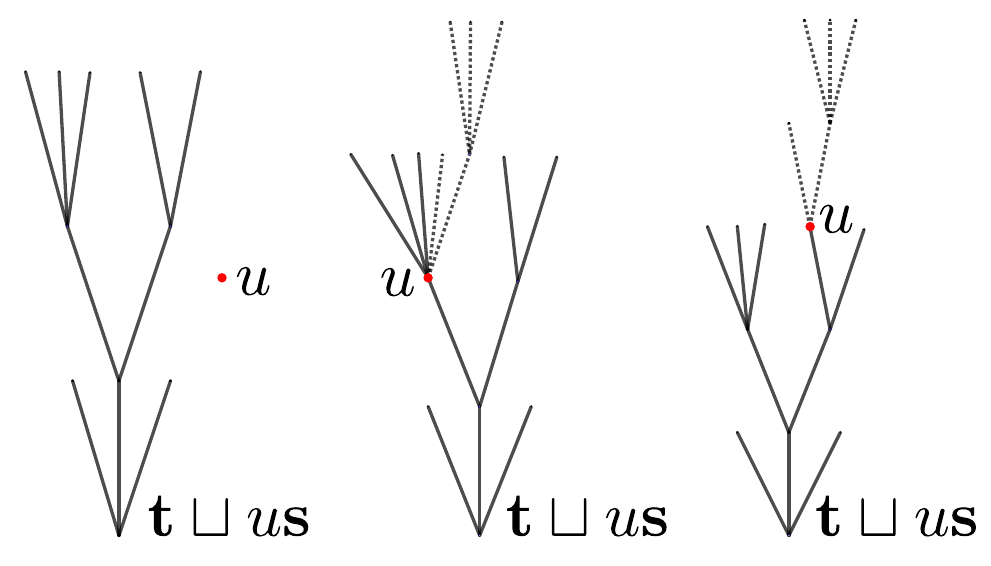}
\caption{Example of the concatenation of $\textbf{t}$ (solid lines) and $\textbf{s}$ (dotted lines) at particle $u$ (red vertex). In the left hand-side picture, $u\notin \textbf{t}$. In the middle picture, $u\in \textbf{t}$ and $u$ is not a leaf. In the right hand-side picture, $u\in \textbf{t}$ and $u$ is a leaf.}
\label{fig:1}
\end{figure}

For our purposes, we introduce   the following operations among rooted trees: the restriction to the first $n$ generations, the subtree attached at  some vertex and the concatenation of trees at a leaf of one of them. More precisely, For $\t\in \mathcal{T}$,  we define its restriction to the first $n$ generation by
$$R_n(\t):=\{u\in\t: |u|\leq n\}, \qquad n\geq 0.$$
Note that $R_n(\t)\in\mathcal{T}_n$. For  $u \in \mathcal{U}$, we define the subtree $S_u(\textbf{t})$ of $\t$ at $u$ (sometimes called \textit{subtree attached at $u$}) as follows
\begin{equation*}
S_u(\textbf{t}):=
\left\{ 
\begin{array}{ll}
\{v\in \mathcal{U}:  uv\in \textbf{t}\} & \text{ if } u\in \textbf{t},\\
\ \emptyset& \text{ if } u\notin \textbf{t}.
\end{array}
\right.
\end{equation*}
See Figure \ref{fig:4}  for an illustrative example of the functionals $R_\cdot$ and $S_\cdot$. We emphasise   the difference between $\emptyset$,  the empty set,  and $\varnothing$,  the root. Additionally,  $S_u(\textbf{t})\in \mathcal{T}$ if $u\in \textbf{t}$. Finally, let $\textbf{t},\textbf{s}\in \mathcal{T}$ and $u\in \t$. We define the \textit{concatenation} of $\textbf{t}$ and $\textbf{s}$ at position $u$ as 
$$\textbf{t}\bigsqcup u \textbf{s}:=\left\{v\in \mathcal{U}: v\in\textbf{t} \mbox{ or }  v=u(l_u(\textbf{t})+w_1, w_2,\dots, w_n) \mbox{ for } w=(w_1, w_2,\dots, w_n) \in \textbf{s}\right\}.$$
If $u$ is a leaf, the concatenation of $\textbf{t}$ and $\textbf{s}$ is much simpler. Indeed, $u$ is a leaf then  $\textbf{t}\bigsqcup u \textbf{s}=\{v\in \mathcal{U}: v\in\textbf{t} \mbox{ or }  v \in u\textbf{s}\}.$
If $u\notin\textbf{t}$, we can still define the concatenation as $\textbf{t}\bigsqcup u \textbf{s}:=\textbf{t}.$ See Figure \ref{fig:1}. The operation of concatenation is  not associative, but we will denote 
$\textbf{t}\bigsqcup u \textbf{s}\bigsqcup v \textbf{r}:=\left(\textbf{t}\bigsqcup u \textbf{s}\right)\bigsqcup v \textbf{r}$.

\begin{figure}
	\includegraphics[width=8cm]{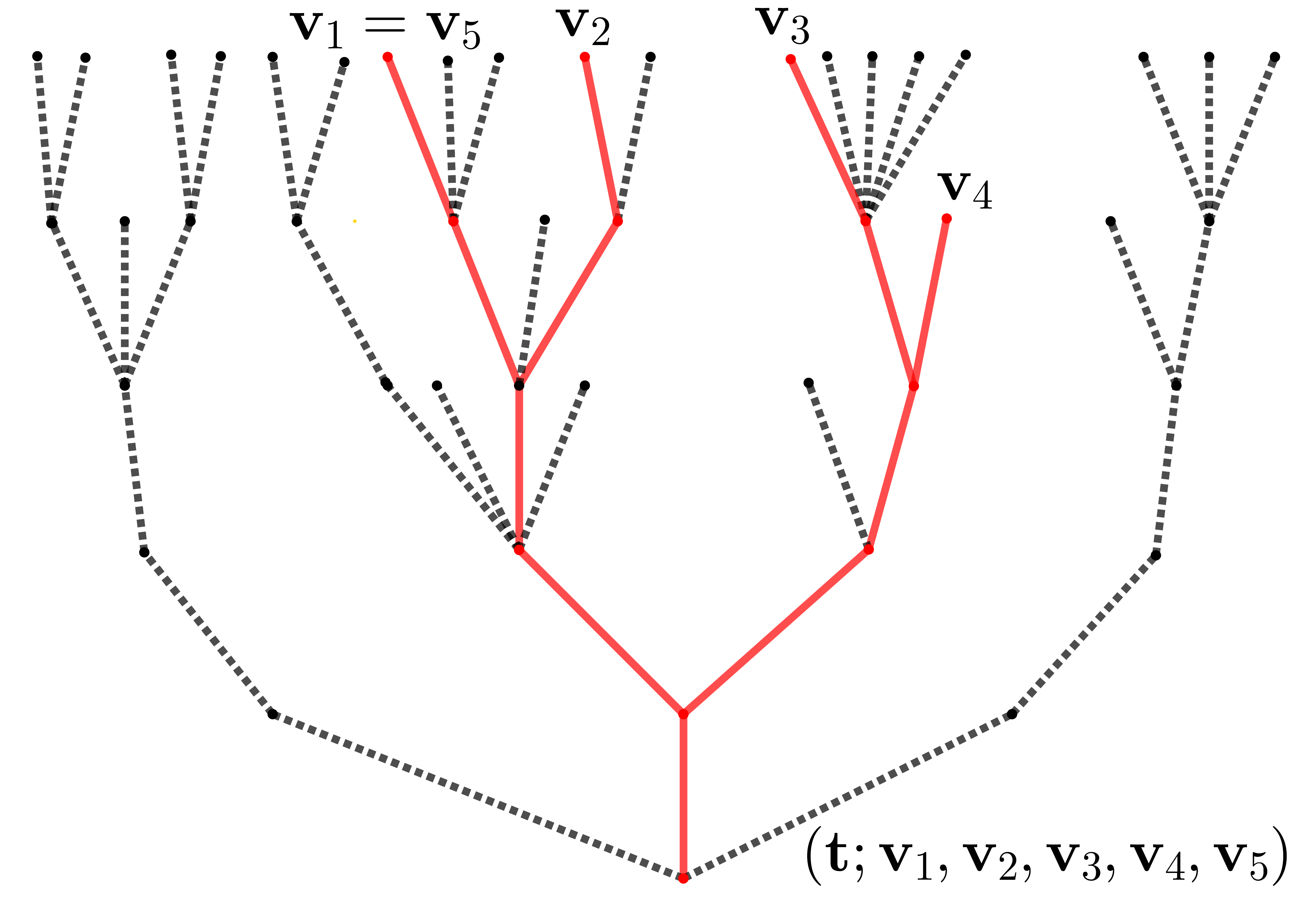}
	\caption{Example of a tree with 5 spines.}
	\label{fig:2}
\end{figure}

Let $\textbf{t}\in \mathcal{T}$, a \textit{spine} $\textbf{v}$  is the genealogical line of a leaf $v\in \textbf{t}$, i.e. $\textbf{v}:=\llbracket\varnothing,v \rrbracket.$ 
The height of $\textbf{v}$ is $|\textbf{v}|=|v|$. For $m\leq |\textbf{v}|$, we denote by $v^{(m)}$ the unique ancestor of $v$ with height $m$. 
Let us denote by 
$$\mathcal{T}^k:= \{(\textbf{t};\textbf{v}_1,\dots, \textbf{v}_k):\ \textbf{t}\ \mbox{is a tree and }  \textbf{v}_i \mbox{ is a spine on } \textbf{t} \mbox{ for all } i\leq k \}, $$ 
the {\it set of trees with $k$ spines} and by   $\mathcal{T}^k_n:=\{(\textbf{t};\textbf{v}_1,\dots, \textbf{v}_k)\in \mathcal{T}^k: |\textbf{t}|=n\}$  its restriction to trees in $\mathcal{T}_n$. See Figure \ref{fig:2}.

We observe that the spines are a special type of rooted trees, the trees with only one leaf. In other words, we can have the following operations: its restriction to the first $n$ generations $R_n(\textbf{v})$ and the subspine $S_u(\textbf{v})$ attached at some $u\in \mathcal{U}$. 
For a $(\textbf{t};\textbf{v}_1,\dots, \textbf{v}_k)\in\mathcal{T}^k$, we define its restriction to the first $n$- generations as follows
$$R_n(\textbf{t};\textbf{v}_1,\dots, \textbf{v}_k):=(R_n(\textbf{t});R_n(\textbf{v}_1),\dots, R_n(\textbf{v}_k))\qquad n\geq 0.$$
Note that $R_n(\textbf{t};\textbf{v}_1,\dots, \textbf{v}_k)\in\mathcal{T}_n^k$. 
Let $(\textbf{t};\textbf{v}_1,\dots, \textbf{v}_k)\in \mathcal{T}^k$, and $u\in \textbf{t}$. We define the marks (or spines) of $u$  by
$$M_u(\textbf{t};\textbf{v}_1,\dots, \textbf{v}_k):=\{\textbf{v}_i: u\in \textbf{v}_i\}.$$
If there is not confusion, we will only use the notation $M_u(\textbf{t})$. We note that $R_{|u|}(\textbf{v})=R_{|u|}(\textbf{w})$ for all $\textbf{v},\textbf{w}\in M_u(\textbf{t})$.
We define the subtree and subspine attached at $u$ by
$$S_u(\textbf{t};\textbf{v}_1,\dots, \textbf{v}_k):=(S_u(\textbf{t});S_{u}(\textbf{v}_{1}), \dots, S_u(\textbf{v}_{k})).$$
More precisely, $M_u(\textbf{t})=\emptyset$, then $S_u(\textbf{t};\textbf{v}_1,\dots, \textbf{v}_k)=S_u(\textbf{t}).$ On the other hand, if $M_u(\textbf{t})=\{\textbf{v}_{i_1},\dots,\textbf{v}_{i_m}\}$, then $S_u(\textbf{t};\textbf{v}_1,\dots, \textbf{v}_k)=(S_u(\textbf{t});S_{u}(\textbf{v}_{i_1}), \dots, S_u(\textbf{v}_{i_m})).$

The  concatenation of two spines $\textbf{v}=\llbracket\varnothing,v\rrbracket$ and $\textbf{w}=\llbracket\varnothing,w\rrbracket$ at $u$ can be performed  for any $u,v,w\in \mathcal{U}$ but  if $u\in \textbf{v} $ is not a leaf, then  the result is a tree with two leaves $v$ and $uw$, not a spine. Thus, we  restrict to concatenate spines only when 
$u\notin \textbf{v}$ or when $u\in \textbf{v}$ is a leaf (i.e. $u=v$). In these cases, we define the \textit{spine  concatenation} as
\begin{equation*}
\textbf{v}\bigsqcup u\textbf{w}:=
\left\{ 
\begin{array}{ll}	
\llbracket\varnothing,v\rrbracket & \text{ if } u\notin \textbf{v},\\
\llbracket\varnothing,uw\rrbracket& \text{ if } u=v.
\end{array}
\right.
\end{equation*}

\begin{figure}
	\includegraphics[width=12cm]{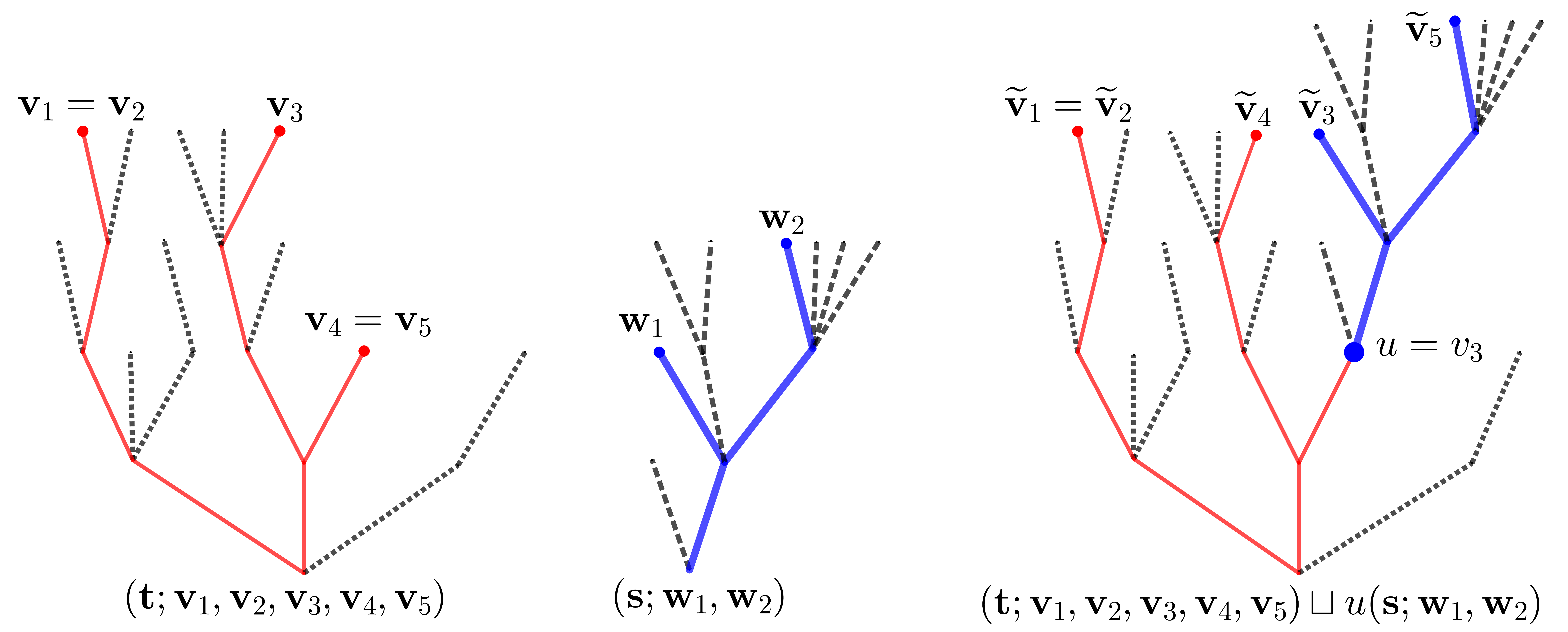}
	\caption{Example of the spine concatenation of $(\textbf{t};\textbf{v}_1,\dots, \textbf{v}_5)$ and $(\textbf{s};\textbf{w}_1, \textbf{w}_2) $  at leaf $v_3$.}
	\label{fig:3}
\end{figure}

Finally, we define the concatenation of two trees with spines at a leaf. 
Let $(\textbf{t};\textbf{v}_1,\dots, \textbf{v}_k)\in\mathcal{T}^k$,  $(\textbf{s};\textbf{w}_1,\dots, \textbf{w}_m)\in\mathcal{T}^m$ and $u\in \textbf{t}$ be a leaf such that $M_u(\textbf{t})=\{\textbf{v}_{i_1},\dots,\textbf{v}_{i_m}\}$.    
We define the concatenation of $(\textbf{t};\textbf{v}_1,\dots, \textbf{v}_k)$ and $(\textbf{s};\textbf{w}_1,\dots, \textbf{w}_m)$ at $u$ as
$$(\textbf{t};\textbf{v}_1,\dots, \textbf{v}_k)\bigsqcup u(\textbf{s};\textbf{w}_1,\dots, \textbf{w}_m):=\left(\textbf{t}\bigsqcup u \textbf{s}; \widetilde{\textbf{v}}_1,\dots, \widetilde{\textbf{v}}_k\right),$$
where $\widetilde{\textbf{v}}_i=\textbf{v}_i$ if $\textbf{v}_i\notin M_u(\textbf{t})$ and $\widetilde{\textbf{v}}_{i_j}=\textbf{v}_{i_j}\sqcup u\textbf{w}_j,$ for $\textbf{v}_{i_j}\in M_u(\textbf{t})$. We observe that $(\textbf{t};\textbf{v}_1,\dots, \textbf{v}_k)\bigsqcup u(\textbf{s};\textbf{w}_1,\dots, \textbf{w}_m)\in \mathcal{T}^k$, see for an illustrative example at Figure \ref{fig:3}.

\subsection{Galton-Watson  trees in varying environments with  spines}

In this subsection, we  study probability measures on the set of  rooted trees which describe  the genealogies of  Galton-Watson processes in  varying  environments. Such trees are known as  Galton-Watson trees in varying  environments (GWTVE). Recall that in a Galton-Watson process in an environment  $e=(q_1,q_2,\ldots)$, any particle in generation  $i$  gives birth to particles in generation  $i+1$  according to  the law $q_{i+1}$. 
Therefore, we say that $(\textbf{T},\p^{(e)})$ is  a {\it Galton-Watson tree in the environment}  $e$ whenever $\textbf{T}$ is a $\mathcal{T}$-valued random variable under law $\p^{(e)}$ satisfying
$$\textbf{G}^{(e)}_n(\textbf{t}):= \p^{(e)}(R_n(\textbf{T})=R_n(\textbf{t})) = \prod_{u\in \textbf{t}:\, |u|<n}  q_{|u|+1}(l_u(\textbf{t})),$$
for all  $n\geq 0$\  and all trees  $\textbf{t}\in\mathcal{T}$. 
Observe that the population size process  $Z=(Z_n:n\geq 0)$  defined by  $Z_n:=X_n(\textbf{T})$ is a \emph{Galton-Watson process in an environment  $e$}. In other words, $((X_m(\textbf{T}), m\leq n), \textbf{G}_n^{(e)})$ has the same law as $((Z_m, m\leq n),\p^{(e)})$.

A \textit{tree in an environment  $e$ with $k$ spines}, under law $\p^{(e)}$, can be constructed as follows:
\begin{enumerate}
\item[(i)] We start with one particle which carries $k$ marks.
\item[(ii)] Particles in generation $m$ give birth to particles in generation  $m+1$  according to  $q_{m+1}$. 
\item[(iii)] If a particle with $j$ marks gives birth to $a>0$ particles, then, the $j$ marks choose a particle to follow independently and uniformly at random from amongst the $a$ available.
\item[(iv)] If a particle with $j$ marks gives birth to $0$ particles, then its marks remain with it as it moves to the graveyard state $\Delta$.
\end{enumerate}

With this construction in hand, we introduce  a {\it Galton-Watson tree in an environment  $e$ with  $k\ge 1$ spines}. The tree with $k$-spines  $(\textbf{T};\textbf{V}_1,\ldots,\textbf{V}_k)$ is a $\mathcal{T}^{k}$-valued r.v. with the following distribution 
\begin{equation*}
\begin{split}
\p_n^{(e,k)}((\textbf{t};\textbf{v}_1,\ldots,\textbf{v}_k)):=&\ \p^{(e)}\left(
R_n(\textbf{T};\textbf{V}_1,\ldots,\textbf{V}_k)=R_n(\textbf{t};\textbf{v}_1,\ldots,\textbf{v}_k)\right)\\ 
=& \prod_{u\in \textbf{t}:\, |u|<n}  q_{|u|+1}(l_u(\textbf{t}))\prod_{i=1}^k\prod_{u\in\textbf{v}_i:\, |u|<|\textbf{v}_i|\wedge n}\frac{1}{l_u(\textbf{t})}\\
=&\prod_{i=1}^k\prod_{u\in\textbf{v}_i:\, |u|<|\textbf{v}_i|\wedge n}\frac{1}{l_u(\textbf{t})}\mathbf{G}_n^{(e)}(\textbf{t}),
\end{split}
\end{equation*}
for any $n\geq 0$ and   $(\textbf{t};\textbf{v}_1,\ldots,\textbf{v}_k)\in \mathcal{T}^k$.  
Observe from the above definition, that the first  and second products follows from 
steps (iii) and  (ii), respectively,  from our construction of   trees in an environment $e$.  It is important to note that we can extend the definition to  the case $k=0$, by just taking a Galton-Watson tree in environment $e$, that is 
$$\p_n^{(e,0)}(\textbf{t}):=\mathbf{G}_n^{(e)}(\textbf{t}),\qquad \textbf{t}\in\mathcal{T}.$$

For an environment $e=(q_1,q_2,\dots )$ and $m\in \N_0$, we set the shift environment as
$$e_m:=(q_{m+1},q_{m+2},\dots ).$$
Our first result, in this section, states that the law $\p_n^{(e,k)}$ satisfies a time-dependent Markov branching property, in that the descendants of any particle behave independently of the rest of the tree.
\begin{proposition}\label{prop:branching property}\emph{\textbf{Markov property for GWTVE with spines under $\p_n^{(e,k)}$.}}
	Let $e$ be an environment and  $(\textbf{T};\textbf{V}_1,\ldots,\textbf{V}_k)$ a r.v. taking values in  $\mathcal{T}_n^k$ with distribution $\p_n^{(e,k)}$. Suppose that a particle $u\in \textbf{T}$ satisfies $|M_u(\textbf{T})|=j$ and  $|u|=m$, for some $0\leq j\leq k$ and $0\leq m<n$. Then, $S_u(\textbf{T};\textbf{V}_1,\ldots,\textbf{V}_k)$ is independent of the rest of the system, 
	that is, from the sigma-algebra  
	$$\mathcal{G}:=\sigma(\{R_m(\textbf{T};\textbf{V}_1,\ldots,\textbf{V}_k)\}\cup\{S_w(\textbf{T};\textbf{V}_1,\ldots,\textbf{V}_k): w \in \textbf{T}, |w|=m, w\neq u\}).$$
	In addition,  the law of $S_u(\textbf{T};\textbf{V}_1,\ldots,\textbf{V}_k)$ is given by $\p_{n-m}^{(e_m,j)}$.
\end{proposition}

\begin{proof}
	Independence  holds true since by construction all the subtrees with (possible zero) spines  starting at generation $m$ are independent. On the other hand, the  environment starting at generation $m$ is the shifted environment $e_m=(q_{m+1},q_{m+2},\dots)$ and from our hypotheses the subtree starting at $u$ has  $j$ marks  and  $n-m$ generations.  Therefore, the distribution of  the subtree and spines starting at $u$ is given by $\p_{n-m}^{(e_m,j)}$ which implies our result.
\end{proof}

Let $k\ge 1$ and denote by 
$$\widehat{\mathcal{T}}_n^k=\{(\textbf{t};\textbf{v}_1,\dots, \textbf{v}_k)\in \mathcal{T}^k_n: |\textbf{v}_i|=n \mbox{ and } \textbf{v}_i\neq \textbf{v}_j \mbox{ for all }  1\leq i\neq j \leq k \}$$
the subspace of $\mathcal{T}^k_n$ where all  spines are distinct and alive at time $n$. Now, we want to define a new probability measure by conditioning $((\textbf{T};\textbf{V}_1,\ldots,\textbf{V}_k), \mathbb{P}_n^{(e,k)})$ to the subspace $\widehat{\mathcal{T}}_n^k$ in such a way that  spines are choosing uniformly from all  distinct particles  alive at time $n$.  Recall that $X_n(\textbf{t})=\#\{u\in\textbf{t}: |u|=n\} $ is the population size at the  $n$-th generation of the tree $\textbf{t}$. 

Let $\theta\in [0,1)$ and define the function $h_{n, k, \theta}:\mathcal{T}^k\rightarrow \mathbb{R}$, as follows
\begin{equation}\label{fctg}
h_{n,k,\theta}(\textbf{t};\textbf{v}_1,\dots, \textbf{v}_k):=\indi_{\{\textbf{v}_i\neq \textbf{v}_j, i\neq j\}}\indi_{\{|\textbf{v}_i|=n , i\leq k\}}e^{-\theta X_{n}(\textbf{t})}\prod_{i=1}^k\prod_{u\in\textbf{v}_i:\, |u|<n}l_u(\textbf{t}).
\end{equation}
We also introduce the law $\mathbf{Q}_n^{(e,k, \theta)}$ by
\begin{equation}\label{cmg}
\begin{split}
\mathbf{Q}_n^{(e,k, \theta)}((\textbf{t};\textbf{v}_1,\ldots,\textbf{v}_k))
:=&\frac{h_{n,k,\theta}(\textbf{t};\textbf{v}_1,\dots, \textbf{v}_k)\p^{(e,k)}_n((\textbf{t};\textbf{v}_1,\ldots,\textbf{v}_k))}{\e^{(e,k)}_n[h_{n,k ,\theta}(\textbf{T};\textbf{V}_1,\ldots,\textbf{V}_k)]},
\end{split}
\end{equation}
for $(\textbf{t};\textbf{v}_1,\dots, \textbf{v}_k)\in \widehat{\mathcal{T}}^k_n$. Recall that $x^{[k]}:=x(x-1) \times\cdots \times(x-k+1)$, for $x\in\mathbb{N}$. 
Then we find the following result:

\begin{proposition}\label{prop:measure Q}
Let $((\textbf{T};\textbf{V}_1,\ldots,\textbf{V}_k),\p^{(e,k)}_n)$  be a Galton-Watson tree up to time  $n$  with  $k$ spines. Then,
$$\e^{(e,k)}_n[h_{n,k,\theta}(\textbf{T};\textbf{V}_1,\ldots,\textbf{V}_k)]=\e^{(e)}\left[ Z_n^{[k]}
e^{-\theta Z_n}\right].$$
In particular,
$$\mathbf{Q}_n^{(e,k,\theta)}((\textbf{t};\textbf{v}_1,\ldots,\textbf{v}_k))
=\frac{\indi_{\{\textbf{v}_i\neq \textbf{v}_j, i\neq j\}}\indi_{\{|\textbf{v}_i|= n, i\leq k\}}e^{-\theta X_{n}(\textbf{t})}}{\e^{(e)}\left[Z_n^{[k]}
e^{-\theta Z_n}\right]}\mathbf{G}_n^{(e)}(\textbf{t}).$$
\end{proposition}

Proposition \ref{prop:measure Q} reveals that if the population size $Z_n$ has finite $k$-th moments then we do not need any discounting term and can simply take $\theta=0$ in the above definition. In such cases, we denote  $\mathbf{Q}_n^{(e,k)}:=\mathbf{Q}_n^{(e,k, 0)}$. However, in general the offspring distributions and population sizes may not have finite $k$-th moments, in which case we must take $\theta>0$ in order to guarantee that the denominator is finite.  Importantly, the introduction of discounting in the change of measure turns out to be very natural. Crucially, it will allow our approach to work elegantly in generality, that is, without any additional moment assumptions or ad hoc truncations (as needed in \cite{harris2020coalescent}).

The new measure $\mathbf{Q}_n^{(e,k, \theta)}$ is defined via a change of measure relative to $\p^{(e,k)}_n$ only at time $n$. Of course, we could determine the martingale corresponding to this change of measure, $M_n=(M_n(m))_{m=0,\dots,n}$, by projecting  $h_{n,k,\theta}(\textbf{T};\textbf{V}_1,\ldots,\textbf{V}_k)$ onto the natural filtration at times $m\in\{0,1,\dots,n\}$. However, for our purposes, we will not need to compute this martingale directly and instead calculate projections only for certain events related to properties of interest (e.g. to describe the first split of the spines in Proposition \ref{prop: offspring spine}).

\begin{proof}[Proof of Proposition \ref{prop:measure Q}]
By definition of $\p^{(e,k)}$ and $h_{n,k,\theta}$, the second identity holds true as soon as we prove the first one. In addition, we also deduce 
$$\e^{(e,k)}_n[h_{n,k,\theta}(\textbf{T};\textbf{V}_1,\ldots,\textbf{V}_k)]=\underset{(\textbf{t};\textbf{v}_1,\dots, \textbf{v}_k)\in\mathcal{T}^k_n}{\sum}
\indi_{\{\textbf{v}_i\neq \textbf{v}_j, i\neq j\}}\indi_{\{|\textbf{v}_i|= n, i\leq k\}}e^{-\theta X_{n}(\textbf{t})}\mathbf{G}_n^{(e)}(\textbf{t}).$$
For a fixed tree $\textbf{t}\in\mathcal{T}_n $, the term associated with  $(\textbf{t};\textbf{v}_1,\dots, \textbf{v}_k)$ only contributes to the sum when all spines $(\textbf{v}_1,\dots, \textbf{v}_k)$ are different and have height equals $n$. In this case, it contributes with the value $e^{-\theta X_{n}(\textbf{t})}$. On the other hand, we observe that there are $X_n(\textbf{t})(X_n(\textbf{t})-1)\times\cdots \times(X_n(\textbf{t})-k+1)$ different possibilities to choose these $k$ spines from those particles alive at time $n$, i.e. $X_{n}(\textbf{t})$. Therefore,
\begin{equation}\label{sum with spines}
\begin{split}
\underset{(\textbf{t};\textbf{v}_1,\dots, \textbf{v}_k)\in\mathcal{T}^k_n}{\sum}&e^{-\theta X_{n}(\textbf{t})}\indi_{\{\textbf{v}_i\neq \textbf{v}_j, i\neq j\}}\indi_{\{|\textbf{v}_i|= n, i\leq k\}}\mathbf{G}_n^{(e)}(\textbf{t})\\
=&\underset{\textbf{t}\in\mathcal{T}_n}{\sum}e^{-\theta X_{n}(\textbf{t})}X_n(\textbf{t})(X_n(\textbf{t})-1) \times\cdots \times(X_n(\textbf{t})-k+1)\mathbf{G}_n^{(e)}(\textbf{t}).
\end{split}
\end{equation}
Since  $\{(X_m(\textbf{T}), m\leq n), \textbf{G}_n^{(e)}\}$ has the same law as $\{(Z_m, m\leq n),\p^{(e)}\}$, we  deduce the first identity in the statement. This completes the proof. 
\end{proof}

From the construction of  $\mathbf{Q}_n^{(e,k,\theta)}$, we see that $\mathbf{Q}_n^{(e,k,\theta)}(X_n(\textbf{t})\geq k)=1$. In other words,  under $\mathbf{Q}_n^{(e,k,\theta)}$, there are at least $k$ different particles at generation $n$. 
Importantly, under 
$\mathbf{Q}_n^{(e,k,\theta)}$, the $k$ spines are a uniform choice without replacement from all particles alive at  generation $n$. 
In addition, under $\mathbf{Q}_n^{(e,k,\theta)}$, the total population size at time $n$, $Z_n$ is the $k$-size biased and $\theta$-discounted transformation of its distribution under $\p^{(e)}$.

As before, we extend the definition of $\mathbf{Q}_n^{(e,k,\theta)}$ for $k=0$, by only performing the discounting. In other words,  $h_{n,0,\theta}:\mathcal{T}\rightarrow \mathbb{R}$ is given by
$h_{n,0,\theta}(\textbf{t}):=e^{-\theta X_{n}(\textbf{t})}$ and 
\begin{equation*}
\begin{split}
\mathbf{Q}_n^{(e,0, \theta)}(\textbf{t})
:=&\frac{e^{-\theta X_{n}(\textbf{t})}}{\e^{(e)}[e^{-\theta Z_{n}}]}\mathbf{G}_n^{(e)}(\textbf{t}).
\end{split}
\end{equation*}
Our next result says that the Markov branching property given in  Proposition \ref{prop:branching property} is inherited by the measure $\mathbf{Q}_n^{(e,k,\theta)}$.  
In particular, it implies  that if we have two subtrees such that neither of their roots are part of the other subtree, then the subtrees are independent.


\begin{proposition}[\emph{\textbf{Markov property for GWTVE with spines under $\mathbf{Q}_n^{(e,k,\theta)}$}}]\label{prop: branching property Q}
For $k\geq 1$, let $(\textbf{T};\textbf{V}_1,\ldots,\textbf{V}_k)$ be a $\mathcal{T}_n^k$-valued random variable with law $\mathbf{Q}_n^{(e,k,\theta)}$. Suppose that a particle $u\in \textbf{T}$ satisfies $|M_u(\textbf{T})|=j$ and $|u|=m$ for some $0\leq j\leq k$ and $0\leq m<n$ and let
$$\mathcal{G}:=\sigma(\{R_m(\textbf{T};\textbf{V}_1,\ldots,\textbf{V}_k)\}\cup\{S_w(\textbf{T};\textbf{V}_1,\ldots,\textbf{V}_k): w \in \textbf{T}, |w|=m, w\neq u\}).$$
Then,  $S_u(\textbf{T};\textbf{V}_1,\ldots,\textbf{V}_k)$ is independent of  $\mathcal{G}$ and its law is given by $\mathbf{Q}_{n-m}^{(e_m,j,\theta)}$.
\end{proposition}

\begin{proof}  We first treat  the case $j\geq 1$. Let $(\textbf{s};\textbf{w}_1,\dots,\textbf{w}_j)\in \mathcal{T}_{n-m}^j$ and denote by $\mathcal{G}$ for the $\sigma$-algebra generated by  the whole information of $(\textbf{T};\textbf{V}_1,\ldots,\textbf{V}_k)$ except for $S_u(\textbf{T};\textbf{V}_1,\ldots,\textbf{V}_k)$.  
By the conditioned version of the Radon Nikodym theorem (see  \cite[Lemma 13]{harris2020coalescent}), we  have
\begin{equation}\label{subtree_at_u}
\begin{split}
\mathbf{Q}_n^{(e,k,\theta)}&\left(S_u(\textbf{T};\textbf{V}_1,\ldots,\textbf{V}_k)=(\textbf{s};\textbf{w}_1,\dots,\textbf{w}_j)\mid \mathcal{G}\right)\\
&\hspace{2cm}=\frac{\e^{(e,k)}_n[h_{n,k,\theta}(\textbf{T};\textbf{V}_1,\ldots,\textbf{V}_k)
\indi_{S_u(\textbf{T};\textbf{V}_1,\ldots,\textbf{V}_k)=(\textbf{s};\textbf{w}_1,\dots,\textbf{w}_j)}\mid \mathcal{G}]}{\e^{(e,k)}_n[h_{n,k,\theta}(\textbf{T};\textbf{V}_1,\ldots,\textbf{V}_k)
\mid \mathcal{G}]}.
\end{split}
\end{equation}
Let us enumerate by $u_1,\dots ,u_{X_m(\textbf{t})-1}$  the particles at generation $m$ which  are different to $u$. Next, we  observe that  $(\textbf{T};\textbf{V}_1,\ldots,\textbf{V}_k)$ can be rewritten as the concatenation of its restriction to the  $m$-th generation  with the subtrees which are attached at each particle in generation $m$ (see Figure \ref{fig:4}),  that is 
\begin{equation*}
(\textbf{T};\textbf{V}_1,\ldots,\textbf{V}_k)=R_m (\textbf{T};\textbf{V}_1,\ldots,\textbf{V}_k)\bigsqcup u S_u(\textbf{T};\textbf{V}_1,\ldots,\textbf{V}_k)
\underset{r=1}{\overset{X_m(\textbf{t})-1}{\bigsqcup}} u_r S_{u_r}(\textbf{T};\textbf{V}_1,\ldots,\textbf{V}_k).
\end{equation*}
	
\begin{figure}
\includegraphics[width=12cm]{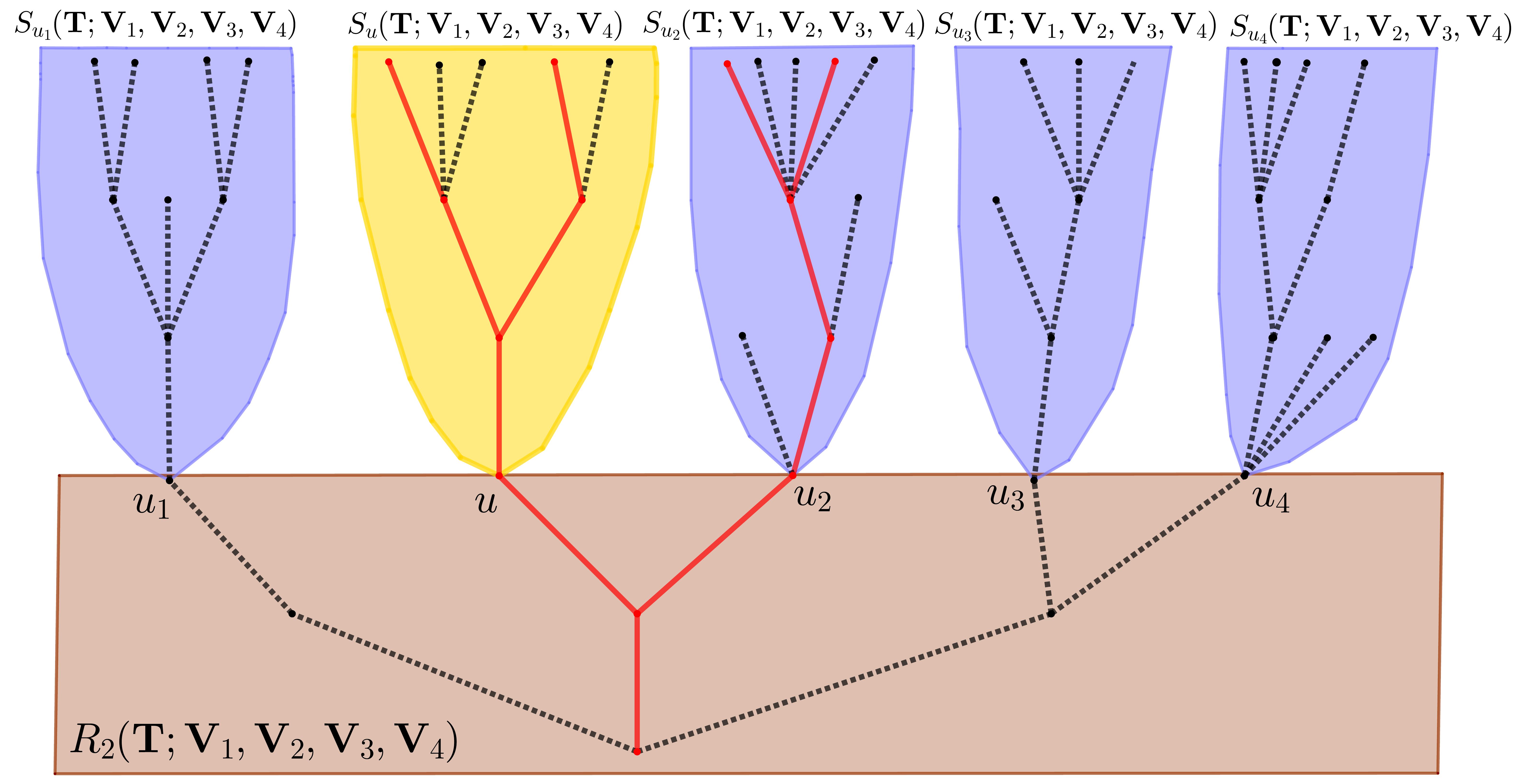}
\caption{Decomposition of the tree with spines $(\textbf{T};\textbf{V}_1,\textbf{V}_2,\textbf{V}_3,\textbf{V}_4)$ at generation $m=2$ into its restriction up to generation $2$ concatenated  with the subtrees attached to particles $u, u_1,u_2,u_3$ and $u_4$. Particles without marks have dotted lines. The spines are the solid red  lines.}
\label{fig:4}
\end{figure}

Therefore, we can decompose the  population size at the  $n$-th generation of the tree  $\textbf{t}$  as follows
\begin{equation*}
X_n(\textbf{t})=X_{n-m}(S_u(\textbf{T}))+\sum_{r=1}^{X_m(\textbf{t})-1}X_{n-m}(S_{u_r}(\textbf{T})).
\end{equation*}
Denote by  $I:=\{i\leq k: \textbf{V}_i\in M_u(\textbf{t})\}$. With the previous decompositions, we observe
$$
h_{n,k,\theta}(\textbf{T};\textbf{V}_1,\ldots,\textbf{V}_k)=h_{n-m,k,\theta}(S_u(\textbf{T};\textbf{V}_1,\ldots,\textbf{V}_k))H,$$
where $ h_{n,k,\theta}$ is  defined in \eqref{fctg} and 
\[
\begin{split}
H&
:=\prod_{r=1}^{X_m(\textbf{t})-1}e^{-\theta X_{n-m}(S_{u_r}(\textbf{T}))}
\left( \prod_{i\in I}
\prod_{w\in\textbf{V}_i: |w|<m}l_w(\textbf{T})\right)\\
&\hspace{4cm} \times\indi_{\{\textbf{V}_j\neq \textbf{V}_\ell, j\neq \ell\in I^c\}}
\prod_{j\in I^c}\indi_{\{|\textbf{V}_j|=n\}}\prod_{v\in\textbf{V}_r: |v|<n}l_v(\textbf{T}),
\end{split}
\]
which is an $\mathcal{G}$-measurable function. Now, we use equation \eqref{subtree_at_u} to obtain
\begin{equation*}
\begin{split}
&\mathbf{Q}_n^{(e,k,\theta)}\left(S_u(\textbf{T};\textbf{V}_1,\ldots,\textbf{V}_k)=(\textbf{s};\textbf{w}_1,\dots,\textbf{w}_j)\mid \mathcal{G}\right)\\	&\hspace{2cm}=\frac{\e^{(e,k)}_n[h_{n-m,k,\theta}(S_u(\textbf{T};\textbf{V}_1,\ldots,\textbf{V}_k))H
\indi_{S_u(\textbf{T};\textbf{V}_1,\ldots,\textbf{V}_k)=(\textbf{s};\textbf{w}_1,\dots,\textbf{w}_j)}\mid \mathcal{G}]}{\e^{(e,k)}_n[h_{n-m,k,\theta}(S_u(\textbf{T};\textbf{V}_1,\ldots,\textbf{V}_k))H
\mid \mathcal{G}]}\\
&\hspace{2cm}=\frac{\e^{(e_m,j)}_{n-m}[h_{n-m,k,\theta}(\textbf{S};\textbf{W}_1,\ldots,\textbf{W}_j)
\indi_{(\textbf{S};\textbf{W}_1,\ldots,\textbf{W}_j)=(\textbf{s};\textbf{w}_1,\dots,\textbf{w}_j)}]}{\e^{(e_m,j)}_{n-m}[h_{n-m,k,\theta}(\textbf{S};\textbf{W}_1,\ldots,\textbf{W}_j)]}\\
&\hspace{2cm}=\mathbf{Q}_{n-m}^{(e_m,j,\theta)}\left(\left(\textbf{s};\textbf{w}_1,\dots,\textbf{w}_j\right)\right),
\end{split}
\end{equation*}
where in the second identity, we have used that $H$ is $\mathcal{G}$ measurable and  Proposition \ref{prop:branching property}, i.e.   the law of $S_u(\textbf{T};\textbf{V}_1,\ldots,\textbf{V}_k)$ under $\p_n^{(e,k)}$ is given by $\p_{n-m}^{(e_m,j)}$ and we can take $S_u(\textbf{T};\textbf{V}_1,\ldots,\textbf{V}_k)=(\textbf{S};\textbf{W}_1,\ldots,\textbf{W}_j)$. In the last identity, we have used the definition of $\mathbf{Q}_{n-m}^{(e_m,j,\theta)}$.
For the case $j=0$, we follow the same arguments. Indeed, we only require to change  $\textbf{s}\in \mathcal{T}_{n-m}$ at the beginning, take  $\textbf{S}:=S_u(\textbf{T};\textbf{V}_1,\ldots,\textbf{V}_k)$ in  the middle and note that an empty  product is equal to one.
\end{proof} 

Now, our aim is to  provide a construction, forward in time, of the tree under $\mathbf{Q}_n^{(e,k,\theta)}$. Our construction follows from a  recursive procedure.  For the cases $k=1$ and $\theta=0$, such construction can be found  in \cite[Section 1.4]{kersting2017discrete} and for  the case  $k=2$ and $\theta=0$, see  \cite[Section 2]{CardonaPalau}.

Before we proceed with our construction, let us recall some notation.
Recall that  $(f_1,f_2,\dots)$ denote the generating functions associated with  the environment $e=(q_1,q_2,\dots).$ 
For each $0\leq m< n$ and $s\in [0,1]$, we define 
$$ f_{m,n}(s):= [f_{m+1} \circ \cdots \circ f_n](s),$$ 
and  $f_{n,n}(s):=s$, where $f\circ g$ is the composition of $f$ with $g$. According to the monograph of Kersting and Vatutin \cite{kersting2017discrete}, we have
\begin{equation}\label{eq: laplace}
\e^{(e)}\left[s^{Z_n}\mid Z_m=1\right]=\e^{(e_m)}\left[s^{Z_{n-m}}\right]=f_{m,n}(s)\qquad 0\leq m< n \quad \mbox{ and }\quad s\in [0,1].
\end{equation}
In particular, for any $m\geq 0$ 
\begin{equation}
\label{moved environment}
\left\{(Z_{n+m},n\geq 0), \p^{(e)}(\ \cdot\mid Z_m=1)\right\}\overset{(d)}{=}\left\{(Z_{n}, n\geq 0), \p^{(e_m)}\right\},
\end{equation}
and, for any $k\geq 1$,  $0\leq m\leq n$ and $\theta\geq 0$,
\begin{equation}
\label{factorial moments}
\e^{(e_m)}\left[e^{-\theta Z_{n-m}}Z_{n-m}^{[k]}\right]=\e^{(e)}\left[e^{-\theta Z_n}Z_n^{[k]}\mid Z_m=1\right]=e^{-\theta k}\left.\frac{\partial^k}{\partial s^k }f_{m,n}(s)\right|_{s=e^{-\theta}}.
\end{equation}

Let $g\in\{1,\dots,k\}$  be the number of spine groups created at a given time  and let  $k_1,\ldots,k_g\geq 1$ be their sizes. Denote by  $g_r$  the number of groups of size $r$. By a combinatorial argument we have 
\begin{equation}\label{eq: combinatorial}
\underset{i=1}{\overset{g}{\sum}}k_i=k,\qquad  \underset{r=1}{\overset{k}{\sum}}\ g_r=g\qquad \mbox{and} \qquad \underset{r=1}{\overset{k}{\sum}}\ rg_r=k.
\end{equation} 
We also introduce   the first time that the spines split apart  and the last time all the spines are together by $\hat\psi_1$ and  $\psi_1$, respectively. We observe that  $\hat\psi_1=\psi_1+1$. We often swap between using $\psi_1$ and $\hat\psi_1$ as most appropriate. For simplicity of exposition, for $g\geq 2$ and $k=k_1+\dots+k_g$, we denote 
\begin{equation}
\label{def: C}
C_{k;k_1,\dots, k_g}:=\{\mbox{Spines split from a group of size $k$  into }g \mbox{ groups of sizes } k_1,\dots, k_g\}.
\end{equation}

\begin{proposition}[\bf A $k$-spine construction under $\mathbf{Q}_n^{(e,k,\theta)}$]\label{prop: forward construction}
A tree with $k\ge 1$ different spines up to generation $n$, under $\mathbf{Q}_n^{(e,k,\theta)}$,  is constructed as follows.
\begin{enumerate}
\item We start with a particle with $k$ marks. 
\item If $k=1$, we consider $\psi_1=n$. Otherwise, select  $\psi_1$, the number of spine groups $g$ and their sizes $k_1,\dots, k_g$ according to 
\begin{equation}
\label{eq: g groups}
\begin{split}
&\mathbf{Q}^{(e,k,\theta)}_n\left(\psi_1=m, C_{k;k_1,\dots, k_g}\right)\\
&=\frac{k!
\Big(\left.\tfrac{\partial^g}{\partial s^g }f_{m+1}(s)\right|_{s=	f_{m+1,n}(e^{-\theta})}\Big)
\e^{(e)}\left[e^{-\theta Z_n}Z_n\right]
\underset{i=1}{\overset{g}{\prod}}\e^{(e_{m+1})}\left[e^{-\theta Z_{n-m-1}}Z_{n-m-1}^{[k_i]}\right]}
{ \left(\underset{r=1}{\overset{k}{\prod}}g_r!\right) \left( \underset{i=1}{\overset{g}{\prod}}k_i!\right)
\e^{(e_{m})}\left[e^{-\theta Z_{n-m}}Z_{n-m}\right]
\e^{(e)}\left[e^{-\theta Z_n}Z_{n}^{[k]}\right]
},
\end{split}
\end{equation}
where $m\in \{0,\dots,n-1\}$,  $g\in\{2,\dots,k\}$,  $k_1,\ldots,k_g\geq 1$ satisfying  $k_1+\cdots +k_g=k$ and $g_r$ is the number of groups of size $r$.
\item An unmarked particle in generation $m\in \{0,\dots, n-1\}$ gives birth to unmarked particles in generation $m+1$, independently of other particles, with probability  
\begin{equation*}
q_{m+1}^{(0,\theta)}\left(\ell \right):=q_{m+1}(\ell) \frac{f_{m+1,n}(e^{-\theta})^\ell}{\e^{(e_m)}\left[e^{-\theta Z_{n-m}}\right]}, \qquad \ell\geq 0.
\end{equation*}
\item A marked particle in generation $m\in \{0,\dots, \psi_1-1\}$ gives birth to particles in generation $m+1$ accordingly to 
\begin{equation*}
\begin{split}
&q_{m+1}^{(1,\theta)}\left(\ell \right):=
\frac{\ell q_{m+1}(\ell) \left(f_{m+1,n}(e^{-\theta})\right)^{\ell-1} }
{f_{m+1}'(f_{m+1,n}(e^{-\theta}))}, \qquad \ell\geq  1.
\end{split}
\end{equation*}
Uniformly, select one of the particles to carry the $k$ marks. All the other particles remain  unmarked.
\item The marked particle at generation $m=\psi_1<n$ gives birth accordingly to 
$$q_{m+1}^{(g,\theta)}\left(\ell \right):=\ell^{[g]}q_{m+1}(\ell)
\frac{(f_{m+1,n}(e^{-\theta}))^{l-g}}{\left.\tfrac{\partial^g}{\partial s^g }f_{m+1}(s)\right|_{s=f_{m+1,n}(e^{-\theta})}},\qquad \ell\geq g.$$
Uniformly, select $g$ of these particles as  marked in generation $m+1$ with $k_1,\dots, k_g$ marks, respectively. All the other particles remain  unmarked.
\item 
Repeat steps (1)-(5) for each of the $g$ marked particles. Next, we have to construct  $g$ independent trees with $k_1, \dots, k_g$ different spines under the measure   $\mathbf{Q}^{(e_{m+1},k_j,\theta)}_{n-m-1}$, respectively.
\end{enumerate}
\end{proposition}

It is important to note that whenever the process possesses  finite $k$-th moments, we can take $\theta=0$. In this case, the distribution of $\psi_1$ and the spine subgroups can be simplified, actually $q_{m+1}^{(0,0)}=q_{m+1}$, $q_{m+1}^{(1,0)}$ is 
the $1$ size-biased transform of $q_{m+1}$ and  $q_{m+1}^{(g,0)}$   is the $g$ size-biased transform of $q_{m+1}$. This guarantees that at $m\neq \psi_1$ the marked particle has at least one child, and at $m= \psi_1$ the marked particle has at least $g$ children. 
In this case, by denoting $u^*$ for the marked particle at generation $m$,  then the probability $\mathbf{Q}^{(e,k,\theta)}_n\left(\psi_1=m, l_{u^*}=\ell, C_{k;k_1,\dots, k_g} \right)$ can be obtained by combining \eqref{eq: g groups} with step 5 above (see also \eqref{tomate}) to deduce
$$\e^{(e)}\left[Z_m\right]
q_{m+1}(\ell) {\ell \choose g}{g \choose g_1,\dots,g_k}{k \choose k_1,\dots, k_g}
\frac{ \underset{i=1}{\overset{g}{\prod}}\e^{(e_{m+1})}\left[Z_{n-m-1}^{[k_i]}\right]}
{	\e^{(e)}\left[Z_{n}^{[k]}\right]
}. 
$$ 
The latter can be interpreted as  follows. The mean   $\e^{(e)}\left[Z_m\right]$ represents the different ways to select $u^*$. Once the particle  $u^*$ is selected, it has $\ell$ offspring with probability $q_{m+1}(\ell)$, from these offspring there are ${\ell \choose g}$ possibilities to select $g$ of them with marks, ${g \choose g_1,\dots,g_k}$ ways of arranging the values $k_1,\ldots,k_g$ with the sizes of each group and ${k \choose k_1,\dots, k_g}$ ways to arrange the $k$ spines into $g$ groups with sizes $k_1,\ldots,k_g$. For each $1\leq i\leq g$, the quantity $\e^{(e_{m+1})}\left[Z_{n-m-1}^{[k_i]}\right]$ represents the ways to select one subtree  with $k_i$ spines with root in generation $m+1$. Finally, we divide by $\e^{(e)}\left[Z_{n}^{[k]}\right]$ because it is the number of the possible trees with $k$ spines alive at generation $n$.

If the process does not  have finite $k$-th  moments, then necessarily  $\theta>0$. In this scenario the previous formulas have the same interpretation but with the discounting factors.

In order to prove Proposition \ref{prop: forward construction}, we require some intermediate results. We first present such results, together with  their proofs,  and at the end of this section we present the proof of Proposition \ref{prop: forward construction}. Actually, it will follow  from  Proposition \ref{prop: branching property Q} (which is step (6)) and the forthcoming  Propositions \ref{prop: first splitting time}, Proposition  \ref{prop: offspring spine} (Steps 2, 4 and 5) and Corollary \ref{offspring without} (Step (3)).

Recall that  $\psi_1$  and $\hat\psi_1$ denote the last time where all marks (or spines) are together and the first spine splitting time, respectively.  In other words, if $\textbf{V}_1,\dots, \textbf{V}_k$ denote the spines, then   
$$\{\psi_1=m\}=\{R_m(\textbf{V}_1)=R_m(\textbf{V}_i),\forall i\leq k\}\bigcap \{\exists i\neq j,\ i,j\leq k: R_{m+1}(\textbf{V}_i)\neq R_{m+1}(\textbf{V}_j)\}.$$
Observe that, unlike the continuous time case $\psi_1$ and $\hat\psi_1$ are different, indeed $\hat\psi_1=\psi_1+1$ and $\hat\psi_1$ is a stopping time in the natural filtration of the branching process whereas $\psi_1$ is not.

The next result provides the distribution of $\psi_1$, hence implicitly that of $\hat\psi_1$. Observe that  since we  have finite second moment, when $k=2$ we may take $\theta=0$ and use the notation $\mathbf{Q}^{(e,2)}$. In this particular case, the next result  explains the choice of the distribution of $\psi_1$ in 
the two spine decomposition in Cardona and Palau \cite{CardonaPalau} which satisfies
\begin{equation}\label{eq_law_Kn}
\mathbf{Q}_n^{(e,2)}(\psi_1= m) = \frac{\nu_{m+1}}{\mu_{m}}
\left(\sum_{k=0}^{n-1}\frac{\nu_{k+1}}{\mu_{k}}\right)^{-1},
\qquad 0 \leq m \leq n-1.
\end{equation}

\begin{proposition}\label{prop: first splitting time}
Let $0\leq m<n$, then
\begin{equation*}
\begin{split}
\mathbf{Q}^{(e,k,\theta)}_n(\psi_1\geq m)
&=\frac{\e^{(e)}\left[e^{-\theta Z_n}Z_n\right]}
{\e^{(e_{m})}\left[e^{-\theta Z_{n-m}}Z_{n-m}\right]}
\frac{\e^{(e_m)}\left[e^{-\theta Z_{n-m}}Z_{n-m}^{[k]}\right]}
{\e^{(e)}\left[e^{-\theta Z_n}Z_n^{[k]}\right]}.
\end{split}
\end{equation*}
In particular, $\mathbf{Q}_n^{(e,2)}(\psi_1=m)$ has distribution \eqref{eq_law_Kn}.
 \end{proposition}
\begin{rem}
Suppose that $Z_n$ has finite $k$-th moments for any $n\geq 0$ and we  take $\theta=0$. By  \eqref{eq_media_Zn},  \eqref{factorial moments} and Proposition  \ref{prop:measure Q}, the previous probability can be interpreted as the ways we select one particle at generation $m$ that has all the marks and from such particle, we consider the subtree with $k$ spines,  that is  
$$\mathbf{Q}^{(e,k,\theta)}_n(\psi_1\geq m)
=\frac{\e^{(e)}\left[Z_m\right] \e^{(e_m,k)}_{n-m}\left[h_{n-m}(\textbf{T};\textbf{V}_1,\ldots,\textbf{V}_k)\right]}
{\e^{(e,k)}_n\left[h_{n}(\textbf{T};\textbf{V}_1,\ldots,\textbf{V}_k)\right]}.$$
\end{rem}

\begin{proof} 
Observe that $\{\psi_1\geq m\}=\{R_m(\textbf{v}_1)=R_m(\textbf{v}_i),  i\leq k\}$ is the event that all $k$ spines are still together at time $m$. From Proposition \ref{prop:measure Q}, we have 
\begin{equation*}
\begin{split}
\mathbf{Q}_n^{(e,k,\theta)}(\psi_1\geq m)
= \sum
\frac{e^{-\theta X_n(\textbf{t})}
\indi_{\{R_m(\textbf{v}_1)=R_m(\textbf{v}_i),  i\leq k\}}
\indi_{\{\textbf{v}_i\neq \textbf{v}_j, i\neq j\}}
\indi_{\{|\textbf{v}_i|= n, i\leq k\}}
\mathbf{G}^{(e)}_n(\textbf{t})}{\e^{(e)}\left[e^{-\theta Z_n}Z_n^{[k]}\right]},
\end{split}
\end{equation*}
where the previous sum is taken over all $(\textbf{t};\textbf{v}_1,\dots, \textbf{v}_k)\in\mathcal{T}^k_n$.

Let $(\textbf{t};\textbf{v}_1,\dots, \textbf{v}_k)\in\mathcal{T}^k_n$ such that 
\[
\indi_{\{R_m(\textbf{v}_1)=R_m(\textbf{v}_i), i\leq k\}}\indi_{\{\textbf{v}_i\neq \textbf{v}_j, i\neq j\}}\indi_{\{|\textbf{v}_i|= n, i\leq k\}}\neq 0.
\]
We also let $u^*$ be the  ancestor of $\textbf{v}_1$ at generation $m$. Then, $|M_{u^*}(\textbf{t})|=k$ .
We enumerate by  $u_1,\dots ,u_{X_m(\textbf{t})-1}$  the particles  at generation $m$ which are different of $u^*$. As in the proof of the previous proposition,  we decompose $(\textbf{t};\textbf{v}_1,\ldots,\textbf{v}_k)$ as the concatenation of its restriction to generation $m$  with all the subtrees which are attached to  particles $u^*, u_1,\dots ,u_{X_m(\textbf{t})-1}$ at generation $m$ (see for instance Figure \ref{fig:4}), that is
\begin{equation}\label{subtrees}
(\textbf{t};\textbf{v}_1,\ldots,\textbf{v}_k)=R_m (\textbf{t};\textbf{v}_1,\ldots,\textbf{v}_k)	\underset{r=1}{\overset{X_m(\textbf{t})-1}{\bigsqcup}} \hspace{-.1cm}u_r S_{u_r}(\textbf{t};\textbf{v}_1,\ldots,\textbf{v}_k)\bigsqcup u^* S_{u^*}(\textbf{t};\textbf{v}_1,\ldots,\textbf{v}_k).
\end{equation}
Since $|M_{u^*}(\textbf{t})|=k$, we  have  $S_{u^*}(\textbf{t};\textbf{v}_1,\ldots,\textbf{v}_k)=(\textbf{t}^\ast;\textbf{v}_1^\ast,\dots, \textbf{v}_k^\ast)\in \mathcal{T}_{n-m}^k$. Moreover, we introduce  $\textbf{t}_r:=S_{u_r}(\textbf{t};\textbf{v}_1,\ldots,\textbf{v}_k)$ which is in $\mathcal{T}_{n-m}$ from our assumption.  Additionally, we have
$R_m(\textbf{t};\textbf{v}_1,\ldots,\textbf{v}_k)=(\textbf{s};\textbf{w}_1,\dots, \textbf{w}_k) \in \mathcal{T}_m^k$  and satisfies 
$X_m(\textbf{t})=X_m(\textbf{s})$. Thus,  we can decompose $X_n(\textbf{t})$ as follows
\begin{equation}\label{subpopulations}
X_n(\textbf{t})=X_{n-m}(\textbf{t}^\ast)+\sum_{j=1}^{X_m(\textbf{t})-1}X_{n-m}(\textbf{t}_j).
\end{equation}
Then, by taking into account the environment and the previous decompositions, we have 
\begin{equation*}
\begin{split}
&\underset{(\textbf{t};\textbf{v}_1,\dots, \textbf{v}_k)\in\mathcal{T}^k_n}{\sum}
e^{-\theta X_n(\textbf{t})}\indi_{\{R_m(\textbf{v}_1)=R_m(\textbf{v}_i),  i\leq k\}}\indi_{\{\textbf{v}_i\neq \textbf{v}_j, i\neq j\}}\indi_{\{|\textbf{v}_i|= n, i\leq k\}}\mathbf{G}^{(e)}_n(\textbf{t})\\
&=\underset{(\textbf{s};\textbf{w}_1,\dots, \textbf{w}_k) \in\mathcal{T}^k_m}{\sum}\indi_{\{\textbf{w}_1=\textbf{w}_i, i\leq k\}}
\indi_{\{|\textbf{w}_1|= m\}}
\prod_{r=1}^{X_m(\textbf{s})-1}
\left(\underset{\textbf{t}_r \in\mathcal{T}_{n-m}}{\sum}
e^{-\theta X_{n-m}(\textbf{t}_r)} \mathbf{G}_{n-m}^{(e_m)}(\textbf{t}_r)\right)\mathbf{G}_m^{(e)}(\textbf{s})\\
&\hspace{1.5cm} \times
\left(\underset{(\textbf{t}^\ast;\textbf{v}_1^\ast,\dots, \textbf{v}_k^\ast)\in \mathcal{T}_{n-m}^k}{\sum}e^{-\theta X_{n-m}(\textbf{t}^\ast)}
\indi_{\{\textbf{v}_i^\ast\neq \textbf{v}_j^\ast, i\neq j\}}\indi_{\{|\textbf{v}_i^\ast|= n-m, i\leq k\}}\mathbf{G}_{n-m}^{(e_m)}(\textbf{t}^\ast)\right).
\end{split}
\end{equation*}
By equation \eqref{sum with spines}, we deduce
$$\underset{(\textbf{t}^\ast;\textbf{v}_1^\ast,\dots, \textbf{v}_k^\ast)\in \mathcal{T}_{n-m}^k}{\sum}
\hspace{-.4cm}e^{-\theta X_{n-m}(\textbf{t}^\ast)}
\indi_{\{\textbf{v}_i^\ast\neq \textbf{v}_j^\ast, i\neq j\}}\indi_{\{|\textbf{v}_i^\ast|= n-m, i\leq k\}}\mathbf{G}_{n-m}^{(e_m)}(\textbf{t}^\ast)=\e^{(e_m)}\left[e^{-\theta Z_{n-m}}Z_{n-m}^{[k]}\right].$$
On the other hand, for each $r\leq X_{m}(\textbf{s})-1$, we have 
\begin{equation}\label{subtree j}
\underset{\textbf{t}_r \in\mathcal{T}_{n-m}}{\sum}
e^{-\theta X_{n-m}(\textbf{t}_r)}\mathbf{G}_{n-m}^{(e_m)}(\textbf{t}_r)=\e^{(e_m)}\left[e^{-\theta Z_{n-m}}\right]=f_{m,n}(e^{-\theta}),
\end{equation}
and observe that for every tree $\textbf{s}$, there are $X_m(\textbf{s})$ possibilities to choose $\textbf{w}_1=\cdots=\textbf{w}_k$. Therefore,
\begin{equation*}
\begin{split}
\underset{(\textbf{s};\textbf{w}_1,\dots, \textbf{w}_k) \in\mathcal{T}^k_m}{\sum}&\indi_{\{\textbf{w}_1=\textbf{w}_i, i\leq k\}}
\indi_{\{|\textbf{w}_1|= m\}}
\prod_{r=1}^{X_m(\textbf{s})-1}
\left(\underset{\textbf{t}_r \in\mathcal{T}_{n-m}}{\sum}
e^{-\theta X_{n-m}(\textbf{t}_r)}\mathbf{G}_{n-m}^{(e_m)}(\textbf{t}_r)\right)
\mathbf{G}_m^{(e)}(\textbf{s})\\
&\hspace{2cm}=\underset{\textbf{s} \in\mathcal{T}_m}{\sum}X_m(\textbf{s}) 
\left( f_{m,n}(e^{-\theta})\right)^{X_m(\textbf{s})-1}
\mathbf{G}_m^{(e)}(\textbf{s}).
\end{split}
\end{equation*}
Putting all  pieces together, we obtain
$$\mathbf{Q}^{(e,k,\theta)}_n(\psi_1\geq m)
=\e^{(e)}\left[Z_m (f_{m,n}(e^{-\theta}))^{Z_m-1}\right]\frac{\e^{(e_m)}\left[e^{-\theta Z_{n-m}}Z_{n-m}^{[k]}\right]}{\e^{(e)}\left[e^{-\theta Z_n}Z_n^{[k]}\right]}.$$
	
In order to get our result,  we differentiate  \eqref{eq: laplace} and  then deduce
\begin{equation}\label{eq Z_mf_{m,n}}
\begin{split}
\e^{(e)}[Z_{m} (f_{m,n}(e^{-\theta}))^{Z_{m}-1}]&=f_{0,m}'(f_{m,n}(e^{-\theta}))=\frac{f_{0,n}'(e^{-\theta})}{f_{m,n}'(e^{-\theta})}=\frac{\e^{(e)}\left[e^{-\theta Z_n}Z_n\right]}{\e^{(e_{m})}\left[e^{-\theta Z_{n-m}}Z_{n-m}\right]},
\end{split}
\end{equation}
where we have use that $f_{0,n}$ is a composition of functions and  the chain rule for its derivative. The latter equality implies the first identity in the statement.

Finally, we compute the case when $k=2$. Recall that in this case, we can take $\theta=0$. Hence from  \eqref{eq_media_Zn} and the shifted environment property \eqref{moved environment}, we see
\begin{equation*}
\begin{split}
\mathbf{Q}^{(e,2)}_n(\psi_1\geq m)
&=\e^{(e)}[Z_m]\frac{\e^{(e_m)}[Z_{n-m}(Z_{n-m}-1)]}{\e^{(e)}[Z_n(Z_n-1)]}\\
&=\mu_m\frac{\left(\frac{\mu_n}{\mu_m}\right)^2\underset{j=m}{\overset{n-1}{\sum}}\frac{\nu_{j+1}}{\mu_j}\mu_m}{\mu_n^2\underset{j=0}{\overset{n-1}{\sum}}\frac{\nu_{j+1}}{\mu_j}}=\frac{\underset{j=m}{\overset{n-1}{\sum}}\frac{\nu_{j+1}}{\mu_j}}{\underset{j=0}{\overset{n-1}{\sum}}\frac{\nu_{j+1}}{\mu_j}},
\end{split}
\end{equation*}
which implies that  $\mathbf{Q}_n^{(e,2)}(\psi_1=m)$ has distribution \eqref{eq_law_Kn}.  This finishes the proof.
\end{proof}
It is important to note that by an analogous  technique, we can obtain the offspring distribution  of particles off spines under $\mathbf{Q}_n^{(e,k,\theta)}$. More precisely, let $u$ be a particle without marks such that $|u|=m$. From Proposition \ref{prop:branching property}, we see that under $\mathbf{Q}_n^{(e,k,\theta)}$, the subtree attached to $u$ has law $\mathbf{Q}_{n-m}^{(e_m,0,\theta)}$. Similarly  to the proof of the previous proposition, we  perform a decomposition by using all the subtrees attached to the offspring of $u$ and using equation \eqref{subtree j}, we deduce the following corollary.
\begin{corollary}\label{offspring without}
Under $\mathbf{Q}_n^{(e,k,\theta)}$, an unmarked particle in generation $m\in \{0,\dots, n-1\}$ gives birth to unmarked particles in generation $m+1$, independently of other particles, with probability  
\begin{equation*}
q_{m+1}(\ell) \frac{f_{m+1,n}(e^{-\theta})^\ell}{\e^{(e_m)}\left[e^{-\theta Z_{n-m}}\right]}, \qquad \ell\geq 0.
\end{equation*}
\end{corollary}

From the change of measure \eqref{cmg}, we notice that under  $\mathbf{Q}^{(e,k,\theta)}_n$,  there are $k$ different spines at generation $n$. The latter can also be deduced from the forward construction.   Indeed by a recursive argument, it is enough to show that if $\psi_1=n-1$, then $g=k$. Thus, by using Proposition \ref{prop: first splitting time} and \eqref{eq: g groups}, we may deduce 
$$ \mathbf{Q}^{(e,k,\theta)}_n\left(g=k \mid \psi_1=n-1  \right)=  \frac{e^{-k\theta}\left.\tfrac{\partial^k}{\partial s^k }f_{n}(s)\right|_{s=f_{n,n}(e^{-\theta})}}{ \e^{(e_{n-1})}\left[e^{-\theta Z_1}Z_{1}^{[k]}\right]}.$$
Since $f_{n,n}(e^{-\theta})=e^{-\theta}$, the claim then follows from \eqref{factorial moments}.

Now, we turn our attention into the  offspring distribution of particles with marks under $\mathbf{Q}_n^{(e,k,\theta)}$. Suppose that at generation $m$ there is only one particle with  $k$ marks, i.e. 
$R_m(\textbf{V}_1)=R_m(\textbf{V}_i)$ for each $i\leq k$. When this particle has offspring, the marks follows one or several of its  offspring. Recall that $g\in\{1,\dots,k\}$  is the number of spine groups created, $k_1,\ldots,k_g\geq 1$ are their sizes and $g_r$ is the number of groups of size $r$. Moreover, the relations \eqref{eq: combinatorial} hold.
Observe that   under the event 
$\{\psi_1=m\}$, we necessarily have $g\geq 2$. 
Also, recall that $\hat\psi_1$ is the first time the spines split apart, $\psi_1$ is the last time all the spines are together, and $\hat\psi_1=\psi_1+1$. We often swap between using $\psi_1$ and $\hat\psi_1$ as most appropriate.

\begin{proposition}\label{prop: offspring spine} 
	
The probability that at generation $m$ the $k$ spines are still all together but at the next generation the spines have split into $g\in\{2,\dots,k\}$ groups of sizes $k_1,\ldots,k_g\geq 1$ is given by \eqref{eq: g groups}, i.e.
\begin{equation}\label{tomatito}
\begin{split}
&\mathbf{Q}^{(e,k,\theta)}_n\left(\psi_1=m, C_{k;k_1,\dots, k_g}\right)\\
&=\frac{k!
\Big(\left.\tfrac{\partial^g}{\partial s^g }f_{m+1}(s)\right|_{s=	f_{m+1,n}(e^{-\theta})}\Big)
\e^{(e)}\left[e^{-\theta Z_n}Z_n\right]
\underset{i=1}{\overset{g}{\prod}}\e^{(e_{m+1})}\left[e^{-\theta Z_{n-m-1}}Z_{n-m-1}^{[k_i]}\right]}
{\left(\underset{r=1}{\overset{k}{\prod}}g_r!\right) \left( \underset{i=1}{\overset{g}{\prod}}k_i!\right)
\e^{(e_{m})}\left[e^{-\theta Z_{n-m}}Z_{n-m}\right]
\e^{(e)}\left[e^{-\theta Z_n}Z_{n}^{[k]}\right]},
\end{split}
\end{equation}
where $C_{k;k_1,\dots, k_g}$ is defined in \eqref{def: C}.

Moreover, we may extend  \eqref{eq: g groups} to include  the case when all $k$ spines remain together at generation $m+1$  (corresponding to $g=1$ and $k_1=k$) and additionally specify the offspring distribution of $V^{(m)}_1$, the particle at generation $m$ with $k$ marks. Indeed, suppose that the spines have not split by time $m$ and at generation $m+1$ the spines follow $g\in\{1,\dots,k\}$ groups of sizes $k_1,\ldots,k_g\geq 1$, then, for $\ell\geq g$, we have
\begin{equation}\label{tomate}
\begin{split}
&\mathbf{Q}^{(e,k,\theta)}_n\left( \hat\psi_1>m;\,\, l_{V^{(m)}_1}=\ell; \mbox{ at time $m+1$ the spines follow }g \mbox{  groups of sizes } k_1,\dots, k_g\right)\\
&=
\frac{q_{m+1}(\ell)\ell!k!\left(f_{m+1,n}(e^{-\theta})\right)^{\ell-g} }{ (\ell-g)!\left(\underset{r=1}{\overset{k}{\prod}}g_r!\right) \left( \underset{i=1}{\overset{g}{\prod}}k_i!\right)}
\frac{\e^{(e)}\left[e^{-\theta Z_n}Z_n\right] \underset{i=1}{\overset{g}{\prod}}\e^{(e_{m+1})}\left[e^{-\theta Z_{n-m-1}}Z_{n-m-1}^{[k_i]}\right]}
{\e^{(e_{m})}\left[e^{-\theta Z_{n-m}}Z_{n-m}\right]\e^{(e)}\left[e^{-\theta Z_n}Z_{n}^{[k]}\right]	 }.
\end{split}
\end{equation}
In particular, if $g=1$ we find that, given all the spines stay together, the offspring distribution is (once) size biased and discounted according to population size over the remaining time, as follows:
\begin{equation}\label{eq: 1 group}
\begin{split}
&\mathbf{Q}^{(e,k,\theta)}_n\left(l_{V^{(m)}_1}=\ell \mid \psi_1>m \right)=
\frac{\ell q_{m+1}(\ell) \left(f_{m+1,n}(e^{-\theta})\right)^{\ell-1} }
{f_{m+1}'(f_{m+1,n}(e^{-\theta}))}.
\end{split}
\end{equation}

\end{proposition}

\begin{proof}
We first prove identity \eqref{tomate}. In order to do so,  we use the same tree decomposition as in the proof Proposition \ref{prop: first splitting time}. Indeed, under the event 
\[
\{m\le \psi_1\}\cap \{ l_{V^{(m)}_1}=\ell\}\cap\{\mbox{the spines follow }g \mbox{ groups of sizes }  k_1,\dots, k_g\},
\]
we  decompose $(\textbf{t};\textbf{v}_1,\dots, \textbf{v}_k)\in\mathcal{T}^k_n$ into (see  \eqref{subtrees}), 
\begin{equation*}
(\textbf{t};\textbf{v}_1,\ldots,\textbf{v}_k)=(\textbf{s};\textbf{w}_1,\dots, \textbf{w}_k)
\underset{r=1}{\overset{X_m(\textbf{t})-1}{\bigsqcup}} \hspace{-.1cm}u_r \textbf{t}_r\bigsqcup u^*(\textbf{t}^\ast;\textbf{v}_1^\ast,\dots, \textbf{v}_k^\ast)
\end{equation*}
where $u^*$ is the  ancestor of $\textbf{v}_1$ at generation $m$,  the other particles at generation $m$ are labelled by $u_1,\dots, u_{X_m(\textbf{t})-1}$; the subtrees  $(\textbf{s};\textbf{w}_1,\dots, \textbf{w}_k) =R_m(\textbf{t};\textbf{v}_1,\ldots,\textbf{v}_k)\in \mathcal{T}_{m}^k$, $\textbf{t}_r=S_{u_r}(\textbf{t};\textbf{v}_1,\ldots,\textbf{v}_k)\in \mathcal{T}_{n-m}$ and  $(\textbf{t}^\ast;\textbf{v}_1^\ast,\dots, \textbf{v}_k^\ast)=S_{u^*}(\textbf{t};\textbf{v}_1,\ldots,\textbf{v}_k)\in \mathcal{T}_{n-m}^k$. 
	
Next, we observe that  a similar  decomposition for $(\textbf{t}^\ast;\textbf{v}_1^\ast,\dots, \textbf{v}_k^\ast)$ can be given at generation $m+1$.  By hypothesis, $u^{*}$ has $\ell-g$ offspring without marks and $g$ with marks. We label  by $ y_1,\dots, y_{\ell-g}$ the offspring without marks and by $y_1^*,\dots, y_{g}^*$  those with $k_1, \dots, k_g$ marks, respectively. We denote by $\textbf{r}_i:=S_{y_i}(\textbf{t}^\ast;\textbf{v}_1^\ast,\dots, \textbf{v}_k^\ast) \in \mathcal{T}_{n-m-1}$, for $i\leq \ell -g $, and $(\textbf{t}^{\ast,j};\textbf{v}_1^{\ast,j},\dots, \textbf{v}_{k_j}^{\ast,j}):=S_{y_j^*}(\textbf{t}^\ast;\textbf{v}_1^\ast,\dots, \textbf{v}_k^\ast) \in \mathcal{T}_{n-m-1}^{k_j}$ for all $j\leq g$. 
Then, the following decomposition holds
\begin{equation}
\label{decomposition2}
(\textbf{t}^\ast;\textbf{v}_1^\ast,\dots, \textbf{v}_k^\ast)=(\textbf{s}^{*};\textbf{w}_1^{*},\dots, \textbf{w}_k^{*})\ \underset{i=1}{\overset{\ell-g}{\bigsqcup}}\  y_i\textbf{r}_i \ \underset{j=1}{\overset{g}{\bigsqcup}}\ 
y_j^* (\textbf{t}^{\ast,j};\textbf{v}_1^{\ast,j},\dots, \textbf{v}_{k_i}^{\ast,j}),	\end{equation}
where $(\textbf{s}^{*};\textbf{w}_1^{*},\dots, \textbf{w}_k^{*}):=R_1(\textbf{t}^\ast;\textbf{v}_1^\ast,\dots, \textbf{v}_k^\ast) \in \mathcal{T}_1^k$ and has $\ell$ leaves where $\ell-g$ has no marks  and  the remaining $g$ leaves have $k_1,\dots, k_g$ marks. 
Moreover, we can decompose $X_{n-m}(\textbf{t}^\ast)$ as  follows
\begin{equation}\label{subpopulations2}
X_{n-m}(\textbf{t}^\ast)=\sum_{i=1}^{\ell -g }X_{n-m-1}(\textbf{r}_i)+\sum_{j=1}^{g }X_{n-m-1}(\textbf{t}^{\ast,j}).
\end{equation}
For simplicity on exposition, we let 
\[
A:=\mathbf{Q}^{(e,k,\theta)}_n\left(\hat\psi_1>m;\,\, l_{V^{(m)}_1}=\ell; \mbox{ at  $m+1$ the spines follow }g \mbox{  groups of sizes } k_1,\dots, k_g\right).
\] 
We proceed similarly as  in the proof of the previous proposition, that is by Proposition \ref{prop:measure Q} and taking into account the environment in decompositions \eqref{subtrees}, \eqref{subpopulations},  \eqref{decomposition2} and \eqref{subpopulations2}, we deduce 
\begin{equation*}
\begin{split}
&A\e^{(e)}\left[e^{-\theta Z_n}Z_n^{[k]}\right]\\
& =\underset{(\textbf{s};\textbf{w}_1,\dots, \textbf{w}_k) \in\mathcal{T}^k_m}{\sum}
\indi_{\{\textbf{w}_1=\textbf{w}_i, i\leq k\}}
\indi_{\{|\textbf{w}_1|= m\}}
\prod_{r=1}^{X_m(\textbf{s})-1}
\left(\underset{\textbf{t}_r \in\mathcal{T}_{n-m}}{\sum}
e^{-\theta X_{n-m}(\textbf{t}_r)}\mathbf{G}_{n-m}^{(e_m)}(\textbf{t}_r)\right)
\mathbf{G}_m^{(e)}(\textbf{s})\\ 
& \times q_{m+1}(\ell)\indi_{\{\mbox{g children has } k_1,\dots,k_g \mbox{ marks}\}}
\prod_{i=1}^{\ell-g}
\left(\underset{\textbf{r}_i\in \mathcal{T}_{n-m-1}}{\sum}
e^{-\theta X_{n-m-1}(\textbf{r}_i)}\mathbf{G}_{n-m-1}^{(e_{m+1})}(\textbf{r}_i)\right)\\
&\times \underset{j=1}{\overset{g}{\prod}}\,\,\sum\indi_{\{\textbf{v}_l^{\ast,j}\neq \textbf{v}_s^{\ast,j}, l\neq s\}}\indi_{\{|\textbf{v}_l^{\ast,j}|= n-m-1, l\leq k_j\}}e^{-\theta X_{n-m-1}(\textbf{t}^{\ast,j})}\mathbf{G}_{n-m-1}^{(e_{m+1})}(\textbf{t}^{\ast,j}), 
\end{split}
\end{equation*}
where the sum runs over all $(\textbf{t}^{\ast,j};\textbf{v}_1^{\ast,j},\dots, \textbf{v}_{k_j}^{\ast,j})\in \mathcal{T}_{n-m-1}^{\, k_j}$.

Observe that there are $X_m(\textbf{s})$ possibilities to choose $\textbf{w}_1$, $u^*$ has $\ell$ offspring with probability $q_{m+1}(\ell)$, from these offspring there are ${\ell \choose g}$ possibilities to select $g$ of them with marks, there are $g!/\underset{r=1}{\overset{k}{\prod}}g_r!$ ways of arranging the values $k_1,\ldots,k_g$ with the sizes of each group and $k!/ \underset{i=1}{\overset{g}{\prod}}k_i!$ ways to arrange the $k$ spines into $g$ groups of sizes $k_1,\ldots,k_g$. Thus using \eqref{subtree j} for every $r\leq X_m(\textbf{s})-1$ and every $i\leq \ell -g$;  and \eqref{sum with spines} for every $j\leq g$,  we obtain that $A$ is equal to 
$$
\frac{\e^{(e)}[Z_{m} (f_{m,n}(e^{-\theta}))^{Z_{m}-1}]	q_{m+1}(\ell)\ell!k!\left(f_{m+1,n}(e^{-\theta})\right)^{\ell-g}}
{ (\ell-g)!\left(\underset{r=1}{\overset{k}{\prod}}g_r!\right) \left( \underset{i=1}{\overset{g}{\prod}}k_i!\right)}
\frac{ \underset{i=1}{\overset{g}{\prod}}\e^{(e_{m+1})}\left[e^{-\theta Z_{n-m-1}}Z_{n-m-1}^{[k_i]}\right]}
{	\e^{(e)}\left[e^{-\theta Z_n}Z_{n}^{[k]}\right]
},$$
which implies \eqref{tomate} after using \eqref{eq Z_mf_{m,n}}. 
	
For identity \eqref{eq: 1 group}, we take $g=1$ in \eqref{tomate} and use Proposition \ref{prop: first splitting time} with $m+1$ instead of $m$  to obtain
\begin{equation*}
\begin{split}
&\mathbf{Q}^{(e,k,\theta)}_n\left(\psi_1>m,  l_{V^{(m)}_1}=\ell \right)\\
&=\ell q_{m+1}(\ell)\left(f_{m+1,n}(e^{-\theta})\right)^{\ell-1} 
\frac{ \e^{(e_{m+1})}\left[e^{-\theta Z_{n-m-1}}Z_{n-m-1}\right]}
{	\e^{(e_{m})}\left[e^{-\theta Z_{n-m}}Z_{n-m}\right]}
\mathbf{Q}^{(e,k,\theta)}_n(\psi_1\geq m+1).
\end{split}
\end{equation*}
By \eqref{eq Z_mf_{m,n}} we have
\begin{equation*}
\begin{split}
\frac{ \e^{(e_{m+1})}\left[e^{-\theta Z_{n-m-1}}Z_{n-m-1}\right]}
{	\e^{(e_{m})}\left[e^{-\theta Z_{n-m}}Z_{n-m}\right]}=\frac{f_{m+1,n}'(e^{-\theta})}{f_{m,n}'(e^{-\theta})}=\frac{1}{f_{m+1}'(f_{m+1,n}(e^{-\theta}))},
\end{split}
\end{equation*}
where we have use that $f_{m,n}$ is a composition of functions and  the chain rule for its derivative.  By observing that $\{\psi_1\geq m+1\}=\{\psi_1> m\}$ and performing  conditional probability, we get \eqref{eq: 1 group}.
	
Finally, we prove \eqref{eq: g groups}. So let $g\geq 2$.  In this case, $\{m\leq \psi_1\}=\{m=\psi_1\}$ and by summing over all possible values of $\ell$ and using the $g$-th derivative of the Laplace transform, we obtain \eqref{eq: g groups}. This completes the proof.
\end{proof}

We are now ready to present the proof of Proposition \ref{prop: forward construction}.
\begin{proof}[Proof of Proposition  \ref{prop: forward construction}]
As we mentioned before, the proof follows from previous intermediate results. More precisely, we start with a particle with $k$ marks. Step (2) follows from Proposition \ref{prop: offspring spine}, equation \eqref{tomate}. Corollary \ref{offspring without} provides the offspring distribution of a particle without marks (Step (3)). For any generation $m<\psi_1$, the offspring distribution of a marked particles is given in 
Proposition \ref{prop: offspring spine}, equation \eqref{eq: 1 group}. This is Step (4). Step (5) is nothing but the probability, under 
$\mathbf{Q}^{(e,k,\theta)}_n$, that the marked particle at generation $m$ gives $\ell$ offspring  given that $m$ is the last time that the spines stay together  and the spines splits into $g$ subgroups of sizes $k_1,\cdots, k_g$. That is 
$$\mathbf{Q}^{(e,k,\theta)}_n\left(l_{V^{(m)}_1}=\ell\mid \psi_1=m, C_{k;k_1,\dots, k_g}\right)=\frac{\mathbf{Q}^{(e,k,\theta)}_n\left(l_{V^{(m)}_1}=\ell \psi_1=m, C_{k;k_1,\dots, k_g}\right)}{\mathbf{Q}^{(e,k,\theta)}_n\left(\psi_1=m, C_{k;k_1,\dots, k_g}\right)}.$$
This follows from by \eqref{tomatito} and \eqref{tomate}. Finally, the last step follows from the Markov property for GWTVE with spines under $\mathbf{Q}_n^{(e,k,\theta)}$, that is Proposition \ref{prop: branching property Q}.
\end{proof}

To conclude this section, we compute the distribution of the total population  size of unmarked particles who are descendants of   the unmarked offspring of  spine particles. It  will be used in the next section, where we analyse the limiting behaviour of functionals of trees with $k$ spines, under $\mathbf{Q}_n^{(e, k, \theta)}$ as $n\rightarrow\infty$.
Let $u$ be a spine particle
and denote by $Y_n(u)$ the total population size of unmarked particles  at generation  $n$ that descend from the unmarked offspring of $u$, see Figure \ref{fig:5} for an illustrative example. In other words, if $u$ has $g$ marked children labelled by $uv_1, \dots,uv_g$,  then
$$Y_n(u):=X_{n-|u|}(S_u(\textbf{T}))-\sum_{j=1}^{g}X_{n-|uv_j|}(S_{uv_j}(\textbf{T})).$$
By Proposition  \ref{prop: branching property Q}, the variables $Y_n(u)$, for $u$ in the spines,  are independent.

\begin{lemma}\label{lem: constructed block}
	Let $(\textbf{T}; \textbf{V}_1,\dots,\textbf{V}_k)$ be a tree with $k\ge 1$ spines, under $\mathbf{Q}_n^{(e,k,\theta)}$. For each $m\leq n-1$, let $u$ be a marked particle at generation $m$ whose offspring distribution is given by $q_{m+1}^{(g,\theta)}$, for $g\geq 1$. Then,
	\begin{equation*}
	\mathbf{Q}_n^{(e,k,\theta)}\left[e^{-(\lambda-\theta) Y_n(u)}\right]=\frac{\left.\tfrac{\partial^g}{\partial s^g }f_{m+1}(s)\right|_{s=f_{m+1,n}(e^{-\lambda})}}{\left.\tfrac{\partial^g}{\partial s^g }f_{m+1}(s)\right|_{s=f_{m+1,n}(e^{-\theta})}}, \qquad \lambda >0.
	\end{equation*}
	In particular, for $k=1$ we have 
	\begin{equation*}
	\mathbf{Q}_n^{(e,1,\theta)}
	\left[e^{-(\lambda-\theta) Z_n}\right]=\frac{e^{-\lambda}}{e^{-\theta}}\prod_{m=0}^{n-1}\frac{f_{m+1}'(f_{m+1,n}(e^{-\lambda}))}{f_{m+1}'(f_{m+1,n}(e^{-\theta}))}, \qquad \lambda >0.
	\end{equation*}
\end{lemma}

\begin{proof}
	From Proposition \ref{prop: branching property Q} and the definition of $\mathbf{Q}^{(e,0,\theta)}$, the descendants at generation $n$ of an unmarked particle at generation $m+1$ is given by
	$$\mathbf{Q}_{n-m-1}^{(e_{m+1},0,\theta)}\left[e^{-(\lambda-\theta) Z_{n-m-1}}\right]=\frac{\e^{(e_{m+1})}\left[e^{-\lambda Z_{n-{m+1}}}\right]}{\e^{(e_{m+1})}\left[e^{-\theta Z_{n-{m+1}}}\right]}=\frac{f_{m+1,n}(e^{-\lambda})}{f_{m+1,n}(e^{-\theta})}, \qquad \lambda>0.$$
	If the particle $u$ has offspring distribution  $q_{m+1}^{g,\theta}(\cdot)$, we know  that it  gives birth at least to $g$ marked particles. Thus, by only considering all  subtrees attached to unmarked  offspring (which are independent),  for $\lambda>0$, we get
	\begin{equation*}
	\begin{split}
	\mathbf{Q}_n^{(e,k,\theta)}\left[e^{-(\lambda-\theta) Y_n(u)}\right]
	&=\underset{\ell=g}{\overset{\infty}{\sum}}
	\frac{\ell^{[g]}q_{m+1}(\ell)(f_{m+1,n}(e^{-\theta}))^{l-g}}
	{\left.\tfrac{\partial^g}{\partial s^g }f_{m+1}(s)\right|_{s=f_{m+1,n}(e^{-\theta})}}
	\frac{(f_{m+1,n}(e^{-\lambda}))^{l-g}}{(f_{m+1,n}(e^{-\theta}))^{l-g}}\\
	&=\frac{\left.\tfrac{\partial^g}{\partial s^g }f_{m+1}(s)\right|_{s=f_{m+1,n}(e^{-\lambda})}}{\left.\tfrac{\partial^g}{\partial s^g }f_{m+1}(s)\right|_{s=f_{m+1,n}(e^{-\theta})}}.
	\end{split}
	\end{equation*}
	The second equality holds by  independence and the decomposition 
	$$Z_n=1+\sum_{m=0}^{n-1}Y_n(v^{(m)}),$$
	where $v^{(m)}$ is the particle at the spine with height $m$.
\end{proof}

\section{Scaling limits of critical GWVE trees}
The  proof of our main result (Theorem \ref{principal}) requires several preliminary results that are presented in this section. Here, we suppose that conditions  \eqref{eq_cond_kersting} and \eqref{def: critical}  are satisfied.
Most of the preliminary results are limiting behaviours  of functionals of trees with $k$ spines, under $\mathbf{Q}_n^{(e, k, \theta)}$, and of the total population size under the original probability measure but  conditional on survival.  In order to do so, we first introduce some notation. Recall that $\hat\psi_1$ denotes the first time that a second spine is created. Hence for $1\leq i\leq k-1$, we introduce recursively $\psi_i$, the  last time where there are at most $i$ marked particles, as follows
$$\psi_i=\max\{m\geq 0: \mbox{ at generation }m \mbox{ there are at most } i \mbox{ marked particles}\}.$$
Additionally, we denote by $\hat\psi_i=\psi_i+1$, the $i$-th spine split time (although note these split times may coincide, for example, the first and second spine split times may be equal which corresponds to marks following three different particles).

Recall that $\rho_n$ is the normalized second factorial moment of $Z_n$ and $\{\tau_n(t), t\in[0,1], \}$ is its right-continuous inverse, see \eqref{rhonDef} and \eqref{SntDef}. We introduce  
\begin{equation}\label{rho and A}
A_{n,m}:=1-\frac{\rho_{m+1}}{\rho_n}=1 - \sum_{j=0}^{m} \frac{\nu_{j+1}}{\mu _{j}} \left(\sum_{k=0}^{n-1} \frac{\nu_{k+1}}{\mu _{k}}\right)^{-1}, \qquad 0\leq m< n. 
\end{equation}
The term $A_{n,m}$ appears in the two-spine decomposition as the correct factor that for the right normalisation of  the descendants attached to the shorter spine, see   \cite[Section 2]{CardonaPalau}.
Note that   $\{A_{n,m}: m= 0, \dots, n-1\}\subset[0, 1]$  is a decreasing sequence with  $A_{n,n-1}=0$.  If we define $A_{n,-1}:=1$, we can introduce a partition of $[0,1]$ as follows
\[
R^{(n)}= \{0= A_{n,n-1} < A_{n,n-2} < \ldots < A_{n,0}<A_{n,-1}\}.
\] 
Note that the norm of the partition $(\rho_m/\rho_n; 0\leq m\leq n)$  is equal to the norm of the partition $R^{(n)}$.
According to  \cite[Lemma 1]{kersting2021genealogical}, under hypothesis \eqref{eq_cond_kersting} and \eqref{def: critical}, the norm of the previous partitions satisfies
\begin{equation}
\label{limit norm Rn}
\underset{n\rightarrow \infty}{\lim}|| R^{(n)}||=\underset{n\rightarrow \infty}{\lim}\max_{0\leq  m \leq n-1} \left\{\frac{\nu_{m+1}}{\mu _{m}} \left(\sum_{k=0}^{n-1} \frac{\nu_{k+1}}{\mu _{k}}\right)^{-1}\right\}=0.
\end{equation}
 Then, the following convergence is uniform on $0\leq t\leq 1$
\begin{equation}\label{limit rho}
1-A_{n,\tau_n(t)-1}=\frac{\rho_{\tau_n(t)}}{\rho_n}\longrightarrow t,\qquad \text{ as } n\rightarrow \infty.
\end{equation}

Recall that  $(f_1,f_2,\dots)$ denote  the generating functions associated to the environment $e=(q_1,q_2,\dots).$ For every $f_k$, we define its {\it shape function}, $\varphi_k:[0,1)\rightarrow \infty$, as 
$$\varphi_k(s):=\frac{1}{1-f_k(s)}-\frac{1}{f_k'(1)(1-s)},\qquad 0\leq s<1.$$
The function $\varphi_k$ can be extend continuously to $s=1$, by writing $\varphi_k(1)=\nu_k/2$.
Recall that $f_{m,n}=f_{m+1} \circ \ldots \circ f_n$. According to \cite[Lemma 5]{kersting2020unifying},  for every $0\leq m<n$,
\begin{equation}\label{lemma 5}
\frac{1}{1-f_{m,n}(s)}=\frac{\mu_m}{\mu_n(1-s)}+\mu_m\underset{k=m+1}{\overset{n}{\sum}}\frac{\varphi_k(f_{k,n}(s))}{\mu_{k-1}}, \qquad 0\leq s<1.
\end{equation}
In addition, by performing a small modification in the proof  of Lemma 8 in \cite{kersting2020unifying}, under hypothesis \eqref{eq_cond_kersting} and \eqref{def: critical}, we have  
\begin{equation}\label{lemma 8}
\underset{0\leq s\leq 1}{\sup}\left|\underset{k=m+1}{\overset{n}{\sum}}\frac{\varphi_k(f_{k,n}(s))}{\mu_{k-1}}-\frac{\rho_n-\rho_m}{2}\right|= o\left(\frac{\rho_n}{2}\right),
\end{equation}
uniformly for all $m<n$. For sake of completeness, we provide the proof of \eqref{lemma 8} in Lemma \ref{lemma:extensionkersting8} in the Appendix.  Kersting also showed in \cite[equation (3.4)]{kersting2021genealogical} that  for any $0<T<1$,
\begin{equation}\label{equation (5)}
\underset{k=m+1}{\overset{n}{\sum}}\frac{\varphi_k(f_{k,n}(0))}{\mu_{k-1}}=\frac{\rho_n-\rho_m}{2}(1+o(1)),
\end{equation}
uniformly for all $m$ such that $\rho_m\leq T \rho_n$. The latter could be deduced from \eqref{lemma 8} by taking $s=0$, since $\rho_n=O(\rho_n-\rho_m)$ for all $\rho_m \leq T \rho_n$. The previous equation was Kersting's primary tool to prove 
\begin{equation}\label{eq: lemma 2}
\p^{(e)}(Z_{n} > 0\mid Z_m=1)= \frac{2+o(1)}{\mu_m(\rho_n-\rho_m)},
\end{equation}
as $n\rightarrow \infty$, uniformly for all $m$ such that $\rho_m\leq T \rho_n$, see \cite[Lemma 2]{kersting2021genealogical}.

With all these ingredients in hand,  we proceed to prove the following limiting result  which is an extension of Yaglom's limit \eqref{eq: Yaglom}. 
Later on, we are going to use it to analyse the asymptotic behaviour of $\{\psi_1,\dots,\psi_{k-1}\}$ under $\mathbf{Q}^{(e,k,\theta)}_{n}$, as $n$ goes to infinity.  
The latter can be deduced through  an  inductive argument. More precisely,  to obtain the limiting probability of the event $\{\psi_{i+1}\leq \tau_n(t)\}=\{\rho_{\psi_{i+1}}/\rho_n \leq t\}$, we will apply the Markov property and we will condition on the information up to time $\psi_i$. Then, we will use the fact that $\psi_i$ is a function that depends on $n$ such that  $\rho_{\psi_{i}}/\rho_n$, under $\mathbf{Q}^{(e,k,\theta)}_{n}$,  converges in distribution to a given random variable. For this reason, we present the following results for general functions $u_n$ such that $\rho_{u_n}/\rho_n$ converges.
\begin{proposition}\label{limit starting Sn}
Assume that conditions  \eqref{eq_cond_kersting}  and \eqref{def: critical} are satisfied. Let $0<T<1$ and $u_n:[0,1]\rightarrow \{0, 1,\dots, n\}$   be a sequence of functions such that $\rho_{u_n(t)}/\rho_n\rightarrow \varrho(t)$, uniformly,  for $0\leq t\leq T$, as $n\rightarrow \infty$, where $\varrho:[0,1]\rightarrow [0,1]$ is a (non-decreasing) function satisfying $\varrho(t)<1$, for all $t<1$. Then,
\begin{equation*}
\left(\left(a_{n-u_n(t)}^{(e_{u_n(t)})}\right)^{-1}Z_{n-u_n(t)}\ ;\  \p^{(e_{u_n(t)})}(\ \cdot\ | Z_{n-u_n(t)} > 0)\right) \stackrel{(d)}{\longrightarrow}  \left(\mathbf{e}; \p\right) , \qquad \mbox{ as }\ n\rightarrow \infty,
\end{equation*}
uniformly for every $t\leq T$,	 where  $\mathbf{e}$  is a standard exponential r.v., under $\mathbb{P}$, and ``$\stackrel{(d)}{\longrightarrow}$'' means convergence in distribution. 
\end{proposition}

\begin{proof}
Let $0<T<1$. We observe that in order to get our  result,  it is enough to prove that  for every $\lambda >0$, the following limit holds
\begin{equation*}\label{eq: limit starting Sn}
\underset{n\rightarrow \infty}{\lim}\,\,\underset{0\leq t\leq T}{\sup}\e^{(e_{u_n(t)})}\left[\left.1-e^{-\lambda \left(a_{n-u_n(t)}^{(e_{u_n(t)})}\right)^{-1}Z_{n-u_n(t)}}\right| Z_{n-u_n(t)}>0\right]=\frac{\lambda}{1+\lambda}.
\end{equation*}
In order to prove our claim, we let $ m<n$ such that $\rho_m<\varrho(T)\rho_n$. By \eqref{lemma 5}, we have
\begin{equation*}
\begin{split}
\e^{(e_{m})}&\left[\left.1-e^{-\lambda \left(a_{n-m}^{(e_{m})}\right)^{-1}Z_{n-m}}\right| Z_{n-m}>0\right]=
\frac{1-f_{m,n}\left(e^{-\lambda/a_{n-m}^{(e_{m})}}\right)}{1-f_{m,n}(0)}\\
&=\frac{\frac{1}{\mu_n}+\underset{k=m+1}{\overset{n}{\sum}}\frac{\varphi_k(f_{k,n}(0))}{\mu_{k-1}}}{
\frac{1}{\mu_n\left(1-e^{-\lambda/ a_{n-m}^{(e_{m})}}\right)}+\underset{k=m+1}{\overset{n}{\sum}}
\frac{\varphi_k\left(f_{k,n}\left(e^{-\lambda/ a_{n-m}^{(e_{m})}}\right)\right)}{\mu_{k-1}}}.
\end{split}
\end{equation*}
We observe that for every $m$ with  $\rho_m\leq \varrho(T)\rho_n$, the following inequality  is satisfied 
\begin{equation*}\label{a_{n,m}}
\frac{1}{a_{n-m}^{(e_{m})}}=\frac{1}{a_n^{(e)}}
\frac{\rho_n}{\rho_n-\rho_m}\leq \frac{1}{a_n^{(e)}}\frac{1}{1-\varrho(T)}.
\end{equation*}
Then, 
\begin{equation}
\label{exponential uniform}
\left(1-e^{-\lambda/ a_{n-m}^{(e_{m})}}\right)=\frac{\lambda}{a_n^{(e)}}
\frac{\rho_n}{\rho_n-\rho_m}(1+o(1)),
\end{equation}
uniformly for all $m$ such that $\rho_m\leq \varrho(T) \rho_n$. Recall that $a_n^{(e)}=\mu_n\rho_n/2$ and  by criticality $\mu_n^{-1}=o(\rho_n)$. Then using \eqref{lemma 8}, \eqref{equation (5)} and \eqref{exponential uniform}, we  get 
\begin{equation*}
\begin{split}
\e^{(e_{m})}&\left[\left.1-e^{\lambda \left(a_{n-m}^{(e_{m})}\right)^{-1}Z_{n-m}}\right| Z_{n-m}>0\right]\\
&=\left(o(1)+\frac{\rho_n-\rho_m}{2\rho_n}(1+o(1))\right)
\left(\frac{\rho_n-\rho_m}{2\rho_n\lambda}+\frac{\rho_n-\rho_m}{2\rho_n}+o(1)\right)^{-1}
\end{split}
\end{equation*}
uniformly for all $m$ such that $\rho_m\leq \varrho(T) \rho_n$. By hypothesis,   $\rho_{u_n(t)}/\rho_n\rightarrow \varrho(t)$ uniformly  for $0\leq t\leq T$. Then, we   conclude that  	uniformly for every $t\leq T$,
$$\e^{(e_{u_n(t)})}\left[\left.1-e^{-\lambda \left(a_{n-u_n(t)}^{(e_{u_n(t)})}\right)^{-1}Z_{n-u_n(t)}}\right| Z_{n-u_n(t)}>0\right]\underset{n\rightarrow\infty}{\longrightarrow}\frac{1-\varrho(t)}{\lambda^{-1}(1-\varrho(t))+(1-\varrho(t))}.$$
We notice that the limit is equal to $\frac{\lambda}{1+\lambda}$, as expected.
\end{proof}

For a given $\theta$, we define $\theta_n:=\theta/a_{n}^{(e)}$. To simplify the notation, for every $0\leq m< n$ and $k\in \N$, we introduce  
\begin{equation}\label{Fknm}
F_{k,n,m}:=	
\frac{\e^{(e_{m})}\left[e^{-\theta_n Z_{n-m}} Z_{n-m}^{[k]}\right](1+\theta A_{n,m-1})^{k+1}\mu_{m}}
{\left(a_{n}^{(e)}\right)^{k-1}k!(A_{n,m-1})^{k-1}\mu_{n}}.
\end{equation}
As an application, we have the following lemma.

\begin{lemma}\label{lem: Yaglom}
	Let $0<T<1$ and $k\in \mathbb{N}$.	Suppose that $Z$ is a regular critical GWVE  or equivalently that conditions  \eqref{eq_cond_kersting} and \eqref{def: critical} are  satisfied, and $u_n:[0,1]\rightarrow \{0, 1,\dots, n\}$  is  a sequence of functions such that $\rho_{u_n(t)}/\rho_n\rightarrow \varrho(t)$, uniformly,  for $0\leq t\leq T$, as $n\rightarrow \infty$, where $\varrho:[0,1]\rightarrow [0,1]$ is a (non-decreasing) function satisfying that $\varrho(t)<1$, for all $t<1$.	Then, $$\lim_{n\rightarrow \infty}F_{j,n,u_n(t)}=1,$$
	uniformly on $0\le  t\leq T$ and $j\leq k.$
\end{lemma}

\begin{proof}For a fixed $k$, we can rewrite $F_{k,n,m}$ as follows
	\begin{equation*}
	\e^{(e_{m})}\left[\left. e^{-\theta_n Z_{n-m}}
	\frac{Z_{n-m}^{[k]}}{\left(a_{n}^{(e)}\right)^{k}}\right| Z_{n-m}>0\right]
	\frac{(1+\theta A_{n,m-1})^{k+1}}{k!\  (A_{n,m-1})^{k} }
	\frac{a_{n}^{(e)}A_{n,m-1}
		\p^{(e_{m})}( Z_{n-m} > 0)\mu_m}{ \mu_n}.
	\end{equation*}
	From our hypothesis and \eqref{rho and A},  $A_{n,u_n(t)-1}$ converges to $1-\varrho(t)$, as $n\rightarrow\infty$, uniformly on $0\le t\le T$. Then, by the previous result and the continuous mapping Theorem \cite[Theorem 5.27]{Kallenberg}, we have
	\begin{equation*}
	\left(A_{n,u_n(t)-1}\left(a_{n-u_n(t)}^{(e_{u_n(t)})}\right)^{-1}Z_{n-u_n(t)}\ ;\  \p^{(e_{u_n(t)})}(\ \cdot\ | Z_{n-u_n(t)} > 0)\right) \stackrel{(d)}{\longrightarrow}  \left((1-\varrho(t))\mathbf{e}; \p\right).
	\end{equation*}
	where  $\mathbf{e}$  is a exponential r.v., under $\mathbb{P}$.		
Note that $A_{n,u_n(t)-1}\left(a_{n-u_n(t)}^{(e_{u_n(t)})}\right)^{-1}=\left(a_n^{(e)}\right)^{-1}$ converges to zero as $n$ goes to infinity.
Additionally, the previous convergence implies that $Z_{n-u_n(t)}$ conditioned on $\{Z_{n-u_n(t)}>0\}$ diverge to infinity, as $n\rightarrow\infty$, and in particular   $\p^{(e_{u_n(t)})}\left(Z_{n-u_n(t)}\geq k\mid Z_{n-u_n(t)}>0\right)\to1.$ Since $$\frac{(Z_{n-u_n(t)}-k)^{k}\indi_{\{Z_{n-u_n(t)}\geq k\}}}{\left(a_{n}^{(e)}\right)^{k}}\leq \frac{Z_{n-u_n(t)}^{[k]}}{\left(a_{n}^{(e)}\right)^{k}}\leq \frac{Z_{n-u_n(t)}^k}{\left(a_{n}^{(e)}\right)^{k}},$$
and  $f(x)=e^{-\theta x}x^{k}$ is a bounded function, we have 
	\begin{equation}
	\label{eq:limit size biased}
	\e^{(e_{u_n(t)})}\left[\left. e^{-\theta  Z_{n-u_n(t)}/a_{n}^{(e)}}
	\frac{Z_{n-u_n(t)}^{[k]}}{\left(a_{n}^{(e)}\right)^{k}}\right| Z_{n-u_n(t)}>0\right]\longrightarrow \e\left[e^{-\theta (1-\varrho(t)) \mathbf{e}}(1-\varrho(t))^{k}\mathbf{e}^{k}\right],
	\end{equation}
	uniformly on $0\leq t\leq T$. On the other hand, it is known that 
	$$\e\left[e^{-\theta (1-\varrho(t)) \mathbf{e}}(1-\varrho(t))^{k}\mathbf{e}^{k}\right]=\frac{k!\ (1-\varrho(t))^k}{(1+\theta(1-\varrho(t)))^{k+1}}=\underset{n\rightarrow \infty}{\lim}\frac{k!\  (A_{n,u_n(t)-1})^{k}}{(1+\theta A_{n,u_n(t)-1})^{k+1}},$$
	uniformly on $[0,T]$. Now, let $ m<n$ such that $\rho_m<\varrho(T)\rho_n$. By definition of $a_{n}^{(e)}$, $A_{n,m-1}$ and equation \eqref{moved environment}, we have 
	\begin{equation*}
	\begin{split}
	\frac{a_{n}^{(e)}A_{n,m-1}
		\p^{(e_{m})}( Z_{n-m} > 0)\mu_m}{ \mu_n}=\frac{\rho_n-\rho_m}{2}\p^{(e)}( Z_{n} > 0\mid Z_m=1)\mu_m.
	\end{split}
	\end{equation*}
	Then, by using that   $\rho_{u_n(t)}/\rho_n\rightarrow \varrho(t)$ uniformly  for $0\leq t\leq T$ and applying the behaviour of the survival probability \eqref{eq: lemma 2}, the result holds for every fixed $k$.
	 Since we have a finite quantity of  $j\leq k$, the uniform convergence on all possible values of $j$ also holds true. 
\end{proof}

Similar arguments allow us to prove the following lemma. 

\begin{lemma}\label{lem: Yaglom for Q}	
	Let $0<T<1$, $k\in \mathbb{N}$, and  $\lambda>0$. Suppose that $Z$ is a regular critical GWVE. We also consider $u_n:[0,1]\rightarrow \{0, 1,\dots, n\}$, a sequence of functions such that
	$\rho_{u_n(t)}/\rho_n\rightarrow \varrho(t)$, uniformly,  for $0\leq t\leq T$, as $n\rightarrow \infty$, where $\varrho:[0,1]\rightarrow [0,1]$ is a (non-decreasing) function such that $\varrho(t)<1$, for all $t<1$.
	Then,
	\begin{equation*}
	\begin{split}
	\underset{n\rightarrow \infty}{\lim}	\mathbf{Q}^{(e_{u_n(t)},j,\theta_n)}_{n-u_n(t)}\left[e^{-(\lambda-\theta) Z_{n-u_n(t)}/a_n^{(e)}
	}\right]=
	& \frac{(1+\theta(1-\varrho(t)))^{j+1}}{(1+\lambda(1-\varrho(t)))^{j+1}},
	\end{split}
	\end{equation*}
	uniformly on $0\le  t\leq T$ and $j\leq k.$
\end{lemma}

\begin{proof}
	By definition of $\mathbf{Q}^{(e_{u_n(t)},k,\theta_n)}_{n-u_n(t)}$, we have
	\begin{equation*}
	\begin{split}
	\mathbf{Q}^{(e_{u_n(t)},k,\theta_n)}_{n-u_n(t)}\left[e^{-(\lambda-\theta) Z_{n-u_n(t)}/a_n^{(e)}}\right]&=
	\frac{\e^{(e_{u_n(t)})}\left[ e^{-\lambda Z_{n-u_n(t)}/a_{n}^{(e)}}Z_{n-u_n(t)}^{[k]}\right]}
	{\e^{(e_{u_n(t)})}\left[ e^{-\theta Z_{n-u_n(t)/a_{n}^{(e)}}}
		Z_{n-u_n(t)}^{[k]}\right]}\\
	&=\frac{\e^{(e_{u_n(t)})}\left[\left. e^{-\lambda Z_{n-u_n(t)}/a_{n}^{(e)}}Z_{n-u_n(t)}^{[k]}\right| Z_{n-u_n(t)}>0\right]}
	{\e^{(e_{u_n(t)})}\left[\left. e^{-\theta Z_{n-u_n(t)/a_{n}^{(e)}}}
		Z_{n-u_n(t)}^{[k]}\right| Z_{n-u_n(t)}>0\right]}.
	\end{split}
	\end{equation*}
	The result then follows from the asymptotic in  \eqref{eq:limit size biased}.
\end{proof}

Recall  that $\hat\psi_1$ and  $\psi_1$ denote the first time the spines split apart and the last time all the spines are together, respectively, and that $\hat\psi_1=\psi_1+1$. Now, we compute the limiting distribution of $\psi_1$ under $\mathbf{Q}^{(e_{\cdot},k,\theta_n)}_{n-\cdot}$, as $n\rightarrow \infty$, and  show  that asymptotically, the first split is binary and that the number of marks following each of those spines is uniformly distributed on $1,\dots, k$.  Moreover, we also deduce that the first spine splitting time and the sizes of the subgroups are asymptotically independent.

We also recall from \eqref{def: C} the definition of $C_{k;k_1,\dots, k_g}$ 
for  $g\geq 2$ and $\{k_i, i\geq g\}$ satisfying $k_1+\cdots+k_g=k$.

\begin{proposition}\label{splitting groups}
	Assume that conditions  \eqref{eq_cond_kersting} and \eqref{def: critical} are  satisfied  and let $0\leq T<1$. Consider $u_n:[0,1]\rightarrow \{0, 1,\dots, n\}$, a sequence of functions such that $\rho_{u_n(r)}/\rho_n\rightarrow \varrho(r)$, uniformly,  for $0\leq r\leq T$, as $n\rightarrow \infty$, where $\varrho:[0,1]\rightarrow [0,1]$ is a (non-decreasing) function such that $\varrho(r)<1$, for all $r<1$.  Let $k_1$ and $ k_2$ such that $k_1+k_2=k$. Then,  uniformly on $0\le  \varrho(r)< t\leq T$, we have
	\begin{equation}\label{psi}
	\underset{n\rightarrow\infty}{\lim}\mathbf{Q}^{(e_{u_n(r)},k, \theta_n)}_{n-u_n(r)}\Big(\psi_1\geq \tau_n(t)-u_n(r)\Big)= \left(	\frac{1+\theta(1-\varrho(r))}{1-\varrho(r)}\right)^{k-1}
	\left(	\frac{1-t}{1+\theta(1-t)}\right)^{k-1},
	\end{equation}
	and
	\begin{equation}\label{groups}
	\begin{split}
	\underset{n\rightarrow\infty}{\lim}\mathbf{Q}^{(e_{u_n(r)},k,\theta_n)}_{n-u_n(r)}\left( C_{k;k_1, k_2}\right)=\frac{1+\indi_{\{k_1\neq k_2\}}}{k-1}.
	\end{split}
	\end{equation}
	Moreover, we have the following asymptotic independence
	\begin{equation}\label{splitting}
	\begin{split}
	\underset{n\rightarrow\infty}{\lim}&\mathbf{Q}^{(e_{u_n(r)},k,\theta_n)}_{n-u_n(r)}\Big(\psi_1\in [0,\tau_n(t)-u_n(r)], C_{k;k_1, k_2}
	\Big)\\
	&=\underset{n\rightarrow\infty}{\lim}\mathbf{Q}^{(e_{u_n(r)},k,\theta_n)}_{n-u_n(r)}\left(  C_{k;k_1, k_2}
	\right)\underset{n\rightarrow\infty}{\lim}\mathbf{Q}^{(e_{u_n(r)},k,\theta_n)}_{n-u_n(r)}\left(\psi_1\in [0,\tau_n(t)-u_n(r)] \right).
	\end{split}
	\end{equation}
\end{proposition}
It is important to note that the  limit in \eqref{splitting} may be given in terms of $\rho$, that is 
\begin{equation}\label{splitting2}
\begin{split}
\underset{n\rightarrow\infty}{\lim}&\mathbf{Q}^{(e_{u_n(r)},k,\theta_n)}_{n-u_n(r)}\Big(\rho_{\psi_1+u_n(r)}/\rho_n \in [0,t] , C_{k;k_1, k_2}
\Big)\\
&\hspace{2cm}=\underset{n\rightarrow\infty}{\lim}\mathbf{Q}^{(e_{u_n(r)},k,\theta_n)}_{n-u_n(r)}\Big(\psi_1\in [0,\tau_n(t)-u_n(r)], C_{k;k_1, k_2}
\Big),
\end{split}
\end{equation}
which follows from the  definition of $\tau_n$, that is
\begin{equation}\label{eq: psi rho}
\{\psi_1\leq \tau_n(t)-u_n(r) \}
=\{\rho_{\psi_1+u_n(r)}\leq t\rho_n\}.
\end{equation}
Additionally, we also note that  $\rho_{\psi_1+u_n(r)}/\rho_n, $ under $\mathbf{Q}^{(e_{u_n(r)},k,\theta_n)}_{n-u_n(r)}$ converges in distribution to a random variable, here denoted by $\tilde{\varrho}(r)$. By differentiating equation  \eqref{psi}, its density satisfies
\begin{equation}\label{density splitting 2}
\left(	\frac{1+\theta(1-\varrho(r))}{1-\varrho(r)}
\right)^{k-1}
\frac{(k-1)(1-t)^{k-2}}{(1+\theta(1-t))^{k}}, \qquad\textrm{for}\quad t\in [\varrho(r),1].
\end{equation}

Equation \eqref{groups} tell us that the probability of  binary splitting in the limit is positive. By summing over $k_1$ and $k_2$, we  discover that this probability is actually one, implying the absence of multiple merges in the limit. The proof establishes this and the presence of asymptotic independence.
Intuitively, we only see binary splitting (or branching) in the sample genealogical tree in finite variance cases due to the following explanation. Imagine the constant environment GW case first for simplicity, as GWVE is essentially the same with some technicalities of course. Note that the probability, starting with a single individual, of a critical GW surviving at time $N$ is of asymptomatic order of magnitude $1/N$. Given a particle survives to some large time $N$ and it gives birth at time $m$, then chance exactly one offspring has a descendant that survives to time $N$ is $1-O(1/(N-m))$, the chance exactly two offspring have surviving offspring is $O(1/(N-m))$, exactly three offspring  with probability $O(1/(N-m))^2$, etc.  Then, over a time period of large length of order $N$, this will approximately give a Poisson number of binary branching events with each subtree surviving to the end, but splits into 3 or more are far too rare to be seen in the limit.  
Notably, in the heavy-tail case  \cite{HJP2022}, this intuitive explanation does not hold true.

\begin{proof}[Proof of Proposition \ref{splitting groups}]
	Let $0\leq \varrho(r)< t\leq T<1$ and define
	$$B_{t,n}:=	\mathbf{Q}^{(e_{u_n(r)},k,\theta_n)}_{n-u_n(r)}\Big(\psi_1\in [0,\tau_n(t)-u_n(r)], C_{k;k_1, k_2}
	\Big).$$
	First,  we use equation \eqref{eq: g groups} with the environment $e_{u_n(r)}$, final generation $n-u_n(r)$ and  $j \in[0,\tau_n(t)-u_n(r)]$. Indeed,  we are using that the generating function at generation $\ell$ is $f_{u_n(r)+j +1}$ to obtain
\begin{equation*}
	\begin{split}
	&B_{t,n}=	
\underset{j=0}{\overset{\tau_n(t)-u_n(r)}{\sum}}
\frac{k!
\left.f_{u_n(r)+j+1}''(s)\right|_{s=	f_{u_n(r)+j+1,n}(e^{-\theta_n})}
	\e^{(e_{u_n(r)})}\left[e^{-\theta Z_{n-u_n(r)}}Z_{n-u_n(r)}\right]}
{ \left(\underset{r=1}{\overset{k}{\prod}}g_r!\right) \left( \underset{i=1}{\overset{2}{\prod}}k_i!\right)
	\e^{(e_{u_n(r)+j})}\left[e^{-\theta Z_{n-u_n(r)-j}}Z_{n-u_n(r)-j}\right]
}\\
&\hspace{5cm}\times
\frac{	\underset{i=1}{\overset{2}{\prod}}\e^{(e_{u_n(r)+j+1})}\left[e^{-\theta Z_{n-u_n(r)-j-1}}Z_{n-u_n(r)-j-1}^{[k_i]}\right]}
{\e^{(e_{u_n(r)})}\left[e^{-\theta_nZ_{n-u_n(r)}}Z_{n-u_n(r)}^{[k]}\right]}.
\end{split}
\end{equation*}
Then, we perform the change of variable $m\mapsto u_n(r)+j$ and take out all the terms that do not depend on $m$,  to get
\begin{equation*}
\begin{split}
&B_{t,n}=
\frac{k!}{\left(\underset{r=1}{\overset{k}{\prod}}g_r!\right) \left( \underset{i=1}{\overset{2}{\prod}}k_i!\right)}
\frac{\e^{(e_{u_n(r)})}\left[e^{-\theta_nZ_{n-u_n(r)}}Z_{n-u_n(r)}\right]}
{\e^{(e_{u_n(r)})}\left[e^{-\theta_nZ_{n-u_n(r)}}Z_{n-u_n(r)}^{[k]}\right]
}\\
	&\hspace{3cm}\times \underset{m=u_n(r)}{\overset{\tau_n(t)}{\sum}}
\left.f_{m+1}''(s)\right|_{s=	f_{m+1,n}(e^{-\theta_n})}	\frac{
\underset{i=1}{\overset{2}{\prod}}\e^{(e_{m+1})}
\left[e^{-\theta Z_{n-m-1}}Z_{n-m-1}^{[k_i]}\right]}{\e^{(e_{m})}\left[e^{-\theta_{n} Z_{n-m}}Z_{n-m}\right]}.
	\end{split}
	\end{equation*}
	Now, from (3.5) in  \cite{kersting2021genealogical}, we have that for  $r_{m,n}\in [f_{m,n}(0),1]$,
	\begin{equation*}
	f_m''(r_{m,n})=f_m''(1)(1+o(1)), \qquad \textrm{as}\quad n\to \infty, 
	\end{equation*}
	uniformly for all $0\leq m\leq n$ such that $\rho_m/\rho_n\leq T<1$. Since $m\leq \tau_n(t)$ and $ t<1$, we use the previous asymptotic and
\eqref{Fknm} but in terms of  expectations 
for different values of $k,n$ and $m$,  to rewrite the above identity as follows 
	\begin{equation*}
	\begin{split}
	B_{t,n}=
	&\frac{(1+o(1))}{\underset{r=1}{\overset{k}{\prod}}g_r!}
\frac{F_{1,n,u_n(r)}}{F_{k,n,u_n(r)}}
\frac{(1+\theta A_{n,u_n(r)-1})^{k-1}	}
{(A_{n,u_n(r)-1})^{k-1}} \\
	&\hspace{1cm}
	\times\underset{m=u_n(r)}{\overset{\tau_n(t)}{\sum}}
\frac{f_{m+1}''(1)
		\mu_m\mu_n}
	{a_{n}^{(e)}\mu_{m+1}^2}	\frac{F_{k_1,n,m+1}F_{k_2,n,m+1}	}
	{F_{1,n,m}}
	\frac{(A_{n,m})^{k-2}(1+\theta A_{n,m-1})^{2} }
	{(1+\theta A_{n,m})^{k+2}}
	,
	\end{split}
	\end{equation*}
	for $n$ large enough. By definition of $\nu_{m+1}$,  $a_{n}^{(e)}$, $A_{n,m-1}$ and $A_{n,m}$,  we have
	$$	\frac{f_{m+1}''(1)\mu_m\mu_n }{\mu_{m+1}^2a_{n}^{(e)}}=2\frac{\nu_{m+1}}{\mu _{m}} \left(\sum_{k=0}^{n-1} \frac{\nu_{k+1}}{\mu _{k}}\right)^{-1}= 2[A_{n,m-1}-A_{n,m}],$$
	and
	$$\left(\frac{1+\theta A_{n,m-1}}{1+\theta A_{n,m}}\right)^2=1+ \frac{2\theta (A_{n,m-1}-A_{n,m})}{1+\theta A_{n,m}}+\frac{\theta^2 (A_{n,m-1}-A_{n,m})^2}{\left(1+\theta A_{n,m}\right)^2}.
	$$
 Replacing the previous identities in $B_{t,n}$ and using the uniform convergence of  \eqref{Fknm} to bound $$\frac{F_{1,n,u_n(r)}}{F_{k,n,u_n(r)}}\frac{F_{k_1,n,m+1}F_{k_2,n,m+1}	}
	{F_{1,n,m}}=1+o(1),$$ we obtain
	\begin{equation}\label{3sums}
	\begin{split}
	B_{t,n}=&\frac{(1+o(1))}{\underset{r=1}{\overset{k}{\prod}}g_r!}\frac{(1+\theta A_{n,u_n(r)-1})^{k-1}	}
	{(A_{n,u_n(r)-1})^{k-1}}
	\underset{m=u_n(r)}{\overset{\tau_n(t)}{\sum}}\left[\frac{2(A_{n,m})^{k-2}\left(A_{n,m-1}-A_{n,m}\right)}{(1+\theta A_{n,m})^{k}}\right.\\
	&\left.+\frac{4\theta(A_{n,m})^{k-2}\left(A_{n,m-1}-A_{n,m}\right)^2}{(1+\theta A_{n,m})^{k+1}}+\frac{2\theta^2(A_{n,m})^{k-2}\left(A_{n,m-1}-A_{n,m}\right)^3}{(1+\theta A_{n,m})^{k+2}}\right],
	\end{split}
	\end{equation}
	for $n$ large enough. 
	
	Next, we recall that $R^{(n)}= \{0= A_{n,n-1} < A_{n,n-2} < \ldots < A_{n,0}<A_{n,-1}=1\}$ is a partition of $[0,1]$ satisfying \eqref{limit norm Rn} and observe that in the limit, $\{ A_{n,\tau_n(t)} < \ldots < A_{n,u_n(r)}\}$ is a partition of $[1-t,1-\varrho(r)]$. Therefore, for any continuous function $f$, the Riemann-Stieltjes integral can be perform as the limiting sum
$$\int_{1-t}^{1-\varrho(r)} f(x)\ud x = \underset{n\rightarrow \infty}{\lim}\underset{m=u_n(r)}{\overset{\tau_n(t)}{\sum}}f(A_{n,m})\left(A_{n,m-1}-A_{n,m}\right).$$
In particular, 
the first sum in the right-hand side of \eqref{3sums} has limit 
\begin{equation*}
	\begin{split}
\underset{n\rightarrow \infty}{\lim}\underset{m=u_n(r)}{\overset{\tau_n(t)}{\sum}}\frac{2(A_{n,m})^{k-2}}{(1+\theta A_{n,m})^{k}}\left(A_{n,m-1}-A_{n,m}\right)=& \int_{1-t}^{1-\varrho(r)} \frac{2x^{k-2}}{(1+\theta x)^{k}}\ud x\\
&= \left.\frac{2(1+\theta x)^{1-k}}{(k-1)x^{1-k}}\right|_{1-t}^{1-\varrho(r)}.
\end{split}
\end{equation*}
For the second sum in the right-hand side  of \eqref{3sums}, we use the definition of the norm of the partition and \eqref{limit norm Rn} to get 
\begin{equation*}
	\begin{split}
&\underset{n\rightarrow \infty}{\lim}\underset{m=u_n(r)}{\overset{\tau_n(t)}{\sum}}\frac{4\theta(A_{n,m})^{k-2}}{(1+\theta A_{n,m})^{k+1}}\left(A_{n,m-1}-A_{n,m}\right)^2\\
&\leq \underset{n\rightarrow \infty}{\lim}\|R^{(n)}\|\underset{m=u_n(r)}{\overset{\tau_n(t)}{\sum}}\frac{4\theta(A_{n,m})^{k-2}}{(1+\theta A_{n,m})^{k+1}}\left(A_{n,m-1}-A_{n,m}\right)\\
&=  0 \cdot \int_{1-t}^{1-\varrho(r)} \frac{4 \theta x^{k-2}}{(1+\theta x)^{k+1}}\ud x= 0.
\end{split}
\end{equation*}
Similarly,
	$$\underset{n\rightarrow \infty}{\lim}\underset{m=u_n(r)}{\overset{\tau_n(t)}{\sum}}\hspace{-.1cm}\frac{2\theta^2
	(A_{n,m})^{k-2}}{(1+\theta A_{n,m})^{k+2}}	\left(A_{n,m-1}-A_{n,m}\right)^3\leq 0\cdot \int_{1-t}^{1-\varrho(r)} \frac{2 \theta^2 x^{k-2}}{(1+\theta x)^{k+2}}\ud x= 0.$$
	By observing that $A_{n,u_n(r)-1}=A_{n,u_n(r)}+(A_{n,u_n(r)-1}-A_{n,u_n(r)})$, we take the limit in \eqref{3sums}  and recalling \eqref{limit norm Rn}, we deduce
	\begin{equation}\label{limit B}
	\underset{n\rightarrow \infty}{\lim}B_{t,n}=\frac{2}{(k-1)\underset{r=1}{\overset{k}{\prod}}g_r!}\left( 1- \left(	\frac{1+\theta(1-\varrho(r))}{1+\theta(1-t)}
	\right)^{k-1}
	\left(	\frac{1-t}{1-\varrho(r)}\right)^{k-1}\right).
	\end{equation}
	uniformly on $0\le  \varrho(r)< t\leq T$. 
	The identity in \eqref{groups} follows by  taking $t\uparrow 1$ above. More precisely, from the   Monotone Convergence Theorem we get
	\[
	\begin{split}
	\underset{n\rightarrow \infty}{\lim}&\mathbf{Q}^{(e_{u_n(r)},k,\theta_n)}_{n-u_n(r)}\left( C_{k;k_1, k_2}\right)=\frac{2}{(k-1)\underset{r=1}{\overset{k}{\prod}}g_r!}=\frac{1+\indi_{\{k_1\neq k_2\}}}{k-1}.
	\end{split}
	\]
	In order to obtain \eqref{psi}, we sum in \eqref{limit B} over all  possible choices of  $k_1$ and $k_2$ satisfying $k_1+k_2=k$.  Finally, from \eqref{psi}, \eqref{groups} and \eqref{limit B}, we see that equation \eqref{splitting} holds. 
\end{proof}

Now, we compute the limiting joint distribution of the spines split times and the splitting groups. Recall that $\psi_i$ denotes the  last time where there are at most $i$ marked particles. 
We will show  that asymptotically, all splittings are binary and that, knowing that there are $a$ marks in the particle that splits,   the number of marks  following one of these spines is uniformly distributed on $\{1,\ldots, a-1\}$. In other words, in the limit, if we 
squeeze or stretch the lines of the spine tree such that each of them has length one, then we obtain a  full (or proper) binary    tree with $k$ leaves. 

Let $\mathcal{B}^k$ be  the set of full unlabelled binary  trees with  $k$ leaves, i.e.   $\mathfrak{t}\in \mathcal{B}^k$ if and only if $l_u(\mathfrak{t})\in\{0,2\}$ for all $u\in \mathfrak{t}$ and $k=|\{u\in  \mathfrak{t}: l_u(\mathfrak{t})=0\}|$. Let
$N(\mathfrak{t})$ be  the set of internal vertices of $\mathfrak{t}$, that is  $N(\mathfrak{t})=\{u\in  \mathfrak{t}: l_u(\mathfrak{t})=2\}$.  For every internal vertex $u$, we define the vector $(k_{u1},k_{u2})$ where $k_{uj}=k_{uj}(\mathfrak{t}):=|\{\mbox{leaves in } S_{uj}(\mathfrak{t})\}|$ denotes the number of leaves in the subtree of $\mathfrak{t}$ with root at $uj$ for $j=1,2$. We endow $\mathcal{B}^k$ with the probability measure 
\begin{equation}
\label{eq: lawbinarytrees} \p_{\mathcal{B}^k}(\mathfrak{t}):=\underset{u\in N(\mathfrak{t})}{\prod}\frac{1+\indi_{\{k_{u1}\neq k_{u2}\}}}{k_{u1}+k_{u2}-1} \qquad \mbox{ for }\quad \mathfrak{t}\in \mathcal{B}^k.
\end{equation}
This probability law is equivalent to start with $k$ leaves at a given particle and each time that the particle branches, we select uniformly without labelling, how many of the leaves goes to each of the two branches.

\begin{figure}
		\includegraphics[width=10cm]{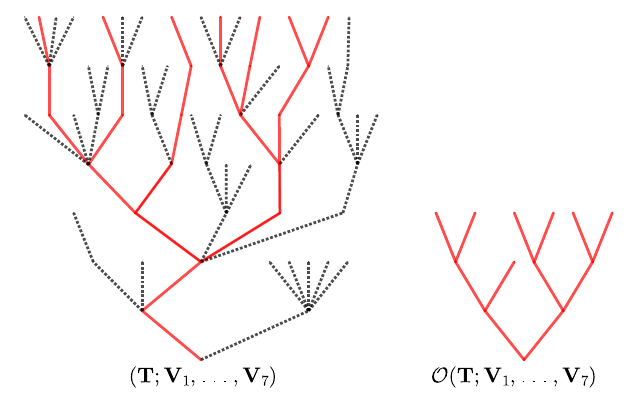}
		\caption{The action  $\mathcal{O}$ applied to $(\mathbf{T};\mathbf{V}_1,\dots,\mathbf{V}_7)$,   squeezes each line in the spines subtree to obtain a tree in $\mathcal{B}^7$.}
		\label{fig:6}
	\end{figure}

Denote by 
\[
\mathcal{O}:\{(\mathbf{t};\mathbf{v}_1,\dots,\mathbf{v}_k)\in \mathcal{T}^k: \mbox{spine splittings are  binary}\}\longrightarrow \mathcal{B}^k,
\]
the operation that squeeze each line  of the tree
$\underset{i=1}{\overset{k}{\bigcup }}\mathbf{v}_i$ in a way that we obtain a tree in $\mathcal{B}^k$, see Figure \ref{fig:6} for an illustrative example. For simplicity of exposition, in the next proposition we denote $\mathcal{O}=\mathcal{O}(\mathbf{T};\mathbf{V}_1,\dots,\mathbf{V}_k)$, the operation applied to a $\mathcal{T}^k$-valued random variable with law $\mathbf{Q}_n^{(e,k,\theta)}$.

\begin{proposition}\label{prop: uniform} 	
Assume that conditions  \eqref{eq_cond_kersting} and \eqref{def: critical} are  satisfied. 
	Let  $\mathfrak{t}\in \mathcal{B}^k$ and $0\leq s_1\leq t_1\leq s_2\leq \dots \leq s_{k-1} \leq  t_{k-1}\leq 1$.  Then 
	\begin{equation*}
	\begin{split}
	\underset{n\rightarrow \infty}{\lim}\mathbf{Q}^{(e,k,\theta_n)}_n&\Big(\psi_1\in (\tau_n(s_1),\tau_n(t_1)), \dots ,\psi_{k-1}\in (\tau_n(s_{k-1}),\tau_n(t_{k-1})), \mathcal{O}=\mathfrak{t} \Big)\\
	&\hspace{-1cm}=\left(\underset{u\in N(\mathfrak{t})}{\prod}\frac{1+\indi_{\{k_{u1}\neq k_{u2}\}}}{k_{u1}+k_{u2}-1}\right)
	(k-1)!\left(1+\theta\right)^{k-1}\underset{i=1}{\overset{k-1}{\prod}}\left(\tfrac{1-s_i}{1+\theta(1-s_i)}-\tfrac{1-t_i}{1+\theta(1-t_i)}\right).
	\end{split}
	\end{equation*} 
	In particular, the sequence $\{\rho_{\psi_i}/\rho_n: i=1,\dots, k-1\}$ converges in distribution, under $\mathbf{Q}^{(e,k,\theta_n)}_n$, to an  ordered sample of $k-1$ random variables with common distribution given by 
	\[
	F(t)=\frac{t}{1+\theta (1-t)}, \qquad \textrm{for} \quad t\in [0,1].
	\] 
\end{proposition}

\begin{proof}
	Consider $u_n:[0,1]\rightarrow \{0, 1,\dots, n\}$, a sequence of functions such that 
	\[
	\rho_{u_n(r)}/\rho_n\xrightarrow[n\to\infty]{} \varrho(r),
	\] uniformly, for $0\leq r\leq T<1$, where $\varrho:[0,1]\rightarrow [0,1]$ is a non-decreasing function such that $\varrho(r)<1$, for all $r<1$. By induction on $k$, we will  prove  below that  uniformly on $0\le  \varrho(r)< s_1$,
	\begin{equation}\label{induction}
	\begin{split}
&	\underset{n\rightarrow\infty}{\lim}\mathbf{Q}^{(e_{u_n(r)},k, \theta_n)}_{n-u_n(r)}
	\left(\underset{i=1}{\overset{k-1}{\prod}}
	\{\psi_i\in (\tau_n(s_i)-u_n(r),\tau_n(t_i)-u_n(r))\}, \mathcal{O}=\mathfrak{t} \right)\\
	&=\left(\underset{u\in N(\mathfrak{t})}{\prod}\frac{1+\indi_{\{k_{u1}\neq k_{u2}\}}}{k_{u1}+k_{u2}-1}\right)
	(k-1)!
	\left(\frac{1+\theta(1-\varrho(r))}{1-\varrho(r)}\right)^{k-1}\underset{i=1}{\overset{k-1}{\prod}}\left(\tfrac{1-s_i}{1+\theta(1-s_i)}-\tfrac{1-t_i}{1+\theta(1-t_i)}\right).
	\end{split}
	\end{equation}
	Therefore, by taking $r=0$ and $u_n(r)=0$, the first claim will follow and  the second one holds by using  identity \eqref{eq: psi rho} and summing over all possible $\mathfrak{t}\in\mathcal{B}^k$. 
	
	Let us prove \eqref{induction}. We start with $k=2$ and  observe that  in this case,  necessarily $\mathfrak{t}_0:=\{\emptyset, 1,2\}$ is the only possible tree in $\mathcal{B}^2$ with $k_{1}=k_{2}=1$. Then, the result holds true by Proposition \ref{splitting groups} with $C_{2,1,1}$. 
	
	Now, we suppose that the result is fulfilled for every $m<k$. 
	Since $\{\mathcal{O}=\mathfrak{t}\}$, at the first splitting time, the marks split in 2 groups of sizes $k_{1}$ and $k_{2}$, 
	respectively. In other words, we have the event $C_{k,k_{1},k_{2}}$ on the first splitting time.  We also observe that we can decompose $\mathfrak{t}$ 
	as its restriction to the first generation  concatenated  with the two subtrees attached at generation $1$,  that is 
	\[\mathfrak{t}=\mathfrak{t}_0\bigsqcup 1 S_{1}(\mathfrak{t})\bigsqcup 2 S_{2}(\mathfrak{t}),
	\]
	see Figure \ref{fig:4} for a similar decomposition.  
	
	Now, denote by $\textbf{i}^{(1)}=\{i^1_1,\dots, i^1_{k_{1}-1}\}$ and $\textbf{i}^{(2)}=\{i^2_1,\dots, i^2_{k_{2}-1}\}$ a partition of $\{2,3,\dots,k-1\}$,  that is   $\textbf{i}^{(1)}\cap \textbf{i}^{(2)}=\emptyset$ and $\textbf{i}^{(1)}\cup \textbf{i}^{(2)}=\{2,3,\dots,k-1\}.$ 
	By conditioning on $\mathcal{G}=\sigma(\psi_1, C_{k,k_{1},k_{2}})$, summing over all  possible partitions $(\textbf{i}^{(1)},\textbf{i}^{(2)})$ such that $\textbf{i}^{(1)}$ and $\textbf{i}^{(2)}$ are associated with $S_{1}(\mathfrak{t})$ and $S_{2}(\mathfrak{t})$, respectively;  and then using Proposition \ref{prop: branching property Q}, we  obtain
	\begin{equation*}
	\begin{split}
	&\mathbf{Q}^{(e_{u_n(r)},k, \theta_n)}_{n-u_n(r)}
	\left[\mathbf{Q}^{(e_{u_n(r)},k, \theta_n)}_{n-u_n(r)}\left(\left.
	\underset{i=1}{\overset{k-1}{\prod}}
	\{\psi_i\in (\tau_n(s_i)-u_n(r),\tau_n(t_i)-u_n(r))\}, \mathcal{O}=\mathfrak{t}\  \right| \mathcal{G}\right) \right]\\
	&=\mathbf{Q}^{(e_{u_n(r)},k, \theta_n)}_{n-u_n(r)}
	\Bigg[\mathbf{1}_{\{\psi_1\in (\tau_n(s_1)-u_n(r),\tau_n(t_1)-u_n(r)), C_{k;k_1,k_2}\}} \underset{(\textbf{i}^{(1)},\textbf{i}^{(2)})}{\sum}  	\prod_{j=1}^2\Bigg.\\
	&\Bigg.
	\mathbf{Q}^{(e_{u_n(r)+\psi_1},k_{j}, \theta_n)}_{n-u_n(r)-\psi_1}\Bigg(
	\underset{\ell=1}{\overset{k_{j}-1}{\prod}}
	\{\widetilde{\psi}_\ell^{(j)}\in (\tau_n(s_{i^j_{\ell}})-u_n(r)-\psi_1,\tau_n(t_{i^j_{\ell}})-u_n(r)-\psi_1)\}, \mathcal{O}=S_{j}(\mathfrak{t})\Bigg) \Bigg],
	\end{split}
	\end{equation*}
	where $\widetilde{\psi}_\ell^{(j)}=\psi_{i_\ell^{(j)}}-\psi_1$ are the splitting times associated with the subtree $S_j(\mathfrak{t})$. 
	
	We recall from  \eqref{splitting2}  that  $\rho_{\psi_1+u_n(r)}/\rho_n$, under $\mathbf{Q}^{(e_{u_n(r)},k,\theta_n)}$, converges in distribution, uniformly on $0\le  \varrho(r)< t\leq T$, towards the r.v. $\widetilde\rho(r)$ whose density function is given by \eqref{density splitting 2}. Thus by the Skorokhod's Representation Theorem \cite[Theorem 6.7]{Billingsley}, we can construct $\widetilde{\varrho}$ and  the sequence  $(\rho_{\psi_1+u_n(r)}/\rho_n)_{n\ge 1}$, on a common probability space in such a way that the convergence is for almost every $\omega$. Let $N=N(S_{1}(\mathfrak{t}))\cup N(S_{2}(\mathfrak{t}))$.
	 For every fixed $(\textbf{i}^{(1)},\textbf{i}^{(2)})$, we use the induction hypothesis \eqref{induction} with $\widetilde{u}_n=\psi_1+u_n(r)$ to obtain,  after rearranging terms,  that 
$$\prod_{j=1}^2
	\mathbf{Q}^{(e_{u_n(r)+\psi_1},k_{j}, \theta_n)}_{n-u_n(r)-\psi_1}\Bigg(
	\underset{\ell=1}{\overset{k_{j}-1}{\prod}}
	\{\widetilde{\psi}_\ell^{(j)}\in (\tau_n(s_{i^j_{\ell}})-u_n(r)-\psi_1,\tau_n(t_{i^j_{\ell}})-u_n(r)-\psi_1)\}, \mathcal{O}=S_{i}(\mathfrak{t})\Bigg)$$
	converges for almost every $\omega$ as $n\rightarrow\infty$ to 
	\[
	\begin{split}
	\left(\underset{u\in N}{\prod}\frac{1+\indi_{\{k_{u1}\neq k_{u2}\}}}{k_{u1}+k_{u2}-1}\right)
	(k_{1}-1)!(k_{2}-1)!
	&\left(\frac{1+\theta(1-\widetilde{\varrho}(r))}{1-\widetilde{\varrho}(r)}\right)^{k-2}\\
	& \underset{i=2}{\overset{k-1}{\prod}}\left(\tfrac{1-s_i}{1+\theta(1-s_i)}-\tfrac{1-t_i}{1+\theta(1-t_i)}\right).
	\end{split}
\]
	Note that this is independent of the choice of the partition $(\textbf{i}^{(1)},\textbf{i}^{(2)})$  and that there are $(k-2)!/(k_{1}-1)!(k_{2}-1)!$ ways to choose them. Therefore, we apply Proposition  \ref{splitting groups} or more precisely the density of the limiting distribution function \eqref{splitting2}, which can be derived from \eqref{splitting} using the explicit identity in  \eqref{groups} and the density in \eqref{density splitting 2}; 
	and  the Dominated Convergence Theorem to obtain
	\begin{equation*}
	\begin{split}
&	\underset{n\rightarrow\infty}{\lim}\mathbf{Q}^{(e_{u_n(r)},k, \theta_n)}_{n-u_n(r)}
	\left(\underset{i=1}{\overset{k-1}{\prod}}
	\{\psi_i\in (\tau_n(s_i)-u_n(r),\tau_n(t_i)-u_n(r))\}, \mathcal{O}=\mathfrak{t} \right)
	\\
	&=\left(\frac{1+\indi_{\{k_{1}\neq k_{2}\}}}{k-1}\right)\left(\underset{u\in N}{\prod}\frac{1+\indi_{\{k_{u1}\neq k_{u2}\}}}{k_{u1}+k_{u2}-1}\right)(k-2)!\underset{i=2}{\overset{k-1}{\prod}}\left(\tfrac{1-s_i}{1+\theta(1-s_i)}-\tfrac{1-t_i}{1+\theta(1-t_i)}\right)\\
	&\hspace{1.5cm} 
	\left(	\frac{1+\theta(1-\varrho(r))}{1-\varrho(r)}
	\right)^{k-1}\int_{s_1}^{t_1}
	\left(\frac{1+\theta(1-v)}{1-v}\right)^{k-2}
	\frac{(k-1)(1-v)^{k-2}}{(1+\theta(1-v))^{k}}
	\ud v.
	\end{split}
	\end{equation*}
	After cancellation by performing the integral and noticing that $N(\mathfrak{t})=\{\emptyset\}\cup N$, we get \eqref{induction}. This completes the proof.
\end{proof}

From the previous result, we may observe that, in the  limit, the tree topology and the split times are asymptotically independent, that is 
\begin{equation*}
\begin{split}
\underset{n\rightarrow \infty}{\lim}\mathbf{Q}^{(e,k,\theta_n)}_n&\left(\underset{i=1}{\overset{k-1}{\prod}}
\{\psi_i\in (\tau_n(s_i),\tau_n(t_i))\}, \mathcal{O}=\mathfrak{t}\right)\\
&=\underset{n\rightarrow \infty}{\lim}\mathbf{Q}^{(e,k,\theta_n)}_n\left( \overset{k-1}{\prod}
\{\psi_i\in (\tau_n(s_i),\tau_n(t_i))\}
\right)\times \underset{n\rightarrow \infty}{\lim}\mathbf{Q}^{(e,k,\theta_n)}_n\left(\mathcal{O}=\mathfrak{t}\right).
\end{split}
\end{equation*}
Moreover,  we can also deduce that  if we start with  $i$ groups of spines of sizes $a_1,\dots, a_{i}$, in the limit, the split times for any group $j$ will be distributed like $a_j-1$ independent random variables, all with the same distribution. 
In particular, this implies that the first group to split will be group $j$ with probability proportional to $a_j-1$, that is, with probability $(a_j-1)/(k-i)$. (In fact, as these limiting probabilities do not depend on $\theta$, this holds starting with the $i+1$ spine groups at any time $\tau_n(t)$ for $t\in[0,1)$ and not just from time $0$, likewise we can condition on knowing the $i+1$ spine group sizes at the time of the $i$-th split, $\psi_i$). We note analogous limiting spine tree topology results appear in \cite{HJRarxiv} (e.g. lemma 30) that are essentially identical.  
We want to emphasize that the previous limit is under $\mathbf{Q}^{(e,k,\theta_n)}_n$, while the limit in our main result is under  $\mathbb{P}^{(e)}$. Roughly speaking, since  $\mathbf{Q}^{(e,k,\theta_n)}_n$ is the  $Z_n^{[k]}e^{-\theta Z_n}$-size biased transform of $\mathbb{P}^{(e)}$, the tree topology and the split times are asymptotically independent under $\mathbb{P}^{(e)}$. However, under $\mathbb{P}^{(e)}$, the limiting split times would not be longer an  ordered sample of $k-1$ independent random variables with common distribution. The formal proof is development in the next section.

Now, let  $(\widetilde{\psi}_1,\dots, \widetilde{\psi}_{k-1})$ be a uniformly random permutation of $({\psi}_1,\dots {\psi}_{k-1})$. As a consequence of the previous proposition we have the next corollary.

\begin{corollary}\label{cor: uniform}
	The spine splitting times $\{\widetilde\psi_i, i=1,\dots,k-1\}$ and $\{\mathcal{O}=\mathfrak{t},\mathfrak{t}\in\mathcal{B}^k\}$ are asymptotically independent. Moreover, for  $\{t_1,\dots,t_{k-1}\}\subset (0,1)$and $\mathfrak{t}\in \mathcal{B}^k$
	\begin{equation*}
	\begin{split}
	\underset{n\rightarrow \infty}{\lim}
	\mathbf{Q}^{(e,k,\theta_n)}_n\left(\underset{i=1}{\overset{k-1}{\prod}}\{\widetilde{\psi}_i \leq \tau_n(t_i)\}  \right)
	&=\underset{i=1}{\overset{k-1}{\prod}}\frac{t_i}{1+\theta (1-t_i) }\\
	\underset{n\rightarrow \infty}{\lim}
	\mathbf{Q}^{(e,k,\theta_n)}_n\left( \mathcal{O}=\mathfrak{t}\right)
	&=\underset{u\in N(\mathfrak{t})}{\prod}\frac{1+\indi_{\{k_{u1}\neq k_{u2}\}}}{k_{u1}+k_{u2}-1}=\p_{\mathcal{B}^k}(\mathfrak{t}). 
	\end{split}
	\end{equation*} 
\end{corollary}

Denote by  $\mathcal{D}_n$ the event that all spines splitting times are different. Moreover, since in the limit there are only binary splittings we  restrict ourselves on trees in $\mathcal{B}^k$. A direct consequence of Proposition \ref{prop: uniform}  is that $\underset{n\rightarrow \infty}{\lim}\mathbf{Q}^{(e,k,\theta_n)}_n\left(\mathcal{D}_n\right)=1$. We can use Lemma \ref{lem: Yaglom for Q} to deduce the next lemma. 

\begin{lemma}\label{lem:limit Q conditioned psi} 
Assume that conditions  \eqref{eq_cond_kersting} and \eqref{def: critical} are  satisfied. 
	Consider $\mathfrak{t}\in \mathcal{B}^k$,  $\lambda>0$ and  $\{t_1,\dots,t_{k-1}\}\subset (0,1)$. Then, 
	\begin{equation*}
	\begin{split}
	\underset{n\rightarrow \infty}{\lim}	\mathbf{Q}^{(e,k,\theta_n)}_{n}&\left[e^{-(\lambda-\theta) Z_n/a_n^e}\Big| \widetilde{\psi}_1=\tau_n(t_1), \dots, \widetilde{\psi}_{k-1}=\tau_n(t_{k-1}),\mathcal{D}_n,\mathcal{O}=\mathfrak{t}\right]\\
	&\hspace{6cm}= \frac{(1+\theta)^{2}}{(1+\lambda)^{2}}\prod_{i=1}^{k-1}\frac{(1+\theta(1-t_i))^2}{(1+\lambda(1-t_i))^2}.
	\end{split}
	\end{equation*}
\end{lemma}

\begin{proof}
	For our purposes, let us  rephrase  Proposition \ref{prop: forward construction} as follows: every marked particle at generation $m\in \{0,\dots,n-1\}\setminus \{\psi_1,\dots,\psi_{k-1}\}$ has offspring according to  $q_{m+1}^{(1,\theta)}$. Moreover, for each $i\leq k-1$, at generation $\psi_i$ there are $i$ different marked particles.
	The marked particle that branches, say $w^{*i}$, has offspring according to $q_{\psi_i+1}^{(2,\theta)}$.  The other $i-1$ marked particles give birth according to $q_{\psi_i+1}^{(1,\theta)}$.
	
	Let $\mathcal{W}$ be the labels  of the different marked particles and $\widetilde{\mathcal{W}}=\mathcal{W}\setminus\{w^{*i}: i=1,\dots k-1\}$. Recall the definition of $Y_n(w)$ in the discussion above Lemma \ref{lem:  constructed block}, that is  the total population size  of unmarked particles at generation  $n$ that descend from the unmarked offspring of  $w$. By Proposition \ref{prop: branching property Q}, we have
	$$Z_n=k+\underset{w\in\mathcal{W}}{\sum}Y_n(w),$$ 
	where $(Y_n(w))_{w\in \mathcal{W}}$ are independent r.v. (see Figure \ref{fig:5} for an example). 
	
	\begin{figure}
		\includegraphics[width=12cm]{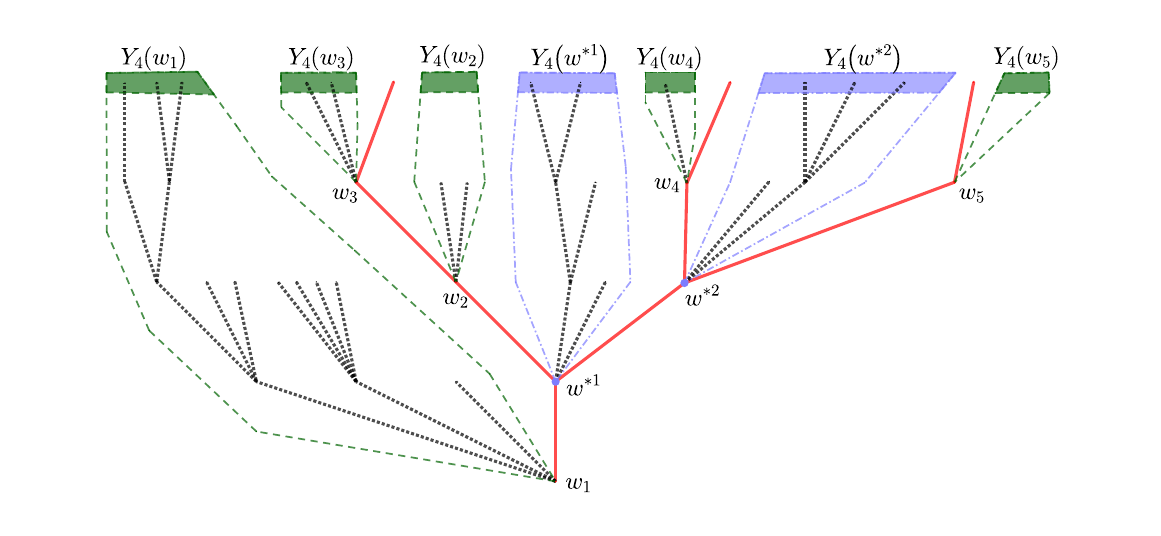}
		\caption{Decomposition of $Z_4$ into  the total population size of unmarked
			particles who are descendants of the unmarked offspring of spine particles. The spines are the red  lines, the particles $w_1,w_2, w_3, w_4$ and $w_5$ give birth according to  $q_{\cdot}^{(1,\theta)}$ and their unmarked progeny  at generation 4 are $(Y_4(w_i))_{1\le i\le 5}$; $w^{*1}$ and $w^{*2}$ give birth according to $q_{\cdot}^{(2,\theta)}$ and their unmarked progeny  at generation 4 are $Y_4(w^{*1})$ and $Y_4(w^{*2})$.}
		\label{fig:5}
	\end{figure}

	Then  Lemma \ref{lem: constructed block} implies
	\begin{equation*}
	\begin{split}
&	\mathbf{Q}^{(e,k,\theta)}_{n}
		\left[e^{-(\lambda-\theta) Z_n}\Big| \psi_1, \dots, \psi_{k-1},\mathcal{D}_n,\mathcal{O}=\mathfrak{t}\right]\\
	&=\frac{e^{-k\lambda}}{e^{-k\theta}}\prod_{w \in \widetilde{\mathcal{W}}} \frac{f_{|w|+1}'(f_{|w|+1,n}(e^{-\lambda}))}{f_{|w|+1}'(f_{|w|+1,n}(e^{-\theta}))}\prod_{i=1}^{k-1} \frac{f_{|w^{*i}|+1}''(f_{|w^{*i}|+1,n}(e^{-\lambda}))}{f_{|w^{*i}|+1}''(f_{|w^{*i}|+1,n}(e^{-\theta}))} \\
	&=\frac{e^{-k\lambda}}{e^{-k\theta}}\prod_{w \in \mathcal{W}} \frac{f_{|w|+1}'(f_{|w|+1,n}(e^{-\lambda}))}{f_{|w|+1}'(f_{|w|+1,n}(e^{-\theta}))}\\
	&\hspace{3cm}\prod_{i=1}^{k-1} \frac{f_{|w^{*i}|+1}''(f_{|w^{*i}|+1,n}(e^{-\lambda}))}{f_{|w^{*i}|+1}''(f_{|w^{*i}|+1,n}(e^{-\theta}))} \frac{f_{|w^{*i}|+1}'(f_{|w^{*i}|+1,n}(e^{-\theta}))} {f_{|w^{*i}|+1}'(f_{|w^{*i}|+1,n}(e^{-\lambda}))}.
	\end{split}
	\end{equation*}
Note that $\{|w|: w\in W \mbox{ or } w\in\widetilde{W}\}$ are in the $\sigma$-algebra generated by the random variables and events on the left-hand side. 
	Observe that at every  $m\in\{\psi_{i-1}+1,\dots, \psi_{i}\}$, there are $i$ marked particles implying that  in the previous decomposition the factor ${f_{m+1}'(f_{m+1,n}(e^{-\lambda}))}/{f_{m+1}'(f_{m+1,n}(e^{-\theta}))}$ appears $i$ times. Thus by Lemma \ref{lem: constructed block}, we obtain
	\begin{equation*}
	\begin{split}
	\mathbf{Q}^{(e,k,\theta)}_{n}
	&	\left[e^{-(\lambda-\theta) Z_n}\Big| \psi_1, \dots, \psi_{k-1},\mathcal{D}_n,\mathcal{O}=\mathfrak{t}\right]\\
	=&
	\mathbf{Q}^{(e,1,\theta)}_{n}\left[e^{-(\lambda-\theta)  Z_n}\right] \prod_{i=1}^{k-1} \mathbf{Q}^{(e_{\psi_i+1},1,\theta)}_{n-\psi_i-1}\left[e^{-(\lambda-\theta)  Z_{n-\psi_i-1}}\right]
	\frac{h(n,\psi_i,\lambda)}{h(n,\psi_i,\theta)},
	\end{split}
	\end{equation*}
	where
	$$h(n,m,\lambda):=\frac{f_{m+1}''(f_{m+1,n}(e^{-\lambda}))} {f_{m+1}'(f_{m+1,n}(e^{-\lambda}))}.$$
	We also note that   the right-hand side of the previous identity does not depend  on the choice $\mathcal{O}=\mathfrak{t}$, in fact we only require to know that  the splitting times are different. 
	
	Now, if we conditioned on the values of $(\widetilde{\psi}_1,\dots, \widetilde{\psi}_{k-1})$  instead of $(\psi_1,\dots,\psi_{k-1})$, we obtain 
	\begin{equation*}
	\begin{split}
&	\mathbf{Q}^{(e,k,\theta_n)}_{n}
		\left[e^{-(\lambda-\theta) Z_n/a_n^e}\Big| \widetilde{\psi}_1=\tau_n(t_1), \dots, \widetilde{\psi}_{k-1}=\tau_n(t_{k-1}),\mathcal{D}_n,\mathcal{O}=\mathfrak{t}\right]\\
	&=
	\mathbf{Q}^{(e,1,\theta_n)}_{n}\left[e^{-(\lambda-\theta)  Z_n/a_n^e}\right] \prod_{i=1}^{k-1} \mathbf{Q}^{(e_{\tau_n(t_i)+1},1,\theta_n)}_{n-\tau_n(t_i)-1}\left[e^{-(\lambda-\theta)  Z_{n-\tau_n(t_i)-1}/a_n^e}\right]
	\frac{h\left(n,\tau_n(t_i),\lambda/a_n^e\right)}{h(n,\tau_n(t_i),\theta/a_n^e)}.
	\end{split}
	\end{equation*}
By taking $\epsilon=1/2$ in Condition \eqref{eq_cond_kersting}, we know that there exists a $c>0$ such that any  $n\ge 1$
\begin{equation*}
\begin{split}
\mathbb{E}\left[\left(\chi_{1}^{(n)}\right)^2\mathbf{1}_{\{\chi_{1}^{(n)}\geq 2\}}\right] &\leq 2 \mathbb{E}\left[\left(\chi_{1}^{(n)}\right)^2\mathbf{1}_{\{2\leq \chi_{1}^{(n)}\leq  c(1+\mathbb{E}[\chi_{1}^{(n)}])\}}\right]\\
&\leq 2 c\left(1+\mathbb{E}\left[\chi_{1}^{(n)}\right]\right)\mathbb{E}\left[\chi_{1}^{(n)}\mathbf{1}_{\{ \chi_{1}^{(n)}\geq 2\}}\right].
\end{split}
\end{equation*} 
Then, there exists a $C>0$ such that 
$$f_n''(1)\leq C f'_n(1)(1+f'_k(1)), \qquad \mbox{ for all } n\geq 1.$$
Therefore, we can apply 
\cite[Lemma 3]{CardonaPalau}, (see equation (25) inside the proof of Lemma 3), to 
deduce 
	\begin{equation}\label{cebolla}
	\underset{n\rightarrow \infty}{\lim}\ \underset{0\leq m\leq n}{\sup}\left|\frac{h\left(n,\tau_n(t_i),\lambda/a_n^e\right)}{h(n,\tau_n(t_i),\theta/a_n^e)}-1
	\right|=0.
	\end{equation}
	Therefore, by \eqref{limit rho} and Lemma \ref{lem: Yaglom for Q}, we get the result.
\end{proof}

With all ingredients in hand, we now proceed with the proof of our main result.
\section{Proof of Theorem \ref{principal}}
\begin{proof}[Proof of Theorem \ref{principal}]
	For a given $\theta\geq 0$, recall that $\mathbf{Q}^{(e,k,\theta_n)}_{n}$ is the $e^{-\theta_n Z_n}Z_n^{[k]}$-sized biased transform of $\p^{(e)}$ with $\theta_n=\theta/a_{n}^{(e)}$. We also recall that   $\mathbf{Q}^{(e,k,\theta_n)}_{n}(Z_n\geq k)=1$ and  the permuted splitting times  are denoted by $(\widetilde{\psi}_1,\dots, \widetilde{\psi}_{k-1})$. 	By using  conditional probability, we  see
	\begin{equation}\label{eq: Q_nZ^k}
	\begin{split}
	&\p^{(e)}\left(\widetilde{B}_1^k(n)\geq \tau_n(t_1), \dots,\widetilde{B}_{k-1}^k(n)\geq \tau_n(t_{k-1}),\mathcal{O}=\mathfrak{t}\mid Z_n\geq k\right)\\
	&=\frac{\e^{(e)}\left[e^{-\theta_n Z_n}Z_n^{[k]}\right]}{(a_n^{(e)})^{k}\p^{(e)}\left(Z_n\geq k\right)}
	\mathbf{Q}^{(e,k,\theta_n)}_{n}\left[\frac{e^{\theta_n Z_n}}{(a_n^{(e)})^{-k}Z_n^{[k]}}\indi_{\{\widetilde{\psi}_1\geq \tau_n(t_1), \dots,\widetilde{\psi}_{k-1}\geq \tau_n(t_{k-1})\},} \mathcal{O}=\mathfrak{t}\right].
	\end{split}
	\end{equation}
	Now, similarly as before we complete the factors with $F_{k,n,0}$ defined in \eqref{Fknm} and use conditional probability, to obtain 
	\begin{equation*}
	\begin{split}
	\frac{\e^{(e)}\left[e^{-\theta_n Z_n}Z_n^{[k]}\right]}{(a_n^{(e)})^{k}\p^{(e)}\left(Z_n\geq k\right)}=
	\left(\frac{k!\ F_{k,n,0}}{(1+\theta)^{k+1}}\right)
	\left(\frac{\mu_n^{(e)}}{a_n^{(e)}\p^{(e)}\left(Z_n> 0\right)}\right)
	\left(\frac{1}{\p^{(e)}\left(Z_n\geq k\mid Z_n>0\right)}\right).
	\end{split}
	\end{equation*}
	According to equation \eqref{limit a}, the second fraction in the right hand side of the above identity goes to $1$, as $n\rightarrow \infty$. Moreover, Proposition \ref{limit starting Sn} implies that $Z_n$ conditioned on $\{Z_n>0\}$ goes to $\infty$, as $n\rightarrow\infty$, which implies  $\p^{(e)}\left(Z_n\geq k\mid Z_n>0\right)\rightarrow 1.$ Then, by Lemma \ref{lem: Yaglom}, we have
	\begin{equation}
	\label{first limit}
	\underset{n\rightarrow \infty}{\lim}
	\frac{\e^{(e)}\left[e^{-\theta_n Z_n}Z_n^{[k]}\right]}{(a_n^{(e)})^{k}\p^{(e)}\left(Z_n\geq k\right)}=
	\frac{k!}{(1+\theta)^{k+1}}.
	\end{equation}
	Therefore, the result holds as soon as we compute the asymptotic behaviour of the second term in the right-hand side of \eqref{eq: Q_nZ^k}. First, we observe   that $(Z_n-k)^{-k}\geq (Z_n^{[k]})^{-1}\geq (Z_n)^{-k}$ and
	$a_n^{(e)}\rightarrow \infty$, as $n\rightarrow \infty$. This implies that the aforementioned limit is the same as
	$$\underset{n\rightarrow \infty}{\lim}\mathbf{Q}^{(e,k,\theta_n)}_{n}\left[
	\frac{e^{\theta_n Z_n}}
	{\left((a_n^{(e)})^{-1}Z_n\right)^{k}}
	\indi_{\{\widetilde{\psi}_1\geq \tau_n(t_1), \dots,\widetilde{\psi}_{k-1}\geq \tau_n(t_{k-1})\}},\mathcal{O}=\mathfrak{t}\right].$$
	Now, if we  use the following identity
	\begin{equation}
	\label{eq: fraction and integral}
	\frac{e^{\theta x}}{x^{k}}=\int_0^{\infty}\frac{\lambda^{k-1}}{(k-1)!}e^{-(\lambda-\theta) x}\ud \lambda,
	\end{equation}
	we obtain
	\begin{equation*}
	\begin{split}
	&\mathbf{Q}^{(e,k,\theta_n)}_{n}
	\left[\frac{e^{\theta_n Z_n}}
	{\left((a_n^{(e)})^{-1}Z_n\right)^{k}}
	\indi_{\{\widetilde{\psi}_1\geq \tau_n(t_1), \dots,\widetilde{\psi}_{k-1}\geq \tau_n(t_{k-1})\}},\mathcal{O}=\mathfrak{t}\right]=\int_0^{\infty}\frac{\lambda^{k-1}}{(k-1)!} \xi_n(\lambda)\ud \lambda,
	\end{split}
	\end{equation*}
	where the function $\xi_n$ is defined by 
	\begin{equation*}
	\begin{split}
	\xi_n(\lambda):&=\mathbf{Q}^{(e,k,\theta_n)}_{n}\left[e^{-(\lambda-\theta) Z_n/a_n^e}\indi_{\{\widetilde{\psi}_1\geq \tau_n(t_1), \dots,\widetilde{\psi}_{k-1}\geq \tau_n(t_{k-1})\}},\mathcal{O}=\mathfrak{t}\right].
	\end{split}
	\end{equation*}
	Next, recall that $\mathcal{D}_n$ denotes the event that all spines splitting times are different.  
	By conditioning on  the values of $(\widetilde{\psi}_1,\dots ,\widetilde{\psi}_{k-1})$ and  $\mathcal{O}=\mathfrak{t}$, we have that 
	\begin{equation*}
	\begin{split}
	&\xi_n(\lambda)=\mathbf{Q}^{(e,k,\theta_n)}_{n}\left[
	e^{-(\lambda-\theta)Z_n/a_n^e}\indi_{\{\widetilde{\psi}_1\geq \tau_n(t_1), \dots,\widetilde{\psi}_{k-1}\geq \tau_n(t_{k-1})\}},\mathcal{O}=\mathfrak{t},\mathcal{D}_n^c\right]+\\
	&\mathbf{Q}^{(e,k,\theta_n)}_{n}\left[
	\mathbf{Q}^{(e,k,\theta_n)}_{n}\left[e^{-(\lambda-\theta)Z_n/a_n^e}\mid \widetilde{\psi}_1, ..., \widetilde{\psi}_{k-1},\mathcal{O},\mathcal{D}_n\right]\right.\\
&\left.	\hspace{7cm}\times\indi_{\widetilde{\psi}_1\geq \tau_n(t_1), \dots,\widetilde{\psi}_{k-1}\geq \tau_n(t_{k-1})},\mathcal{O}=\mathfrak{t},\mathcal{D}_n\right].
	\end{split}
	\end{equation*}
	Since $\underset{n\rightarrow \infty}{\lim}\mathbf{Q}^{(e,k,\theta_n)}_n\left(\mathcal{D}_n\right)=1$, we may use   Corollary \ref{cor: uniform}  and Lemma \ref{lem:limit Q conditioned psi}, to  see
	\[
	\begin{split}
	\xi_n(\lambda)\underset{n\rightarrow\infty}{\longrightarrow } &\p_{\mathcal{B}^k}(\mathfrak{t})\frac{(1+\theta)^{k+1}}{(1+\lambda)^{2}}\prod_{i=1}^{k-1}
	\int_{t_i}^{1}
	\frac{1}{(1+\lambda(1-r_i))^2}\ud r_i \\
	&=\p_{\mathcal{B}^k}(\mathfrak{t}) \frac{(1+\theta)^{k+1}}{(1+\lambda)^{2}\lambda^{k-1}}\prod_{i=1}^{k-1}\left(1-\frac{1}{1+\lambda (1-t_i)}\right).
	\end{split}
	\]
	Now, we define the function $\varsigma_n(\lambda):=\mathbf{Q}^{(e,k,\theta_n)}_{n}\left[e^{-(\lambda-\theta) Z_n/a_n^e}\right]$. Then, by the definition of $\xi_n$ and Lemma \ref{lem: Yaglom for Q}, we deduce
	$$\xi_n(\lambda)\leq \varsigma_n(\lambda) \underset{n\rightarrow\infty}{\longrightarrow }\frac{(1+\theta)^{k+1}}{(1+\lambda)^{k+1}}.$$
	Additionally, if we use again equation \eqref{eq: fraction and integral}, we get 
	$$\int_0^\infty \frac{\lambda^{k-1}}{(k-1)!}\varsigma_n(\lambda)\ud \lambda
	=\mathbf{Q}^{(e,k,\theta_n)}_{n}\left[\frac{e^{\theta_n Z_n}}
	{\left((a_n^{(e)})^{-1}Z_n\right)^{k}}\right]
	\sim \mathbf{Q}^{(e,k,\theta_n)}_{n}\left[\frac{e^{\theta_n Z_n}}{(a_n^{(e)})^{-k}Z_n^{[k]}}\right],$$
	as $n\rightarrow \infty$. Therefore, by using \eqref{eq: Q_nZ^k} with $0=t_1=\cdots=t_{k-1}$ and \eqref{first limit}, we get 
	\begin{equation}\label{density}
\int_0^\infty \frac{\lambda^{k-1}}{(k-1)!}\varsigma_n(\lambda)\ud \lambda\longrightarrow \frac{(1+\theta)^{k+1}}{k!}=\int_0^\infty \frac{\lambda^{k-1}}{(k-1)!}\frac{(1+\theta)^{k+1}}{(1+\lambda)^{k+1}}\ud \lambda.
\end{equation}
	Then, we apply the extended Dominated Convergence Theorem (see for instance \cite[Theorem 1.23]{Kallenberg},) with the sequence of functions $\{\lambda^{k-1}\xi_n(\lambda): n\geq 1\}$ and $\{\lambda^{k-1}\varsigma_n(\lambda): n\geq 1\}$ in order to obtain 
	\begin{equation*}
	\begin{split}
	&\underset{n\rightarrow\infty}{\lim}
	\mathbf{Q}^{(e,k,\theta_n)}_{n}
	\left[\frac{e^{\theta_n Z_n}}
	{\left((a_n^{(e)})^{-1}Z_n\right)^{k}}
	\indi_{\{\widetilde{\psi}_1\geq \tau_n(t_1), \dots,\widetilde{\psi}_{k-1}\geq \tau_n(t_{k-1})\}},\mathcal{O}=\mathfrak{t}\right]\\
	&\hspace{2cm}=\p_{\mathcal{B}^k}(\mathfrak{t})\frac{(1+\theta)^{k+1}}{(k-1)!}\int_0^{\infty}\frac{1}{(1+\lambda)^{2}}\prod_{i=1}^{k-1}\left(1-\frac{1}{1+\lambda (1-t_i)}\right)\ud \lambda.
	\end{split}
	\end{equation*}
	Thus, combining with equations \eqref{eq: Q_nZ^k} and \eqref{first limit},
	we have
	\begin{equation}\label{mixturerep}
	\begin{split}
	\underset{n\rightarrow\infty}{\lim}
	\p^{(e)}&\left(\widetilde{B}_1^k(n)\geq \tau_n(t_1), \dots,\widetilde{B}_{k-1}^k(n)\geq \tau_n(t_{k-1}),\mathcal{O}=\mathfrak{t}\mid Z_n\geq k\right)\\
	&=\p_{\mathcal{B}^k}(\mathfrak{t})\int_0^{\infty}\frac{k}{(1+\lambda)^{2}}\prod_{i=1}^{k-1}\left(1-\frac{1}{1+\lambda (1-t_i)}\right)\ud \lambda.\\
&=\p_{\mathcal{B}^k}(\mathfrak{t})\int_0^{\infty}\frac{k\lambda^{k-1}}{(1+\lambda)^{k+1}}
\prod_{i=1}^{k-1}\left(1-\frac{1}{1+\lambda (1-t_i)}\right)\frac{1+\lambda}{\lambda}\ud \lambda.\\
	\end{split}
	\end{equation}

Observe that the previous integral is 
the corresponding joint distribution of the density found in \cite[Theorem 4]{harris2020coalescent}. By computing the first integral explicitly, i.e. taking $e_0=1$ and $e_i=1-t_i$ for $i=1,\dots , k-1$ in \cite[Lemma 34]{HJRarxiv}, we deduce equation \eqref{eqtheorem}. 

Note that equation \eqref{mixturerep} implies, from the definition of  $\p_{\mathcal{B}^k}$, see \eqref{eq: lawbinarytrees}, that asymptotically all splittings are binary. Moreover, again from \eqref{mixturerep}, the tree topology $\{\mathcal{O}=\mathbf{t},\mathbf{t}\in\mathcal{B}^k\}$  and the splitting times are asymptotically independent.
Statement 2 follows directly from the definition of the probability $\p_{\mathcal{B}^k}$. 

We observe  that \eqref{mixturerep} implies that  
we can view the coalescent times as a mixture of independent random variables. More precisely, with density function $k\lambda^{k-1}(1+\lambda)^{-k-1}$, we select a parameter $\lambda$. (Note that it is a distribution by \eqref{density}). Given  $\lambda$, we take $k-1$ independent random variables with tail distribution given by $(1+\lambda)(1-t)(1+\lambda(1-t))^{-1}$. Finally, for Statement 1, we use the mixture of independent random variables.  If at some time, there are 
 $i$ blocks of sizes $b_1,\dots, b_{i}$, then, conditionally on $\lambda$, block $j$ has associated $b_j-1$ of the previous random variables. Therefore, the first group to split will be group $j$ with probability proportional to $b_j-1$, that is, with probability $(b_j-1)/(k-i)$. 
 
\end{proof}	

\section{Appendix}
The following lemma is an extension of Lemma 8 in Kersting \cite{kersting2020unifying},  which provides the  case $m=0$.  
\begin{lemma}\label{lemma:extensionkersting8}
Assume that conditions  \eqref{eq_cond_kersting} and \eqref{def: critical} are  satisfied. Then, 
\begin{equation*}
\underset{0\leq s\leq 1}{\sup}\left|\underset{k=m+1}{\overset{n}{\sum}}\frac{\varphi_k(f_{k,n}(s))}{\mu_{k-1}}-\frac{\rho_n-\rho_m}{2}\right|= o\left(\frac{\rho_n}{2}\right),
\end{equation*}
uniformly for all $m<n$, as $n\rightarrow\infty$.
\end{lemma}

\begin{proof}
For any $m<n$, we  apply the triangle inequality to obtain
\begin{equation*}
\begin{split}
\left|\underset{k=m+1}{\overset{n}{\sum}}\frac{\varphi_k(f_{k,n}(s))}{\mu_{k-1}}-\frac{\rho_n-\rho_m}{2}\right|
\leq 
\left|\underset{k=1}{\overset{n}{\sum}}\frac{\varphi_k(f_{k,n}(s))}{\mu_{k-1}}-\frac{\rho_n}{2}\right|
+\left|\underset{k=1}{\overset{m}{\sum}}\frac{\varphi_k(f_{k,n}(s))}{\mu_{k-1}}-\frac{\rho_m}{2}\right|.
\end{split}
\end{equation*}
Observe that for any $s\in [0,1]$ and $k\leq m$, we can decompose $f_{k,n}(s)=f_{k,m}(f_{m,n}(s))$ where $f_{m,n}(s)\in[0,1]$. Then, by applying supremum on both sides of the previous inequality, we deduce
\begin{equation*}
\begin{split}
&\underset{0\leq s\leq 1}{\sup}\left|\underset{k=m+1}{\overset{n}{\sum}}\frac{\varphi_k(f_{k,n}(s))}{\mu_{k-1}}-\frac{\rho_n-\rho_m}{2}\right|\\
&\leq 
\underset{0\leq s\leq 1}{\sup}
\left|\underset{k=1}{\overset{n}{\sum}}\frac{\varphi_k(f_{k,n}(s))}{\mu_{k-1}}-\frac{\rho_n}{2}\right|
+\underset{0\leq s\leq 1}{\sup}\left|\underset{k=1}{\overset{m}{\sum}}\frac{\varphi_k(f_{k,m}(s))}{\mu_{k-1}}-\frac{\rho_m}{2}\right|.
\end{split}
\end{equation*}
Since the sequence $\{\rho_n,n\geq 0\}$ is increasing and $\rho_n\rightarrow \infty$ as $n\rightarrow \infty$, the result holds by \cite[Lemma 8]{kersting2020unifying}.
\end{proof} 

\bibliographystyle{abbrv}

\end{document}